\documentclass[12pt,a4paper]{report}
\usepackage{fancyhdr}
\usepackage{etex}
\usepackage[mathscr]{eucal}
\usepackage{amsfonts,amsmath,mathtools,amsthm,amscd,amssymb,color}
\allowdisplaybreaks

\usepackage{latexsym,gensymb}
\usepackage{graphicx}
\usepackage{graphics}
\graphicspath{{images}}
\usepackage[font=small,skip=7pt]{caption}
\usepackage{xcolor}
\usepackage{dcpic,pictex}
\usepackage{multicol}
\usepackage{tikz,tikz-3dplot}
\usepackage{pgf,pgfplots}
\usetikzlibrary{shapes.geometric,decorations.markings,arrows,shadings,patterns}
\usepackage{textcomp}
\usepackage{epsfig}
\usepackage[T1]{fontenc}
\usepackage{makeidx}
\makeindex
\usepackage{pdflscape}
\usepackage{filecontents}
\usepackage{blkarray}
\usepackage{titlesec}

\usepackage[nottoc]{tocbibind}
\usepackage{bm}
\usepackage{listings}
\usepackage{float}
\usepackage{cite}

\numberwithin{equation}{chapter}

\tdplotsetmaincoords{70}{165}

\titleformat{\chapter} [display]
{\normalfont\large\bfseries}{\filcenter\chaptertitlename\ \thechapter}
{3pt}{\LARGE\filcenter}
\titleformat{\section} 
{}{ \filcenter\bfseries
	\thesection}{0pt} { 
	\large\bfseries\filcenter} 

\newtheorem{theorem}{Theorem}[section]
\newtheorem{corollary}[theorem]{Corollary} 
\newtheorem{lemma}[theorem]{Lemma} 

\newtheorem{remark}[theorem]{Remark} 

\newtheorem{proposition}[theorem]{Proposition}

\newtheorem{definition}[theorem]{}
\newtheorem{Definition}[theorem]{Definition}

\newtheorem{conjecture}[theorem]{Conjecture}

\newsavebox{\smallblockbox}

\usepackage[top=1.55in, bottom=1.55in, left=1.16in, right=1.17in]{geometry}
\begin{document}
\thispagestyle{empty}
\enlargethispage{1cm}
\thispagestyle{empty}
\begin{center}
{\Large {\bf\scshape Title}}\\
\rule{\textwidth}{1.0pt}\vspace*{-\baselineskip}\vspace*{2pt} 
\rule{\textwidth}{0.8pt}\\[\baselineskip]
{\bf {\LARGE \scshape \textcolor{blue}{ ON DISTRIBUTION OF LAPLACIAN EIGENVALUES OF GRAPHS} }}
\rule{\textwidth}{0.8pt}\vspace*{-\baselineskip}\vspace{3.2pt} 
\rule{\textwidth}{1.0pt}\\ 
\vspace{1\baselineskip}
\begin{figure}[!ht]
\centering
\includegraphics[width=140pt]{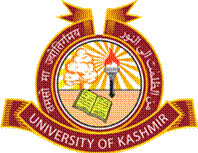}
\end{figure}
{\large A thesis submitted in partial fulfilment of the requirements for the degree of}\\
\vspace*{6mm}
{\Large{\bf\scshape Doctor of Philosophy}}\\
\vspace*{3mm}
{\large in}\\
\vspace*{3mm}
{\Large {\bf\scshape Mathematics}}\\
\vspace*{3mm}
{\large by}\\
\vspace*{3mm}
{\Large {\bf\scshape Bilal Ahmad Rather}}\\
\vspace{3mm}
{\large under the joint supervision of\\
\vspace*{3mm}
{\Large\bf\scshape  Prof. S. Pirzada \quad and\quad Prof. T. A. Chishti}}
\end{center} 

\begin{center}
\textbf{\Large\scshape Department of Mathematics}\\ \vspace{-.1\baselineskip}
\textbf{\scshape University of Kashmir}\\ \vspace{-.1\baselineskip}
\textbf{\small\scshape Srinagar,  Jammu and Kashmir 190006}\\ { NAAC Accredited Grade A$ ^{+} $ University}\\\vspace{-.1\baselineskip}
{{\bf October-2020}}
\end{center}
\newpage
\pagenumbering{roman}

\thispagestyle{plain}
\addcontentsline{toc}{chapter}{Declaration}
\begin{center}
  \textbf{\LARGE \underline{DECLARATION}}
\end{center}
\vspace*{1cm}
\noindent
I, {\bf\scshape  Bilal Ahmad Rather}, hereby declare that this thesis entitled, {\bf\scshape On Distribution of  Laplacian Eigenvalues of Graphs} and the work presented in it are my own. I confirm that:
\begin{itemize} 
\item This work was done wholly while in candidature for a research degree at the University of Kashmir, Srinagar.
\item Where any part of this thesis has not previously been submitted for a degree or any other qualification at this University or any other institution, this has been clearly stated.
\item Where I have consulted the published work of others, this is always clearly attributed.
\item I have also fulfilled the requirements of the UGC regulations for carrying out research work in Ph.D.\\
\end{itemize}
\vspace{1.5cm}

\noindent Signed:\\
\rule[0.5em]{15em}{0.5pt} 
 
\noindent Date:\\
\rule[0.5em]{15em}{0.5pt}\\~\\

\noindent \hfill{\bf Bilal Ahmad Rather}.\\

\newpage
\thispagestyle{plain}
\addcontentsline{toc}{chapter}{Certificate}
\begin{center}
\includegraphics[width=4.8cm,height=3.7cm]{KULogo1}\\
\vspace*{5mm}
\end{center}
\begin{center}
\textbf{\Large\scshape Department of Mathematics}\\ \vspace{-.1\baselineskip}
\textbf{\scshape University of Kashmir}\\ \vspace{-.1\baselineskip}
\textbf{\small\scshape Srinagar, jammu and Kashmir 190006}\\ { NAAC Accredited Grade A$ ^{+} $ University}\\ \vspace{-.1\baselineskip}
{{\bf October-2020}}
\end{center}
\vspace*{5mm}
{\qquad \qquad \qquad \qquad\quad~ \qquad \color {red} {\bf \Huge \underline {CERTIFICATE}}}
\bigskip
\bigskip
\bigskip

Certified that the thesis entitled {\color {blue} {\bf\scshape"on distribution of Laplacian eigenvalues of graphs"}} being submitted by Bilal Ahmad Rather, in partial fulfilment of the requirements for the award of \textbf{Doctor of Philosophy in Mathematics, School of Physical and Mathematical Sciences, University of Kashmir, Srinagar}, is his own work carried out by him under our supervision and guidance. The content of this thesis, in full or in parts, has not been submitted to any Institute or University for the award of any degree or diploma.\\~\\~\\~\\

\noindent Professor S. Pirzada.\qquad \qquad\qquad\qquad \qquad\qquad\qquad\qquad Professor T. A. Chishti\\
\indent Supervisor \qquad\qquad\qquad \qquad\qquad\qquad\qquad\qquad\qquad\qquad Co-supervisor\\

\begin{center}
 Professor B. A. Zargar\\
Head of the Department
\end{center}
\newpage

\newpage
\thispagestyle{plain}
\begin{center}
  \textbf{\Huge Acknowledgement}
\end{center}
\vspace*{1cm}

Words do have shortcoming of not being able to express the immense of my gratitude to the Almighty Allah for guiding me all the way in my life and for helping me accomplish this task. I would like to warmly acknowledge and express my deep sense of
gratitude and indebtedness to my guides Prof. S. Pirzada, Department of Mathematics (Dean of Physical and Mathematical Sciences) and Prof. T. A. Chishti, Director of Distance Education, University of Kashmir, for their keen guidance, constant encouragement and prudent suggestions during the course of my study and preparation of the final manuscript of this work.\\
\indent I am grateful to Prof. B. A. Zargar, Head Department of Mathematice, University of Kashmir for providing the required infrastructure in the department and facilitating us in the official procedures.\\
\indent I also take this opportunity to acknowledge the help and support of all the teaching staff and non teaching members of the department.\\
\indent It gives me great pleasure to convey my sincere thanks to all the research scholars in my department. I am highly thankful to Dr Mushtaq Ahmad Bhat, Department of Higher Education,  Kashmir for his research guidance, motivation and sincere support. I also acknowledge Prof. Vilmar Trevisan, Institute de Matem\'atica, UFRGS, Brazil and  Dr. Hilal Ahmad Ganie, Department of School Education Kashmir,  for their advice and crucial support regarding the joint research work carried during this period. I would also like to thank all my friends, associates and well-wishers for their motivation and support.\\
\indent I also wish to acknowledge the support rendered by my family in the form of encouragement, personal support and for their prayers and concern along my way.\\~\\

\vskip .5cm
\vskip .5cm
\noindent {October 2020} \hfill{\bf Bilal Ahmad Rather}
\bigskip
\bigskip
\bigskip
\bigskip

\chapter*{List of Publications}
\begin{enumerate}

\item S. Pirzada, Bilal A. Rather, M. Aijaz and T. A. Chishti, On distance signless Laplacian spectrum of graphs and spectrum of zero-divisor graphs of $ \mathbb{Z}_{n} $,\emph{ Linear and Multilinear Algebra}, (2020) DOI:10.1080/03081087.2020.1838425.\\
\textbf{SCI, Scopus.}

\item S. Pirzada, H. A. Ganie, Bilal A. Rather and R. U. Shaban, On generalized distance energy of graphs, \emph{ Linear Algebra and its Applications} \textbf{603} (2020) 1-19. \\
\textbf{SCI, Scopus.}

\item	H. A. Ganie, S. Pirzada, Bilal A. Rather and R. U. Shaban, On Laplacian eigenvalues of graphs and Brouwer's conjecture, \textit{ Journal of the Ramanujan Mathematical Society} \textbf{36(1)} (2021) 1--9.\\
\textbf{SCI, Scopus.} 

\item H. A. Ganie, S. Pirzada, Bilal A. Rather and V. Trevisan, Further development on Brouwer's conjecture for the sum of Laplacian eigenvalues of graphs, \emph{ Linear Algebra and its Applications}  \textbf{558} (2020) 1-18. \\
\textbf{SCI, Scopus.}

\item H. A. Ganie, Bilal A. Rather and S. Pirzada, Laplacian energy of trees of diameter four and beyond, communicated (2020).

\item Bilal A. Rather, S. Pirzada, T. A. Chishti and M. A. Alghamdi, On normalized Laplacian eigenvalues of power graphs of finite cyclic group, communicated (2020).

\end{enumerate}
\chapter*{Abstract}
\addcontentsline{toc}{chapter}{Abstract}
\qquad The work in this thesis concerns the investigation of eigenvalues of the Laplacian matrix, normalized Laplacian matrix, signless Laplacian matrix and distance signless Laplacian matrix of graphs.  The thesis consists of five chapters.

 In Chapter 1, we present a brief introduction of spectral graph theory and we list some definitions required throughout the thesis. 

Chapter $2$ deals with the sum of $ k $ largest Laplacian eigenvalues $ S_{k}(G) $ of graph $ G $ and Brouwer's conjecture. We obtain the upper bounds for  $ S_{k}(G) $  for some classes of graphs and use them to verify Brouwer's conjecture for these classes of graphs. Also, we  prove Brouwer's conjecture for more general classes of graphs.  

In Chapter $3$, we investigate the Laplacian eigenvalues of graphs and the Laplacian energy conjecture for trees. We prove the Laplacian energy conjecture completely for trees of diameter $ 4 $. Further, we prove this conjecture for all trees having at most $ \frac{9n}{25}-2 $ non-pendent vertices. Also, we obtain the sufficient conditions for the truth of conjecture for trees of order $ n $. 

In Chapter $4$, we determine the normalized Laplacian spectrum of the joined union of regular graphs and obtain the spectrum of some well known graphs. As consequences of joined union, we obtain the normalized Laplacian spectrum of power graphs associated to finite cyclic groups.

In Chapter $5$, we find the distance signless Laplacian spectrum of regular graphs and zero-divisor graphs associated to finite commutative ring. Also, we find the bounds for spectral radius of generalized distance matrix. Further, we obtain the generalized distance energy for bipartite graphs and trees. We prove that the complete bipartite graph has minimum generalized distance energy among all connected bipartite graphs. Besides, for $ \alpha\in \big(0, \frac{2n}{3n-2}\big) $, we show that the star graph has minimum generalized distance energy among all trees.

\chapter*{NOTATIONS}

Here we give some symbols and notation which will be used throughout this work. Many other definitions and results will be found wherever required.
\begin{center}
\begin{tabular}{l l}
\textbf{Notation} & \textbf{Definition}\\
$G$ & Graph\\
$ \overline{G} $ & Complement of $ G. $\\
V(G) & Vertex set of $ G $\\
$E(G)$ & Edge set of $ G $\\
$N(v)$ & Neighborhood of vertex $ v $\\
$d_{v}$ & Degree of vertex $ v $\\
$ \overline{d} $ & Average degree\\
$d_v^{*}$ & Conjugate degree of vertex $ v $\\
$ \omega $ & Clique number\\
$ d(u,v) $ & Distance between $ u $ and $ v $\\
$ Tr_{G}(v) $ & Transmission of vertex $ v $\\
$ W(G) $ & Transmission number or Wiener index\\
$ u\sim v $ & $ u $ is adjacent to $ v $\\
$ g $ & Girth\\
$ T $ & Tree\\
 $ SNS $ & Star-NonStar tree\\
$K_{n}$ & Complete graph \\
$K_{a,b}$ & Complete bipartite graph \\
$S_{n} $ & Star graph \\
$CS_{\omega,n-\omega}$ & Complete split graph \\
$P_{n}$ & Path \\

\end{tabular}
\end{center}

\begin{center}
\begin{tabular}{l l}
\textbf{Notation} & \textbf{Definition}\\
$C_{n}$ & Cycle \\
$W_{n}$ & Wheel \\
$F_{n}$ & Friendship graph \\
$W_{a,b}$ & Generalized wheel \\
$C_{a,b}$ & Cone graph \\
$F_{p,n-p}$ & Firefly graph \\
$ K_{n_{1},n_{2},\dots,n_{p}} $ & Complete $ p $-partite graph\\
$\Upsilon_{n}$ & Proper divisor graph \\
$G[G_{1},\dots,G_{n}]$ & Joined union of graphs\\
$ G_{1}\triangledown G_{2} $ & Join of two graphs\\
$ \mathcal{G} $ & Group\\
$ \mathbb{Z}_{n} $ & Integer modulo $ n $\\
$ \mathcal{P}(\mathcal{G}) $ & Power graph of group $ \mathcal{G} $\\
$ R $ & Ring\\
$ Z(R) $ & Zero-divisors of $ R $\\
$ Z^{*}(R) $ & Non-zero zero-divisors of $ R $\\
$ \Gamma(R) $ & Zero-divisor graph of ring $ R $\\
$ \Upsilon_{n} $ & Proper divisor graph\\
$ \phi $ & Euler's totient function\\
$ |~| $ & Cardinality\\
$ \mathbb{Q}(\lambda) $ & Polynomial ring\\
$ d|n $ & $ d $ Divides $ n $\\
$ D(n) $ & Proper divisor set of $ n $\\
$ A(G) $ & Adjacency matrix\\
$ L(G) $ & Laplacian matrix\\
$ \mathcal{L}(G) $ & Normalized Laplacian matrix\\
$ Q(G) $ & Signless Laplacian matrix\\
$\mathcal{D}(G) $ & Distance matrix\\
$\mathcal{D}^{L}(G) $ & Distance Laplacian matrix\\
$\mathcal{D}^{Q}(G) $ & Distance signless Laplacian matrix\\
$\mathcal{D}_{\alpha}(G) $ & Generalized distance matrix\\

\end{tabular}
\end{center}
\begin{center}
\begin{tabular}{l l}
\textbf{Notation} & \textbf{Definition}\\
$ \mathcal{Q} $ & Quotient matrix\\
$ \textbf{J} $ & Vector whose each entry equals $ 1 $\\
$ e_{n_{i}} $ & $ (\underbrace{1,1,\dots,1}_{n_{i}})^{T} $\\
$ \textbf{O}_{n_{i}\times n_{j}} $ & Zero matrix of order $ n_{i}\times n_{j} $\\
$ \lambda $ & Eigenvalue of a matrix or eigenvalues of $ A(G) $\\
$ \mu $ & Eigenvalue of $ L(G) $\\
$ \rho $ & Eigenvalue of $ \mathcal{L}(G) $\\
$ q $ & Eigenvalue of $ Q(G) $\\
$\rho^{\mathcal{D}} $ & Distance  eigenvalue\\
$\rho^{Q} $ & Distance signless Laplacian eigenvalue\\
$ \partial $ &  Spectral radius  of $\mathcal{D}_{\alpha}(G) $\\
$ \partial_{i} $ &  Eigenvalues of $\mathcal{D}_{\alpha}(G) $\\
$ s_{i} $ &  Singular value of a matrix\\
$ Spec_{Q}(G) $ & Signless Laplacian spectrum of $ G $\\
$ S_{k}(G) $ & Sum of $ k $ largest $ \lambda_{i} $'s\\
$ S_{k}^{+}(G) $ & Sum of $ k $ largest $ q_{i} $'s\\
$ \sigma $ & Positive integer satisfying $ \lambda_{\sigma}\geq \overline{d} $\\
$ \sigma^{'} $ & Positive integer satisfying $ q_{\sigma^{'}}\geq \overline{d} $\\

$ \mathcal{E}(G) $ & Energy of $ G $\\
$ LE(G) $ & Laplacian energy of $ G $\\
$ QE(G) $ & Signless Laplacian energy of $ G $\\
$ E_{\alpha}^{\mathcal{D}}(G) $ & Generalized distance energy of $ G $\\
$ \| ~~ \|_{p} $ & Schatten $ p $-norm\\
$ \| ~~ \|_{F} $ & Frobenius norm\\
$ \| ~~ \|_{k} $ & Ky Fan k-norm\\

\end{tabular}
\end{center}

\newpage
\listoffigures
\tableofcontents
\markboth{}{CONTENTS}
\newpage
\newpage
\pagenumbering{arabic}


\chapter{Introduction}
\section{Background}
Spectral graph theory (Algebraic graph theory) is the study of spectral properties of matrices associated to graphs. The spectral properties include the study of characteristic polynomial, eigenvalues and eigenvectors of matrices associated to graphs. This also includes the graphs associated to algebraic structures like groups, rings and vector spaces. The major source of research in spectral graph theory has been the study of relationship between the structural and spectral properties of graphs. Another source has research in mathematical chemistry (theoretical/quantum chemistry). One of the major problems in spectral graph theory lies in finding the spectrum of matrices associated to graphs completely or in terms of spectrum of simpler matrices associated with the structure of the graph.  Another problem which is worth to mention is to characterise the extremal graphs among all the graphs or among a special class of graphs with respect to a given graph, like spectral radius, the second largest eigenvalue, the smallest eigenvalue, the second smallest eigenvalue, the graph energy and multiplicities of the eigenvalues that can be associated with the graph matrix. The main aim is to discuss the principal properties and structure of a graph from its eigenvalues. It has been observed that the eigenvalues of graphs are closely related to all graph parameters, linking one property to another. \\

Spectral graph theory has a wide range of applications to other areas of mathematical science and to other areas of sciences which include Computer Science, Physics, Chemistry, Biology, Statistics, Engineering etc. The study of graph eigenvalues has rich connections with many other areas of mathematics. An important development is the interaction between spectral graph theory and differential geometry. There is an interesting connection between spectral Riemannian geometry and spectral graph theory. Graph operations help in partitioning of the embedding space, maximising inter-cluster affinity and minimising inter-cluster proximity. Spectral graph  theory plays a major role in deforming the embedding spaces in geometry. Graph spectra helps us in making conclusions that we cannot recognize the shapes of solids by their sounds. Algebraic spectral methods are also useful in studying the groups and the rings in a new light. This new developing field investigates the spectrum of graphs associated with the algebraic structures like groups and rings. The main motive to study these algebraic structures graphically using spectral analysis is to explore several properties of interest.\\

In 2010 monograph `An Introduction to the Theory of Graph Spectra' by Cvetkovi\'{c}, Rowlinson and Simi\'c \cite{cds1} summarised all the results to the date and included the results of 1980 text. For some more books on spectral graph theory, we refer to \cite{BH1,bp1,nb,chung,cds,cds1,gr,jason}.\\

The adjacency matrix $ A(G) $ of graph $ G $ is a $ (0,1) $ square matrix of order $ n $ whose $ (i,j) $-entry is $ 1 $ if $ v_{i} $ is adjacent to $ v_{j} $ and $ 0 $, otherwise. This matrix is real symmetric, so its eigenvalues are real. The set of all eigenvalues including multiplicities is known as adjacency spectrum (or simply spectrum) of $ G $ and can be ordered as $ \lambda_{1}\geq \lambda_{2}\geq \dots\geq \lambda_{n} $. The largest eigenvalue $ \lambda_{1} $ is known as the spectral radius of $ A(G) $. The energy of a graph $ G $ is  denoted by $ \mathcal{E}(G) $, is defined as the sum of the absolute values of eigenvalues of $ A(G) $, that is \[ \mathcal{E}(G)=\sum\limits_{i=1}^{n}|\lambda_{i}|. \]
The concept of energy of a graph have its origin in the pioneering work of H\"{u}ckel \cite{ggt,h} who made simplification of Schrodinger's wave equations \cite{cds}. Chemists are interested in finding the wave functions and energy levels of a given molecule.  The spectra of graphs can be used to calculate the energy levels of a conjugated hydrocarbon as calculated with H\"{u}ckel Molecular Orbital (HMO) method. The details of H\"{u}ckel theory and how it is related to spectral graph theory can be found in \cite{cds}. An important quantum-chemical characteristic of a conjugated molecule is its total $\pi$-electron energy. Within the H\"{u}ckel Molecular Orbital (HMO) theory this quantity can be reduced to the energy of adjacency matrix.  Several problems in Chemistry were modelled by graphs and were solved using spectral graph theory \cite{glaa}. Almost all the results which were published on graph energy up to 2012 are summarised in the book \cite{graphenergy} by Xueliang Li, Yongtang Shi and Ivan Gutman. But this energy concept is actually the \emph{Ky Fan k-norm}, which in case of symmetric matrix is the sum of absolute values of the eigenvalues of matrix.\\

If $ D(G)=diag(d_{1},d_{2},\dots,d_{n}) $ is the diagonal matrix of vertex degrees $ d_{i} $, then the positive semi-definite matrix $ L(G)=D(G)-A(G) $ is known as the Laplacian matrix of $ G $. Since row  and column sums of $ L(G) $ are zero, so $ 0 $ must be eigenvalue of $ L(G) $. The set of all eigenvalues (with multiplicities) of $ L(G) $ is called the Laplacian spectrum of $ G $ and is denoted by $ \mu_{1}\geq\mu_{2}\geq \dots\geq \mu_{n-1}\geq\mu_{n}=0. $ The second smallest Laplacian eigenvalue $ \mu_{n-1} $ is positive if and only if $ G $ is connected and is known as algebraic connectivity of $ L(G) $. \emph{Schur's inequality} states that eigenvalues of $ L(G) $ are majorized by the degree sequence of $ G $.\emph{ Grone-Merris-Bai theorem} says that the sum of the first $ k $ conjugate degrees of $ G $ are bounded below by the sum of $ k $ largest Laplacian eigenvalues. \emph{Brouwer's Conjecture} states that the sum of the first $ k $ largest Laplacian eigenvalues of $ G $ are bounded above by $ \binom{k+1}{2}+m $, where $ m $ is size of $ G $. It is well known that Brouwer's bound is sharper than Grone-Merris-Bai bound for at least one $ k. $ If $ \overline{d} $ is the average of the Laplacian eigenvalues of $ G $, the Laplacian energy of $ G $, denoted by $ LE(G) $, is defined as
\[ LE(G)=\sum\limits_{i=1}^{n-1}\left |\mu_{i}-\overline{d}  \right |. \]
In the literature, several research papers can be found on  Laplacian energy and on characterizing the graphs attaining the maximum and minimum Laplacian energy. But the problem is hard and is not fully solved. One such problem is Laplacian energy conjecture for trees, which states that the Laplacian energy of a tree is bounded below by Laplacian energy of a path.\\

The normalized Laplacian matrix of a graph $ G $ is defined by $ \mathcal{L}(G)=D(G)^{-\frac{1}{2}}L(G)D(G)^{-\frac{1}{2}} $, where $ D(G)^{-\frac{1}{2}} $ is the diagonal matrix whose $ i $-the diagonal entry is $ \frac{1}{\sqrt{d_{i}}} $. This matrix is real symmetric and positive semi-definite, so its eigenvalues are arranged as $ \rho_{1}\geq \rho_{2}\geq \dots\geq \rho_{n} $ and are known as normalized Laplacian eigenvalues of $ G. $ The smallest normalized Laplacian eigenvalue $ \rho_{n} $ is always $ 0 $.\\

The signless Laplacian matrix of a graph $ G $ is defined as $ Q(G)=D(G)+A(G) $. This matrix is real symmetric and positive semi-definte, so eigenvalues of $ Q(G) $ are real, denoted by $ q_{i} $ and arranged as $ q_{1}\geq q_{2}\geq\dots\geq q_{n} $. In case $ G $ is bipartite, the Laplacian eigenvalues coincide with signless Laplacian eigenvalues. The set of all eigenvalues (including multiplicities) of $ Q(G) $ are known as signless Laplacian spectrum of $ G $. The largest signless Laplacian eigenvalues $ q_{1} $ is known as $ Q $-index or the signless Laplacian spectral radius of $ G. $ The signless Laplacian energy of $ G $ is denoted by $ QE(G) $ and is defined by 
\[ QE(G)=\sum\limits_{i=1}^{n}|q_{i}-\overline{d}|. \]
The invariants $ LE(G) $ and $ QE(G) $ are the trace norms of the symmetric matrices $ L(G)-\frac{2m}{n}I_{n} $ and $ Q(G)-\frac{2m}{n}I_{n} $, where $ I_{n} $ is the identity matrix.\\

\indent In a graph $G$, the \textit{distance} between any two vertices $u,v\in V(G),$ denoted by $d(u,v)$, is defined as the length of a shortest path between $u$ and $v$. The \textit{distance matrix} of $G$, denoted by $\mathcal{D}(G)$, is defined as $\mathcal{D}(G)=(d_{uv})$, where $u,v\in V(G)$. The matrix $\mathcal{D}(G)$ is real symmetric and positive definite for $ n\geq 3 $, so its eigenvalues can be arranged as $\rho^{\mathcal{D}}_{1}(G)\geq \dots\geq\rho^{\mathcal{D}}_{n-1}(G)\geq \rho^{\mathcal{D}}_{n}(G)$, where $ \rho^{\mathcal{D}}_1(G) $ is called the distance spectral radius of $ G. $ The distance energy of $ G $ is defined as the absolute sum of distance eigenvalues of $ \mathcal{D}(G) $, that is
\[ E^{\mathcal{D}}(G)=\sum\limits_{i=1}^{n}|\rho^{\mathcal{D}}_{i}(G)|. \]

\indent The \textit{transmission degree} $Tr_{G}(v)$ of a vertex $v$ is defined to be the sum of the distances from $v$ to all other vertices in $G$, that is, $Tr_{G}(v)=\sum\limits_{u\in V(G)}d(u,v)$. If $Tr_G(v_i)$ (or simply $Tr_{i}$) is the transmission degree of the vertex $v_i\in V(G)$, the sequence $\{Tr_{1},Tr_{2},\ldots,Tr_{n}\}$ is called the \textit{transmission degree sequence} of the graph $G$. Let $Tr(G)=diag (Tr_1,Tr_2,\ldots,Tr_n) $ be the diagonal matrix of vertex transmissions of $G$. The matrix $ \mathcal{D}^L(G)=Tr(G)-\mathcal{D}(G) $ is called the \textit{distance Laplacian matrix} of $G$, while the matrix  $\mathcal{D}^{Q}(G)=Tr(G)+\mathcal{D}(G)$ is called the \textit{distance signless Laplacian matrix of $G$.} The matrix $\mathcal{D}^{Q}(G)$ is real symmetric and positive definite for $ n\geq 3 $, so its eigenvalues can be arranged as $\rho^{Q}_{1}(G)\geq \dots\geq\rho^{Q}_{n-1}(G)\geq \rho^{Q}_{n}(G)$, where $ \rho^{Q}_1(G) $ is called the distance signless Laplacian spectral radius of $ G. $ Since the matrix $\mathcal{D}^{Q}(G)$ is non-negative and irreducible, so by \emph{Perron-Frobenius Theorem}, $ \rho^{Q}_1(G) $ is positive and simple zero of the characteristic polynomial of $\mathcal{D}^{Q}(G)$. There corresponds a unique eigenvector with positive entries having unit length, which is called the distance signless Laplacian \emph{Perron vector} of graph $ G. $ The distance signless Laplacian energy of $ G $ is defined as the absolute sum of  eigenvalues of $ \mathcal{D}^{Q}(G)-\frac{2W(G)}{n}I_{n} $, that is
\[ E^{Q}(G)=\sum\limits_{i=1}^{n}\left |\rho^{Q}_{i}(G)-\frac{2W(G)}{n}\right |, \]
where $ I_{n} $ is the identity matrix.

The \textit{generalized distance matrix} $\mathcal{D}_{\alpha}(G)$ of $ G $ introduced by Cui et al. \cite{CHT} as a convex combinations of $Tr(G)$ and $\mathcal{D}(G)$, and is defined as $\mathcal{D}_{\alpha}(G)=\alpha Tr(G)+(1-\alpha)\mathcal{D}(G)$,  for $0\leq \alpha\leq 1$. Clearly, $\mathcal{D}_{0}(G)=\mathcal{D}(G)$ (distance matrix), $2\mathcal{D}_{\frac{1}{2}}(G)=\mathcal{D}^{Q}(G)$ (distance signless Laplacian matrix), and  $\mathcal{D}_{\alpha}(G)-\mathcal{D}_{\beta}(G)=(\alpha-\beta)\mathcal{D}^{L}(G)$ (distance laplacian matrix).
The matrix $ \mathcal{D}_{\alpha}(G) $ reduces to merging the distance spectral, distance Laplacian spectral and distance signless Laplacian spectral theories. Since the matrix $ \mathcal{D}_{\alpha}(G)$ is real symmetric, all its eigenvalues are real. Therefore, we can arrange them  as $ \partial_{1}\geq \partial_{2}\geq \dots \geq \partial_{n}$. The largest eigenvalue  $ \partial_{1} $ of the matrix $\mathcal{D}_{\alpha}(G)$ is called the \textit{generalized distance spectral radius} of $ G$. As $ \mathcal{D}_{\alpha}(G) $ is non-negative and irreducible, by the Perron-Frobenius theorem, $ \partial(G)$ is the unique eigenvalue and there is a unique positive unit eigenvector $X$ corresponding to $ \partial_{1},$ which is called the \textit{generalized distance Perron vector} of $G.$ 

The generalized distance energy of $G$ is defined as the mean deviation of the  values of the generalized distance eigenvalues of $G$, that is,
\begin{equation*}
E^{\mathcal{D}_{\alpha}}(G)=\sum_{i=1}^{n}\left|\partial_i-\frac{2\alpha W(G)}{n}\right|=\displaystyle\sum_{i=1}^{n}|\Theta_{i}|.
\end{equation*}
Clearly, $ \sum\limits_{i=1}^{n}
|\Theta_{i}| $ is the trace norm of the matrix $ \mathcal{D}_{\alpha}(G)-\frac{2W(G)}{n}I_{n} $, where $ I_{n} $ is identity matrix. Further, $\sum_{i=1}^{n}\Theta_{i}=0$, $E^{\mathcal{D}_{0}}(G)=E^{\mathcal{D}}(G)$ and $2E^{\mathcal{D}_{\frac{1}{2}}}(G)=E^{Q}(G).$ This shows that the concept of generalized distance energy of $G$ merges the theories of distance energy and the distance signless Laplacian energy of $G$. Therefore, it will be interesting to study the quantity $E^{\mathcal{D}_{\alpha}}(G)$ and explore some properties like the bounds, the dependence on the structure of graph $G$ and the dependence on the parameter $\alpha$.

\section{Basic Definitions}\label{sec1}
In this section we give some basic definitions of graphs. Our graph notations are standard and can be taken from \cite{ping,p1}, for Abstract Algebra we follow \cite{roman,Atiyah} and for Matrix Theory and Linear Algebra we follow \cite{hj,fz, golub, jhonson, roman1}.
\begin{Definition}\textbf{ Graph}. A graph $G$ is a pair $(V, E)$, where ${V}$ is a nonvoid set of objects called vertices and $E$ is a subset of ${V}^{(2)}$, (the set of distinct unordered pairs of different elements of ${V}$). The elements of $E$ are called edges of $G$.
\end{Definition}

\begin{Definition}\textbf{Multigraph}. A multigraph $G$ is a pair $({V}, E)$, where ${V}$ is a nonvoid set of vertices and $E$ is a multiset of unordered pairs of distinct elements of ${V}$. The number of times an edge occurs in $G$ is called its multiplicity and edges with multiplicity more than one are called multiple edges.
\end{Definition}

\begin{Definition}\textbf{General graph}.
A general graph $G$ is a pair $({V}, E)$, where ${V}$ is a nonvoid set of vertices and $E$ is a multiset of  unordered pairs of elements of ${V}$. We denote by $uv$ an edge from the vertex $u$ to the vertex $v$. An edge of the form $uu$ is called loop of $G$ and edges which are not loops are known as proper edges. The cardinalities of $ V $ and $ E $ are known as order and size of $ G$, respectively.
\end{Definition}

\begin{Definition}\textbf{Subgraph of a graph}.
Let $G=({V}, E)$ be a graph, $H=(W, E^{\prime})$ is the subgraph of $G$ whenever $W\subseteq V$ and $E^{\prime}\subseteq E$. If $W=V$ the subgraph is said to be spanning subgraph of $ G $. An induced subgraph $\langle W\rangle$ is the subset of $V$ together with all the edges of $G$ between the vertices in $W$.
\end{Definition}
\begin{Definition}\textbf{Complement}. Let $ G $ be a graph, the complement of $ G $, denoted by $ \overline{G} $, is a graph on same set of vertices such that two vertices of $ \overline{G} $ are adjacent if and only if they are not adjacent in $ G. $

\end{Definition}

\begin{Definition}\textbf{Bipartite graph}.
A graph $G$ is said to be bipartite, if its vertex set $V$ can be partitioned into two subsets, say $V_1$ and $V_2$ such that each edge has one end in $V_1$ and other in $V_2$. It is denoted by $ K_{a,b} $, where $ a $ are $ b $ are the cardinalities of $ V_{1} $ and $ V_{2} $, respectively.
\end{Definition}

\begin{Definition}
$\bf (i) $ \textbf{Degree}. Degree of a vertex $v$ in a graph $G$ is the number of edges incident on $v$ and is denoted by $d_v$ or $d(v)$.\\
$\bf (ii) $ \textbf{Conjugate degree}. Let $ d_{v_{k}}$ be the degree of vertex $ v_{k}\in V $. Then conjugate degree of $ v_{k} $ is denoted by $ d_{v_{k}}^{*} $ and is defined as $ d_{v_{k}}^{*}=|\{v_{i}: d_{v_{i}}\geq k \}|$, where $ |~.~| $ denotes cardinality of set.
\end{Definition}

\begin{Definition}\textbf{ $\textbf{r}$-Regular graph}. A graph $G(V,E)$ is said to be $r$-regular if for every vertex $v\in V$, $d_v=r$.
\end{Definition}

\begin{Definition}\textbf{Path}. A path of length $n-1$ $(n\geq 2)$, denoted by $P_n$, is a graph with $n$ vertices ${v_1,v_2,\ldots,v_n}$ and with $n-1$ edges $(v_j, v_{j+1})$, where $j=1,2,\ldots,n-1$.
\end{Definition}

\begin{Definition}\textbf{Cycle}. A cycle of length $n$, denoted by $C_n$, is the graph with vertex set ${v_1,v_2,\ldots,v_n}$ having edges $(v_j, v_{j+1})$, $j=1,2,\ldots,n-1$ and $(v_n, v_1)$.
\end{Definition}

\begin{Definition} \textbf{Connectedness in graphs}. A graph $G(V,E)$ is said to be connected if for every pair of vertices $u$ and $v$ there is a path from $ u $ to $ v $.
\end{Definition}
\begin{Definition}\textbf{Complete graph} A graph in which every pair of distinct vertices is connected by an edge. It is denoted by $ K_{n}. $ Complement of a connected graph is an empty graph.

\end{Definition}

\begin{Definition}\textbf{Tree}. A Tree $ T $ is a connected acyclic graph. 
\end{Definition}

\begin{Definition}\textbf{Matching}. Let $G$ be a graph of order $n$ and size $m$. A $k$-matching of $G$ is a collection of $k$ independent edges (that is, edges which do not share a vertex) of $G$. A $ k $-matching is called perfect if $ n=2k. $
\end{Definition}

\begin{Definition}\textbf{Clique}. Let $G$ be a graph with $n$ vertices and $m$ edges. A maximal complete subgraph of $G$ is called clique of $G$ and its order is known as clique number of $G$, denoted by $\omega(G)$ (or simply $ \omega $).
\end{Definition}

\begin{Definition}\textbf{Independent set}. A vertex subset $W$ of $G$ is said to be an independent set of $G$, if the induced subgraph $<W>$ is a void graph. An independent set of $G$ with maximum number of vertices is called a maximum independent set of $G$ and number of vertices in such a set is called the independence number of $G$ and is denoted by $\alpha(G)$.
\end{Definition}

\begin{Definition} \textbf{Girth of graph}. The girth of a graphs $G$ is the length of smallest cycle of $ G $ and is denoted by $ g $.
\end{Definition}

\begin{Definition}\textbf{C-cyclic graph}. If $ G $ contains $ n $ vertices and $ n+c-1 $ edges, then $ G $ is called a c-cyclic graph.
\end{Definition}

\begin{Definition}\textbf{Diameter}. The distance between two vertices in $ G $ is the number of edges in a shortest path connecting them. The eccentricity of $ v\in V $ is the greatest distance between $ v $ and any other vertex. The minimum and maximum eccentricity of any vertex in $ G $ are known as radius and diameter, respectively.
\end{Definition}

\begin{Definition} \textbf{Joined union}. Let $G(V,E)$ be a graph of order $n$ and $ G_{i}(V_{i}, E_{i})$ be graphs of order $n_i,$ where $i=1,\ldots, n $. The  \textit{joined union} $ G[G_{1},\ldots, G_{n}] $ is the graph $ H(W, F) $ with
\begin{equation*}
W=\bigcup_{i=1}^{n}V_{i}~~~\text{and}~~F=\bigcup_{i=1}^{n}E_{i}~\bigcup\bigcup_{\{v_{i}, v_{j}\}\in E}V_{i}\times V_{j}.
\end{equation*}
In other words, the joined union is obtained from the union of graphs $ G_{1},\ldots, G_{n} $ by joining an edge between each pair of  vertices from $ G_{i} $  and $ G_{j} $ whenever $ v_{i} $ and $ v_{j} $ are adjacent in $ G. $ Thus, the usual join of two graphs $ G_1 $ and $G_2$ is a special case of the joined union $ K_{2}[G_1, G_2]=G_{1}\triangledown G_{2} $ where $ K_{2} $ is the complete graph of order $2$. From the definition, it is clear that the diameter of the graph $G[G_{1},\ldots, G_{n}] $ is same as the diameter of the graph $G$, for every choice of the connected graphs  $G_1,G_2,\dots,G_n$, $n\geq 3$ and $G \ncong K_{n}$.
\end{Definition}

\begin{Definition} \textbf{Group}. A non-empty set $\mathcal{G}$ together with a binary operations $*$, is called a group if it satisfies the following properties.
\begin{itemize}
\item [\bf(i)]  $*$ is associative on $ \mathcal{G} $, that is $ a*(b*c)=(a*b)*c $ for each $ a,b $ and $ c $ in $ \mathcal{G}. $
\item [\bf(ii)]  There is a unique element $ e\in \mathcal{G} $, such that for any $ a\in \mathcal{G} $, $ a*e=a=e*a. $ Such an element $ e $ is known as \textit{identity} of $ \mathcal{G}. $
\item [\bf(iii)] For every $ a\in \mathcal{G} $, there exists a unique $ b\in \mathcal{G} $ such that $a*b=e=b*a$. This element $ b $ is known as inverse of $ a $ and is usually denoted by $ a^{-1} $.

\end{itemize}
\end{Definition}
\begin{Definition} \textbf{Abelian group}. A group $\mathcal{G}$ is called abelian if $ a*b=b*a $ for every $ a $ and $ b $ in $ \mathcal{G} $, otherwise $ \mathcal{G} $ is called non-abelian.
\end{Definition}
\begin{Definition} \textbf{Finite group}. A group $\mathcal{G}$ having finite number of elements is called a finite group. The order of group is the number of elements in $ \mathcal{G} $ and is denoted by $ |\mathcal{G}| .$
\end{Definition}
\begin{Definition} \textbf{Genrating set}. Let $\mathcal{G}$ be a finite group with identity element $ e $. A subset $ S $ of $ \mathcal{G} $ is called the generating set of $ \mathcal{G} $ if every element of $ \mathcal{G} $ can be expressed as a combination of elements of $ S $ and in such case we write $ \mathcal{G}=\langle S\rangle. $ If $ S $ is a generating set of $ \mathcal{G} $, then the elements of $ S $ are called generators of $ \mathcal{G} $ and we say $ S $ generates $ \mathcal{G} $.
\end{Definition}
\begin{Definition} \textbf{Cyclic group}. If a singleton $\{x\} $ is a generating set of $ \mathcal{G} $ then $ \mathcal{G} $ is called a cyclic group and we write $ \mathcal{G}=\langle x\rangle =\{x^{k}:k\in \mathbb{Z}\}$.
\end{Definition}
\begin{Definition} \textbf{Integer modulo group}. For any positive integer $ n\geq 2 $, let $ \mathbb{Z}_{n} $ denote the set $ \{\overline{0},\overline{1},\dots,\overline{n-1}\} $ of all congruences classes of integers modulo $ n $. It is well known that $ \mathbb{Z}_{n} $ is a cyclic group under addition modulo $ n $. Also every cyclic group of order $ n $ is isomorphic copy of $ \mathbb{Z}_{n}. $
\end{Definition}

\begin{Definition} \textbf{Ring}. A non-empty set $R$ together with two binary operations $(+)$ and $(\cdot)$, called addition and multiplication, respectively is called a ring if it satisfies the following properties.
\begin{itemize}
\item[ {\bf(i)}]  $(R,+)$ is an abelian group.
\item[ {\bf(ii)}]  $(R,\cdot)$ is a semi-group.
\item[ {\bf(iii)}]  Distributive laws hold, that is for all $a,b,c\in R$, 
\[\left\{
	\begin{array}{ll}
		a.(b+c)=a.b+a.c, \\
		(a+b).c=a.c+b.c.
	\end{array}
\right. \]
\end{itemize}
\end{Definition}

\begin{Definition}\textbf{Commutative ring}. A ring $R$ is said to be commutative if the semi-group $(R,\cdot)$ is commutative, that is, $a\cdot b=b\cdot a$, for all $a,b\in R$.
\end{Definition}

\begin{Definition}\textbf{ Ring with identity}.  A ring $R$ is said to have an identity (or contains 1) if there is an element $1\in R$ such that $$1\cdot a=a\cdot 1=a, ~\text{for~all}~a\in R.$$
We simply write $ab$ rather than $a\cdot b$ for $a,b\in R$.
\end{Definition}

\begin{definition}\textbf{Unit element}. Let $R$ be a ring with identity $1\ne 0$. An element $u$ in $R$ is called a unit in $R$ if there is some $v$ in $R$ such that $uv=vu=1$. The set of all units of $R$ is denoted by $R^{\times}$.
\end{definition}

\begin{Definition} \textbf{Zero divisor}. Let $ R $ be a commutative ring. A nonzero element $a$ of $R$ is called a zero divisor if there is a nonzero element $b$ in $R$ such that $ab=0$.
\end{Definition}

\newpage
\chapter{Laplacian eigenvalues of graphs and Brouwer's conjecture}

In this chapter, we obtain  upper bounds for the sum of $ k $ largest Laplacian eigenvalues of certain families of graphs. We use these bounds to verify the truth of Brouwer's conjecture for these family of graphs. Besides, we also generalize the results of  Chen \cite{chen} and validate Brouwer's conjecture for graphs with conditions on the size of graph.
\section{Introduction}
Let $G$ be a graph with $n$ vertices $v_1,v_2,\ldots,v_n$ and $m$ edges. If vertex $ v_{i} $ is adjacent to vertex $ v_{j} $, then we write $ v_{i}\sim v_{j} $, otherwise we write $ v_{i}\nsim v_{j} $.  The adjacency matrix of $G$ is the $n\times n$ matrix $A=A(G)=(a_{ij})$, where\\
\begin{equation*}
a_{ij}=\left\{\begin{array}{lr}1, &\mbox{if  $v_i\sim v_j,$}\\
0, &\mbox {if  $v_{i}\nsim v_{j} $}.
\end{array}\right.
\end{equation*} \\
\indent The diagonal matrix associated to $G$ is $D(G)={diag}(d_1, d_2, \dots, d_n)$, where $d_i=\deg(v_i),$ for all $i=1,2,\dots,n$. The matrix $L(G)=D(G)-A(G)$ is known as Laplacian matrix and its spectrum is the Laplacian spectrum of the graph $G$. $ L(G) $ can also be defined as
$L(G)=(l_{ij})$, where\\
\begin{equation*}
l_{ij}=\left\{\begin{array}{lr} d_{i}, &\mbox{if  $v_i= v_j,$}\\
-1, &\mbox{if  $v_i\sim v_j,$}\\
0, &\mbox {if  $v_{i}\nsim v_{j} $}.
\end{array}\right.
\end{equation*} \\
Clearly, $ L(G) $ is real symmetric and positive semi-definite matrix, so its eigenvalues are denoted and ordered as $0=\mu_n\leq\mu_{n-1}\leq\dots\leq\mu_1$, which is the Laplacian spectrum of $G$. It is well known that $\mu_{n-1}=0$ with multiplicity equal to the number of connected components of $G$ and $\mu_{n-1}>0$ if and only if $G$ is connected. For $k=1,2,\dots,n$, let
\begin{align}\label{sk}
S_k(G)= & \sum\limits_{i=1}^{k}\mu_i,
\end{align}
be the sum of the $k$ largest Laplacian eigenvalues of $G$.\\
\indent We note that the sum $S_k(G)$ defined by \eqref{sk} is of much interest by itself and some exciting details, extensions and open problems about it can be found in Nikiforov \cite{nik15}. Further, we observe that the investigations of the parameter $S_k(G)$ may turn out to be useful in the study of several fundamental problems in spectral graph theory. In particular, $S_k(G)$ has a strong relation with the extensively studied spectral parameter
\begin{align*}
LE(G) = \sum_{i=1}^n \left |\mu_i -\frac{2m}{n}\right |,
\end{align*}
defined by Gutman and Zhou \cite{gz} as the \textit{Laplacian energy} of $G$ (discussed in chapter $ 3 $). 

Let $\mathbb{M}_n(\mathbb{C})$ be the set of all square matrices of order $n$ with entries from complex field $ \mathbb{C} $. For $ M\in \mathbb{M}_n(\mathbb{C}) $, the square roots of the eigenvalues of $MM^*$ or $M^{*}M$, where $M^*$ is the complex conjugate of $M$ are known as \emph{singular values}. As $ MM^* $ is positive semi-definite, so singular values of $ M $ are non negative, denoted by $s_1(M)\ge s_2(M)\ge\cdots\ge s_n(M)$.
The \emph{Schatten p-norm} of a matrix $M\in \mathbb{M}_n(\mathbb{C})$ is the $ p $-th root of the sum of the $ p $-th powers of the singular values, that is  
\[ \|M\|_{p}=(s_{1}^{p}(M)+s_{2}^{p}(M)+\dots+s_{n}^{p}(M))^{\frac{1}{p}}, \]
where $n \geq p\geq 1 $. For $ p=2 $, Schatten $ 2 $-norm is \emph{Frobenius norm}, denoted by $ \| M\|_{F} $, defined by 
\[ \|M\|_{F}^{2}=s_{1}^{2}(M)+s_{2}^{2}(M)+\dots+s_{n}^{2}(M). \]
The Schatten $ 1 $-norm is the sum of all singular values and is known as \emph{trace norm} of $ M. $
The sum of first $ k $ singular values is the \emph{Ky Fan k-norm or nuclear norm}, that is for $ p=1 $ and $ n\geq k\geq 1 $, we have 
\[ \|M\|_{k}=s_{1}(M)+s_{2}(M)+\dots+s_{k}(M). \]
$ \| M\|_{1} $ is the largest singular value of $ M $ and is called \emph{spectral norm}.
It is well known that for a \emph{Hermition matrix} $M$, $ s_i(M)=|\lambda_i(M)|$, and for positive semi-definite matrix $ M $, $ s_i(M)=\lambda_i(M)$  where $\lambda_i(M)$, $ i=1,2,\dots,n,$  are the eigenvalues of $ M. $ Thus, $ S_{k}(G) $ is actually the Ky Fan k-norm of $ G $. It is important problem in Linear Algebra to characterize the linear operators having maximum and minimum norms. Similarly, $ LE(G) $ is actually the trace norm of the symmetric matrix $ L(G)-\frac{2m}{n}I_{n} $, where $ I_{n} $ is the identity matrix.

The parameter $S_k(G)$ is  also of great importance in the well known theorem by Grone-Merris \cite{gm}, a nice proof of which is due to  Bai in \cite{b}, its statement is given in the following theorem.
\begin{theorem} \label{Theorem 1.1.} (Grone-Merris-Bai Theorem) If $G$ is any graph of order $n$ and $k$, $1\leq k\leq n,$ is any positive integer, then
$S_k(G)\leq \sum\limits_{i=1}^{k}d_i^*(G),$ where $d_i^*(G)=|\{v\in V(G):d_{v}\geq i\}|$, for $i=1,2,\dots,n$.
\end{theorem}

A strong observation from Theorem \ref{Theorem 1.1.} is that $LE(G)\leq \sum_{i=1}^n \left |d_i^* - \frac{2m}{n}\right |$. We note that equality is attained by threshold graphs and therefore it is worth to find the threshold graph having the largest Laplacian energy. Indeed this has been done by Helmberg and Trevisan \cite{hel15} with the main tool being Theorem \ref{Theorem 1.1.}. These observations have increased the evidence to a belief that threshold graphs have large Laplacian energy. In this family, it supported that the candidate graph for the largest Laplacian energy is a particular pineapple graph. The \emph{pineapple} graph is obtained by appending pendent edges to a vertex of a complete graph.

Analogous to Grone-Merris-Bai Theorem, Brouwer \cite{bh} conjectured the upper bound for $ S_{k}(G) $, which is known as Brouwer's conjecture, and is stated as follows.

\begin{conjecture}[Brouwer's conjecture]\label{Conjecture 1.2.} If $G$ is any graph with order $n$ and size
$m$, then
\[ S_k(G)\leq m+\binom{k+1} {2}, ~~\text{for any} ~~k \in \{1,2,\dots,n\}. \]
\end{conjecture}

In order to move forward, it is important to know the progress of  Brouwer's conjecture. Mayank \cite{Mayank} (see also \cite{ht}) proved that the split graphs and the cographs satisfy Brouwer's conjecture. By using computer computations, Brouwer \cite{bh} checked this  conjecture for all graphs with at most $10$ vertices. For $k=1$, Brouwer's conjecture follows from the well-known inequality $\mu_1(G) \leq n $ and the cases $k=n$ and $k=n-1$ are straightforward. Haemers et al. \cite{hmt} showed that Brouwer's conjecture is true for all graphs when $k=2$ and is also true for trees. Du et al. \cite{dz} obtained various upper bounds for $S_k(G)$ and proved that Brouwer's conjecture is also true for unicyclic and bicyclic graphs.  Rocha and Trevisan \cite{rt} obtained various upper bounds for $S_k(G)$, which improve the upper bounds obtained in \cite{dz} for some cases and proved that Brouwer's conjecture is true for all $k$ with $1\leq k\leq \lfloor\frac{g}{5}\rfloor$, where $g$ is the girth of the graph $G$. They also showed that Brouwer's conjecture is true for a connected graph $G$ having maximum degree $\Delta$ with $p$ pendant vertices and $c$ cycles for all $\Delta\geq c+p+4$.  Ganie, Alghamdi and Pirzada, \cite{hap} obtained upper bounds for $S_k(G)$ in terms of various graph parameters, which improve some previously known upper bounds and showed that Brouwer's conjecture is true for some new families of graphs. For the further progress on Brouwer's  Conjecture, we refer to \cite{chen,hap,ht} and the references therein. However, Brouwer's conjecture still remains open at large.

A very interesting and useful lemma due to Fulton \cite{ful} is as follows.

\begin{lemma}\label{Lemma 2.1.}
Let $A$ and $B$ be two real symmetric matrices both of order $n$. If $k$, $1\leq k \leq n,$ is a positive integer, then
$\sum\limits_{i=1}^{k}\lambda_i(A+B)\leq \sum\limits_{i=1}^{k}\lambda_i(A) + \sum\limits_{i=1}^{k}\lambda_i(B),$ where $\lambda_i(X)$ is the $i^{th}$ eigenvalue of $X$.
\end{lemma}

The following upper bound for the sum of the $k$ largest Laplacian eigenvalues of a tree $T$ can be found in \cite{fhrt_1}.

\begin{lemma}\label{Lemma 2.2.}
Let $T$ be a tree with $n\geq 2$ vertices. If $S_k(T)$ is the sum of the $k$ largest Laplacian eigenvalues of $T$, then $S_k(T)\leq n-2+2k-\frac{2k-2}{n}$, for $1\leq k\leq n$. For $k=1$, equality   occurs when $G\cong S_{n}$.
\end{lemma}

Assume that a graph $G$ has a kind of symmetry so that its associated matrix is written in the form
\begin{align}\label{7}
M=\begin{pmatrix}
X & \beta & \cdots &\beta & \beta \\
\beta^T & B & \cdots & C & C\\
\vdots & \vdots & \cdots & \vdots &\vdots\\
\beta^T & C & \cdots &B & C\\
\beta^T & C & \cdots & C & B
\end{pmatrix},\end{align}
where $X\in  R^{t\times t},~\beta \in R^{t\times s}$ and $B,C\in R^{s\times s}$, such that $n=t+ cs$, where $c$ is the number of copies of $B$.   Then the spectrum of this matrix can be obtained as the union of the spectrum of smaller matrices using the following technique given in \cite{ft}. In the statement of the theorem, $\varrho^{(k)}(Y)$ indicates the multi-set formed by $k$ copies of the spectrum of $Y,$ denoted by $\varrho(Y)$.

\begin{lemma}\label{Lemma 3.4.}
Let $M$ be a matrix of the form given in (\ref{7}), with $c\geq 1$ copies of the
block $B$. Then
\begin{enumerate}
\item[\bf{(i)}] $\varrho(B-C)\subseteq \varrho(M)$ with multiplicity $c-1$,
\item[\bf (ii)]  $\varrho(M)\setminus \varrho^{(c-1)}(B-C)=\varrho(M^{'})$ is the set of the remaining $t+s$ eigenvalues of $M$, where
$M^{'}=\begin{pmatrix}
X & \sqrt{c}.\beta  \\
\sqrt{c}.\beta^T & B+(c-1)C
\end{pmatrix}.$
\end{enumerate}
\end{lemma}

\indent The next useful lemma found in \cite{w}.

\begin{lemma}\label{Lemma 3.5.}
Let $X$ and $Y$ be Hermitian matrices of order $n$ such that $Z=X+Y$. Then
\begin{align*}
\lambda_t(Z)\leq \lambda_j(X)+\lambda_{t-j+1}(Y),~~ n\geq t\geq j\geq 1,\\
\lambda_t(Z)\geq \lambda_j(X)+\lambda_{t-j+n}(Y), ~~n\geq j\geq t\geq 1,
\end{align*}
where $\lambda_i(M)$ is the $i^{th}$ largest eigenvalue of the matrix $M$.
\end{lemma}

\section{Upper bounds for $S_{k}(G)$ }
Now, we obtain an upper bound for $S_k(G)$ in terms of the clique number $\omega$, the order $n$ and integers $p\geq 0,~r\geq 1,~s_1\geq s_2\geq 2$ associated to the structure of the graph $G$ for certain types of graphs.

\begin{theorem}\label{Theorem 2.3.}
Let $G$ be a connected graph of order $n\geq 4$ and size $m$ having clique number
$\omega\geq 2$. If $H=G\setminus K_{\omega}$ is a graph having $r$ non-trivial components $C_1,C_2,\ldots, C_r$, each of which is a c-cyclic graph and $p\geq 0$ trivial components, then
\begin{align}
S_k(G)\leq\left\{\begin{array}{lr}\omega(\omega-1)+n-p+2r(c-1)+2k, &\mbox{if $k\geq \omega -1$},\\
k(\omega+2)+n-p+2r(c-1), &\mbox{if $k\leq \omega-2$}.
\end{array} \right.
\end{align}

\end{theorem}
\noindent{\bf Proof.} Consider the connected graph $G$ of order $n\geq 4$ and size $m$. Let $\omega\geq 2$
be the clique number of $G$. Clearly, $K_{\omega}$ is a subgraph of $G$. If we remove
the edges of $K_{\omega}$ from $G$, the Laplacian matrix of $G$ can be
decomposed as $Q(G)=Q(K_{\omega}\cup(n-\omega)K_1)+Q(H)$, where
$H=G\setminus K_{\omega}$ is the graph obtained from $G$ by removing the edges of
$K_{\omega}$.  Applying Lemma \ref{Lemma 2.1.} and using the fact that
$S_k(K_{\omega}\cup(n-\omega)K_1)=S_k(K_{\omega})$, for $1\leq k\leq n$, we have 
\begin{align*}
S_{k}(G)=\sum\limits_{i=1}^{k}\mu_i(G)\leq \sum\limits_{i=1}^{k}\mu_i(K_{\omega})+\sum\limits_{i=1}^{k}\mu_i(H)
=S_k(K_{\omega})+S_{k}(H).
\end{align*}
\noindent It is well known that the Laplacian spectrum of $K_{\omega}$ is
$\{\omega^{[\omega-1]},0\}$, so that $S_k(K_{\omega})=k\omega$. For
$k\geq \omega-1 $, it is better to consider $S_k(K_{\omega})\leq \omega(\omega-1)$
and for $k\leq \omega-2$, it is better to consider $S_k(K_{\omega})\leq k\omega$.
Therefore, it follows that
\begin{align}\label{2}
S_k(K_{\omega})\leq\left\{\begin{array}{lr}\omega(\omega-1), &\mbox{if $k\geq \omega-1$},\\
k\omega, &\mbox{if $k\leq \omega-2$}.
\end{array} \right.
\end{align}
In the hypothesis, it is given that $H=G\setminus K_{\omega}$ is a  graph
having $r\geq 1$ non-trivial components, each of which is a c-cyclic graph and $p$ trivial components.
Let $C_1,C_2,\dots,C_r$ be the non-trivial components  of $H$ with $|C_i|=n_i\geq 2$ and $p$
be the number of isolated vertices in $H$. That is, $H=H^{'}\cup p K_1$, where
$H^{'}=C_1\cup C_2\cup\cdots\cup C_r$. Clearly $S_k(H)=S_k(H^{'})$ and
$\sum\limits_{i=1}^{r}n_i=n-p$. Let $k_i$ be the number of the first $k$
largest Laplacian eigenvalues of $H^{'}$ that are Laplacian eigenvalues of $C_i$,
 where $0\leq k_i\leq k,~~1\leq i\leq r$ and $\sum\limits_{i=1}^{r}k_i=k$. Since $C_i$ is a c-cyclic graph on $n_i$ vertices and $m_i=n_i+c-1$ edges, it follows that the Laplacian matrix of $C_i$ can be written as
 $L(C_i)=L(T_i)+L(cK_2\cup (n_i-2c)K_1)$, where $T_i$ is the spanning tree of $C_i$. Therefore, applying Lemma \ref{Lemma 2.2.} to $C_i$, and using Lemma \ref{Lemma 2.1.}, we get
 \begin{align}\label{3}
 S_{k_i}(C_i)\leq S_{k_i}(T_i)+2c\leq n_i-2+2k_i+2c.
 \end{align}
\indent Now, applying Lemma \ref{Lemma 2.2.} to the graph $H^{'}=C_1\cup C_2\cup\cdots\cup C_r$ and using Inequality \eqref{3}, for all $k_i,~1\leq i\leq r$, it follows that
\begin{align*}
S_k(H^{'})&=S_k(C_1\cup C_2\cup\cdots\cup C_r)\leq\sum\limits_{i=1}^{r}S_{k_i}(C_i)
\\& \leq \sum\limits_{i=1}^{r}(n_i-2+2k_i+2c)\\&
=\sum\limits_{i=1}^{r}n_i+\sum\limits_{i=1}^{r}(2c-2)+\sum\limits_{i=1}^{r}2k_i\\
&=n-p+2r(c-1)+2k.
\end{align*}
This shows that
\begin{align}\label{4}
S_k(H)=S_k(H^{'})\leq n-p+2r(c-1)+2k.
\end{align}

The result now follows from  Inequalities \eqref{2} and \eqref{4}.\qed

\indent In particular, taking $c=0$, we have the following observation, which was obtained in \cite{hpv}.

\begin{corollary}\label{Corollary 2.4.}
Let $G$ be a connected graph of order $n\geq 4$ and size $m$ having clique number
$\omega\geq 2$. If $H=G\setminus K_{\omega}$ is a forest having $r$ non-trivial components $C_1,C_2,\ldots, C_r$ and $p\geq 0$ trivial components, then
\begin{align}
S_k(G)\leq\left\{\begin{array}{lr}\omega(\omega-1)+n-p-2r+2k, &\mbox{if $k\geq \omega-1$},\\
k(\omega+2)+n-p-2r, &\mbox{if $k\leq \omega-2$}.
\end{array} \right.
\end{align}
\end{corollary}

Now, we obtain an upper bound for $S_k(G)$ in terms of the positive integers
$s_1,~s_2$ and the order $n$ of graph $G$.

\begin{theorem}\label{Theorem 2.5.}
Let $G$ be a connected graph of order $n\geq 4$ and size $m$. Let $K_{s_1,s_2}$,
$s_1\leq s_2\geq 2$, be the maximal complete bipartite subgraph of the graph $G$.
If $H=G\setminus K_{s_1,s_2}$  is a graph having $r$ non-trivial components $C_1,C_2,\ldots, C_r$, each of which is a $c$-cyclic graph and $p\geq 0$ trivial components, then
\begin{align}
S_k(G)\leq\left\{\begin{array}{lr}
2s_1s_2+n-p+2r(c-1)+2k, &\mbox{if $k\geq s_1+s_2-1$},\\
s_2+ks_1+n-p+2r(c-1)+2k, &\mbox{if $k\leq s_1+s_2-2$}.
\end{array} \right.
\end{align}
\end{theorem}
\noindent{\bf Proof.}
Consider the connected graph $G$ with $K_{s_1,s_2}$, $(s_1\geq s_2)$, as its maximal complete bipartite subgraph. If we remove the edges of $K_{s_1,s_2}$ from $G$, the Laplacian matrix of $G$ can be decomposed as $L(G)=L(K_{s_1,s_2}\cup(n-s_1-s_2)K_1)+L(H)$, where $H=G\setminus K_{s_1,s_2}$ is the graph obtained from $G$ by removing the edges of $K_{s_1,s_2}$.  Applying Lemma \ref{Lemma 2.1.} and using the fact that $S_k(K_{s_1,s_2}\cup(n-s_1-s_2)K_1)=S_k(K_{s_1,s_2})$, for $1\leq k\leq n$, we have
\begin{align*}
S_{k}(G)\leq \sum\limits_{i=1}^{k}\mu_i(K_{s_1,s_2})+\sum\limits_{i=1}^{k}\mu_i(H)
=S_k(K_{s_1,s_2})+S_{k}(H).
\end{align*}
\indent  Now, proceeding similarly as in Theorem \ref{Theorem 2.3.} and using the fact that the Laplacian spectrum of $K_{s_1,s_2}$ is $\{s_1+s_2,s_1^{[s_2-1]},s_2^{[s_1-1]},0\}$, the result follows.\qed

\indent In particular, if $s_1=s_2$, we have the following consequence of Theorem \ref{Theorem 2.5.}.

\begin{corollary}\label{Corollary 2.6.}
Let $G$ be a connected graph of order $n\geq 4$ and size $m$ and let $K_{s,s},~~s\geq 2$, be the maximal complete bipartite subgraph of graph $G$. If $H=G\setminus K_{s,s}$ is a graph having $r$ non-trivial components $C_1,C_2,\ldots, C_r$, each of which is a c-cyclic graph and $p\geq 0$ trivial components, then 
\begin{align*}
S_k(G)\leq\left\{\begin{array}{lr}2s^2+n-p+2r(c-1)+2k, &\mbox{if $k\geq 2s-1$},\\
s+k(s+2)+n-p+2r(c-1)+2k, &\mbox{if $k\leq 2s-2$}.
\end{array} \right.
\end{align*}
\end{corollary}

\section{Brouwer's conjecture for some classes of graphs}

\label{sec:brou}

\indent This section is devoted to verify the truth of Brouwer's conjecture for new families of graphs.

\begin{theorem}\label{Theorem 3.1.}
Let $G$ be a connected graph of order $n\geq 4$ and size $m$ having clique number
$\omega\geq 2$. If $H=G\setminus K_{\omega}$ is a graph having $r$ non-trivial components $C_1,C_2,\ldots, C_r$, each of which is a c-cyclic graph, then $S_k(G)\leq m+\dfrac{k(k+1)}{2}$, for all $k\in [1, \Delta_{1}]$ and $k\in [\beta_1,n]$, where $\Delta_1=min\{\omega-2,\gamma_{1}\}, \gamma_{1}=\frac{2\omega+3-\sqrt{16\omega+8r(c-1)+9}}{2}$ and $\beta_1=\frac{3+\sqrt{4\omega^{2}-4\omega+8r(c-1)+9}}{2}$.
\end{theorem}
\noindent{\bf Proof.}
Let $G$ be a connected graph of order $n$ and size $m$ having clique number $\omega\geq 2$. If $H=G\setminus K_{\omega}$ is a graph having $r$ non-trivial components $C_1,C_2,\ldots, C_r$, each of which is a c-cyclic graph and $p\geq 0$ trivial components, then $m=\frac{\omega(\omega-1)}{2}+n-p+r(c-1)$. For $k\geq \omega-1$, from Theorem \ref{Theorem 2.3.}, we have
\begin{align*}
S_k(G)&\leq \omega(\omega-1)+n-p+2r(c-1)+2k\\
& \leq m+\frac{k(k+1)}{2}=\frac{\omega(\omega-1)}{2}+n-p+r(c-1)+\frac{k(k+1)}{2},
\end{align*} if
\begin{align}\label{5a}
k^2-3k-(\omega(\omega-1)+2r(c-1))\geq 0.
\end{align}
Now, consider the polynomial $f(k)=k^2-3k-(\omega(\omega-1)+2r(c-1))$. The roots of this polynomial are $\beta_1=\frac{3+\sqrt{4\omega^2-4\omega+8r(c-1)+9}}{2}$ and $\beta_2=\frac{3-\sqrt{4\omega^2-4\omega+8r(c-1)+9}}{2}$. This shows that $f(k)\geq 0$, for all $k\geq \beta_1$ and for all $k\leq \beta_{2}$. Since $1\leq r\leq \omega$, it can be seen that $2-\omega<\beta_2 < 3-\omega$. Thus, it follows that (\ref{5a}) holds for all $k\geq \beta_1$. It is easy to see that $ \beta_{1}\geq \omega-1 $. This completes the proof in this case.\\

For $k\leq \omega-2$, from Theorem \ref{Theorem 2.3.}, we have
\begin{align*}
S_k(G)&\leq (\omega+2) k+n-p+2r(c-1) \leq m+\frac{k(k+1)}{2}\\
&=\frac{\omega(\omega-1)}{2}+n-p+r(c-1)+\frac{k(k+1)}{2},
\end{align*} if
\begin{align}\label{6}
k^2-(2\omega+3)k-2r(c-1)+\omega(\omega-1)\geq 0.
 \end{align}
Proceeding, similarly as above, it can be seen that \eqref{6} holds for all $k\leq\gamma_{1}= \frac{2\omega+3-\sqrt{16\omega+9+8r(c-1)}}{2}$. Indeed, $ \gamma_{1}\leq \omega-2 $, holds for any $ c\geq 1 $, completing the proof in this case as well.\qed

Evidently, if $ c=0 $, then $ \Delta_{1}=\omega-2 $ for $ 2\leq \omega \leq 5 $; and $ \Delta_{1}=\gamma_{1} $ for $ \omega\geq 6 $. Further, if $ c=1 $, then $ \Delta_{1}=w-2 $ for $ \omega=2 $; and $ \Delta_{1}=\gamma_{1} $ for $ \omega\geq 3 $. For $ c\geq 2 $, clearly $ \Delta_{1}=\gamma_{1} $.

\indent From Theorem \ref{Theorem 3.1.}, we have the following observation.
\begin{corollary}\label{corollary 3.2.}
Let $G$ be a connected graph of order $n\geq 4$ and size $m$ having clique number
$\omega\geq 2$. Let $H=G\setminus K_{\omega}$ be a graph having $r$ non-trivial components $C_1,C_2,\ldots, C_r$, each of which is a $c$-cyclic graph and $\alpha=\sum\limits_{i=1}^{r}\frac{2k_i-2}{n_i}.$\\
{\bf (i)} If $c=0$, that is, each $C_i$ is a tree, then Brouwer's conjecture holds for all $k\in [1, \Delta_{1}]$ and $k\in [\omega +1, n]$, where $ \Delta_{1}=min\{\omega-2,\gamma_{1}\} $ and $\gamma_{1}=\frac{2\omega+3-\sqrt{16\omega+9-8r}}{2}$.\\
{\bf (ii)} If $c=1$, that is, each  $C_i$ is a unicyclic graph, then Brouwer's conjecture holds for all $k\in [\omega +2, n]$ and $k\in [1,\triangle_{1}],$ where $ \triangle_{1}=\min\{\omega-2,\gamma_{1}\} $ and  $\gamma_{1}=\frac{2\omega+3-\sqrt{16\omega+9}}{2} $.\\
{\bf (iii)} If $c=2$, that is, each  $C_i$ is a bicyclic graph, then  Brouwer's conjecture holds for all $k\in [\omega +3, n]$ and $k\in \left [1,\frac{2\omega+3-\sqrt{16\omega+9+8r}}{2}\right ]$.\\
{\bf (iv)} If $c\geq 3$, that is, each  $C_i$ is a c-cyclic graph, then Brouwer's conjecture holds for all $k\in [\omega +c, n]$ and $k\in \left [1,\frac{2\omega+3-\sqrt{16\omega+9+8r(c-1)}}{2}\right ]$.
\end{corollary}
\noindent{\bf Proof.} {\bf (i)}. If $c=0$, then $\beta_1=\frac{3+\sqrt{4\omega^2-4\omega-8r+9}}{2}$ and $\gamma_1=\frac{2\omega+3-\sqrt{16\omega+9-8r}}{2}$. Using the fact that $1\leq r \leq \omega$, we have
\begin{align*}
\beta_1=\frac{3+\sqrt{4\omega^2-4\omega-8r+9}}{2}\leq \frac{3+\sqrt{4\omega^2-4\omega+1}}{2}=\omega+1.
\end{align*}
{\bf (ii)}. If $c=1$, then $\beta_1=\frac{3+\sqrt{4\omega^2-4\omega+9}}{2}$ and $\gamma_1=\frac{2\omega+3-\sqrt{8\omega+9}}{2}$. We have
\begin{align*}
\beta_1=\frac{3+\sqrt{4\omega^2-4\omega+9}}{2}\leq\frac{3+\sqrt{4\omega^2+4\omega+1}}{2}=\omega+2.
\end{align*}
{\bf (iii)}. Proceeding similarly as in part (i) and (ii), we can prove part (iii). \\
{\bf (iv)}.  If $c\geq 3$, then using $r\leq \omega$, we have
\begin{align*}
\beta_1&=\frac{3+\sqrt{4\omega^2-4\omega+8r(c-1)+9}}{2}\leq \frac{3+\sqrt{4\omega^2+4\omega(2c-3)+9}}{2}\\&
= \frac{3+\sqrt{(2\omega+(2c-3))^2-4c(c-3)}}{2} \leq \omega+c.
\end{align*}\qed

Now, we consider some special classes of graphs satisfying the hypothesis of Theorem \ref{Theorem 3.1.}. Let $C_{\omega}(a,a,\dots,a),~~a\geq 1$, be the family
of connected graphs of order $n=\omega(a+1)$ and size $m$ obtained by identifying one of the  vertex of a $c$-cyclic graph $C$ of order $a+1$  to each vertex of the clique $K_{\omega}$. For the family of graphs $C_{\omega}(a,a,\dots,a)$, we see that Brouwer's conjecture is true for various subfamilies depending upon the value of $c$, the order of the $c$-cyclic graphs, and clique number of the graph.

\begin{theorem}\label{Theorem 3.8.}
For the graph $G\in C_{\omega}(a,a,\dots,a), ~~a\geq 1$, Brouwer's conjecture holds for all $k$, if $c=0$. If $c=1$, Brouwer's conjecture holds for all $k, k\neq \omega+1$, provided $a\leq \omega+1 $. If $ c=2, $ Brouwer's conjecture holds for all $ k, k\notin \{\omega+1,\omega+2\} $, holds for $ k=\omega+1 $, provided $ a\leq \omega +2 $ and holds for $ k=\omega +2, $ provided $ a\leq 2\omega +\frac{1}{2} $. If $ c\geq 3, $ Brouwer's conjecture holds for all $ k $, provided $a\leq \omega -\frac{1}{2}+\sqrt{(2c-1)\omega}$.
\end{theorem}
 \noindent{\bf Proof.}
Consider a connected graph $G$ with order $n$ and size $m$ and let $G$ belong to
the family $C_{\omega}(a,a,\dots,a),~~a\geq 1$. For $a\ne 1$, let $C$ be a $c$-cyclic graph on
$a+1$ vertices having one of the pendent vertices fused with a vertex of the
clique $K_{\omega}$. Let $L(C)$ be the Laplacian matrix of $C$ and let
$0=\mu_{a+1}\leq \mu_{a}\leq\cdots\leq \mu_1$ be the eigenvalues of the Laplacian matrix $L(C)$.
Let  $C_{q\times q}$ be the matrix defined as
\[ C_{q\times q}=\begin{pmatrix}
-1 & 0 & \cdots  & 0 \\
0 & 0 & \cdots  & 0\\
\vdots & \vdots & \cdots &\vdots\\
0 & 0 & \cdots & 0\\
0 & 0 & \cdots  & 0
\end{pmatrix}_{a+1}, \]  where $q=a+1$.\\
\indent By a suitable labelling of vertices of $G$, it can be seen that the
 Laplacian matrix of $G$ can be written as
\begin{align*}L(G)=\begin{pmatrix}
F & C_{q\times q} & \cdots &C_{q\times q} \\
C_{q\times q} & F & \cdots & C_{q\times q} \\
\vdots & \vdots & \cdots & \vdots \\
 C_{q\times q}& C_{q\times q} & \cdots & C_{q\times q}\\
C_{q\times q} & C_{q\times q} & \cdots &  F
\end{pmatrix}_{\omega},\end{align*}
where $ F=L(C)-(\omega-1) C_{q\times q} $.
\noindent Taking $X$ and $\beta$ to be a matrix of order zero,
$B=L(T)-(\omega-1) C_{q\times q}$ and $C=C_{q\times q}$ in (\ref{7}) and using Lemma
\ref{Lemma 3.4.}, we have $\sigma(L(G))=\sigma^{(\omega-1)}(L(C)-\omega
C_{q\times q})\cup \sigma(L(C)$. The eigenvalues of the matrix
$-\omega C_{q\times q}$ are $\omega$ with multiplicity one and $0$ with
multiplicity $a$. If $\mu_{a+1}^{'}\leq \mu_{a}^{'}\leq \cdots\leq \mu_{1}^{'}$ are
the eigenvalues of the matrix $L(C)-\omega C_{q\times q}$, we see that the
eigenvalues of the matrix $L(G)$ are $\mu_{i}^{'},~~1\leq i\leq a+1$,
each with multiplicity $\omega-1$ and $\mu_{i},~~~1\leq i\leq a+1$.\\
\indent Clearly the graph $G\setminus K_{\omega}$ has $\omega$  non-trivial components, each of which is a $c$-cyclic graph. If $c=0$, then each of the $\omega$ components of $G\setminus K_{\omega}$ is a tree on $a+1$ vertices. 
Therefore, by  Theorem 3.5 in \cite{hpv},  Brouwer's conjecture holds for all $k$, completing the proof in this case.\\  By Lemma \ref{Lemma 2.1.}, we first observe that
\begin{align}\label{11}
\mu_{1}^{'}\leq \omega+ \mu_{1}.
\end{align}
\indent If $c=1$, then each of the $\omega$ components of $G\setminus K_{\omega}$ is a unicyclic graph on $a+1$ vertices. So, by Corollary \ref{corollary 3.2.}, Brouwer's conjecture holds for all $k\geq \omega+2$. Now, we need to show that Brouwer's conjecture holds for all $k\leq \omega+1$. For $1\leq k\leq \omega
+1$, from (\ref{11}), it follows that
 \begin{align*}
 S_k(G)&\leq k \mu_1^{'}\leq k(\omega+ \mu_1)\leq m+\frac{k(k+1)}{2}\\&
 =\frac{\omega(\omega-1)}{2}+n+\frac{k(k+1)}{2},
 \end{align*}
  \begin{align}\label{13}\text{if} \quad
   k^2-(2\mu_1+2\omega-1)k+\omega^2-\omega+2n\geq 0.
  \end{align}
 \noindent Now, consider the polynomial
 \begin{align*}
 f(k)=k^2-(2 \mu_1+2\omega-1)k+\omega^2-\omega+2n.
\end{align*}
The discriminant of this polynomial is $d=(2 \mu_1-1)^2+8\omega \mu_1-8n$ and its roots are
$x=\frac{(2 \mu_1+2\omega-1)+\sqrt{d}}{2}$ and $y=\frac{(2 \mu_1+2\omega-1)-
\sqrt{d}}{2}$. This shows that $f(k)\geq 0$, for all $k\leq x$. We have
\begin{align*}
y=\frac{(2 \mu_1+2\omega-1)-\sqrt{(2\mu_1-1)^2+8\omega \mu_1-8n}}{2}\geq \omega,
\end{align*}
which implies that $2 \mu_1-1\geq  \sqrt{(2\mu_1-1)^2+8\omega q_1-8n}$, further it implies that $\mu_1\leq\frac{n}{\omega}$, which is always true as $\mu_1\leq a+1=\frac{n}{\omega}$. This shows that \eqref{13} holds for all $k,~~~k\leq \omega.$ \\
\indent For $t=2, j=1$ and $i=2$, from Lemma \ref{Lemma 3.5.}, it follows that $ \mu^{'}_2\leq q_1$. Therefore,  for $k=\omega+1$, we have
\begin{align*}
S_{\omega+1}(G)&=(\omega-1) \mu_1^{'}+2\mu_2^{'}\leq (\omega-1)(\omega+ \mu_1)+2\mu_1=\omega(\omega-1)+(\omega+1) \mu_1\\&
\leq m+\frac{(\omega+1)(\omega+2)}{2}
=\frac{\omega(\omega-1)}{2}+n+\frac{(\omega^2+3\omega+2)}{2},
\end{align*}
which is true if $\mu_1\leq \frac{n+2\omega+1}{\omega+1}$. Since $ \frac{n}{\omega}\leq \frac{n+2\omega+1}{\omega+1}$ for $ a\leq \omega+1 $, it follows that Brouwer's conjecture holds in this case for $ a\leq \omega +1 $, completing the proof in this case.\\
\indent   If $c=2$, each of the $\omega$ components of $G\setminus K_{\omega}$ is a bicyclic graph on $a+1$ vertices. So, by Corollary \ref{corollary 3.2.}, Brouwer's conjecture holds for all $k\geq \omega+3$. To complete the proof in this case, we need to show that Brouwer's conjecture holds for all $k\leq \omega+2$. For $1\leq k\leq \omega+2$, from \eqref{11}, it follows that
\begin{align*}
S_k(G)&\leq k \mu_1^{'}\leq k(\omega+ \mu_1)\leq m+\frac{k(k+1)}{2}\\&
=\frac{\omega(\omega+1)}{2}+n+\frac{k(k+1)}{2},
\end{align*}
\begin{align}\label{14}\text{if} \quad
k^2-(2 \mu_1+2\omega-1)k+\omega^2+\omega+2n\geq 0.
\end{align}
\indent Now, consider the polynomial
\begin{align*}
f(k)=k^2-(2\mu_1+2\omega-1)k+\omega^2+\omega+2n.
\end{align*}
The discriminant of this polynomial is $d=(2 \mu_1-1)^2+8\omega \mu_1-8\omega-8n$ and its roots are $x=\frac{(2 q_1+2\omega-1)+\sqrt{d}}{2}$ and $y=\frac{(2 q_1+2\omega-1)-\sqrt{d}}{2}$.
This shows that $f(k)\geq 0$, for all $k\leq y$ and $f(k)\geq 0$, for all
$k\geq x$. We have
\begin{align*}
y=\frac{(2 \mu_1+2\omega-1)-\sqrt{(2 \mu_1-1)^2+8\omega(\mu_1-1)-8n}}{2}\geq \omega,
\end{align*}
which implies that $2\mu_1-1\geq  \sqrt{(2\mu_1-1)^2+8\omega(\mu_1-1)-8n}$, further it implies that $\mu_1\leq\frac{n}{\omega}+1$, which is always true as $\mu_1\leq a+1=\frac{n}{\omega}$. This shows that \eqref{13} holds for all $k,~~~k\leq \omega.$ \\
\indent  For $t=2, j=1$ and $i=2$, from Lemma \ref{Lemma 3.5.}, it follows that $ \mu^{'}_2\leq \mu_1$. Therefore, for $k=\omega+1$, we have
\begin{align*}
S_{\omega+1}(G)&= (\omega-1)\mu_1^{'}+2\mu_2^{'}\leq (\omega-1)(\omega+\mu_1)+2\mu_1=\omega(\omega-1)+(\omega+1)\mu_1\\&
\leq m+\frac{(\omega+1)(\omega+2)}{2}
=\frac{\omega(\omega+1)}{2}+n+\frac{(\omega^2+3\omega+2)}{2},
\end{align*}
which is true if $\mu_1\leq \frac{n+3\omega+1}{\omega+1}$. Since $\frac{n}{\omega}\leq \frac{n+3\omega+1}{\omega+1} $, for $ a\leq \omega+2 $, it follows that Inequality \eqref{13} hods for $ k=\omega+1 $, provided that $ a\leq \omega+2 .$\\
\indent Also, for $k=\omega+2$, we have
\begin{align*}
S_{\omega+2}(G)&= (\omega-1)\mu_1^{'}+3\mu_2^{'}\leq (\omega-1)(\omega+q_1)+3\mu_1=\omega(\omega-1)+(\omega+2)\mu_1\\&
\leq m+\frac{(\omega+1)(\omega+2)}{2}
=\frac{\omega(\omega+1)}{2}+n+\frac{(\omega^2+5\omega+6)}{2},
\end{align*}
which is true if $\mu_1\leq \frac{n+4\omega+3}{\omega+2}$. Proceeding similarly as above it can be seen that $\mu_1\leq \frac{n+4\omega+3}{\omega+2}$ holds, provided that $ a\leq 2\omega+\frac{1}{2} $ and the proof is complete in this case.\\
\indent   Now, suppose that $c\geq 3$. Then each of the $\omega$ components of $G\setminus K_{\omega}$ is a $c$-cyclic graph on $a+1$ vertices. So, by Corollary \ref{corollary 3.2.}, Brouwer's conjecture holds for all $k\geq \omega+c$. To complete the proof, we need to show that Brouwer's conjecture holds for all $k\leq \omega+c-1$.  We consider the cases  $1\leq k\leq \omega$ and $\omega+1 \leq k \leq \omega+c-1$ separately.  For $1\leq k\leq \omega$, for the case $c=2$, we have earlier seen that
\begin{align*}
S_k(G)&\leq k \mu_1^{'}\leq \frac{\omega(\omega-1)}{2}+n+\omega+\frac{k(k+1)}{2}\\&
< \frac{\omega(\omega-1)}{2}+n+\omega(c-1)+\frac{k(k+1)}{2}= m+\frac{k(k+1)}{2},
\end{align*} as $c\geq 3$. This shows that Brouwer's Conjecture holds for $1\leq k\leq \omega$. So, suppose that $\omega+1\leq k \leq \omega+c-1$.  For $\omega+1\leq k \leq \omega+c-1$, from \eqref{11}, it follows that
\begin{align*}
S_k(G)&\leq  (\omega-1)\mu_1^{'}+(k-\omega+1)\mu_2^{'}\leq (\omega-1)(\omega+\mu_1)+(k-\omega+1)\mu_1\\&
\leq m+\frac{k(k+1)}{2}
=\frac{\omega(\omega-1)}{2}+n+\omega(c-1)+\frac{k(k+1)}{2},
\end{align*}
\begin{align}\label{15}\text{if} \quad
k^2-(2\mu_1-1)k+(2c-1)\omega+2n-\omega^2\geq 0.
\end{align}
Since $\omega+1\leq k \leq \omega+c-1$, from \eqref{15}, it follows that
\begin{align*}
(\omega+1)^2-(2\mu_1-1)(\omega+c-1)+(2c-1)\omega+2n-\omega^2\geq 0,
\end{align*}
which implies that \[ 2\mu_1\leq \frac{2n+2c\omega+2\omega+c}{\omega+c-1}. \]
Using the fact that $\mu_1\leq a+1$, it follows that
\begin{align*}
2a+2\leq \frac{2n+2c\omega+2\omega+c}{\omega+c-1},
\end{align*}
if and only if \[  ~~~a\leq \frac{(2\omega-1)c+2\omega+2}{2(c-1)}. \]
Since $\frac{(2\omega-1)c+2\omega+2}{2(c-1)}> \omega-1$, it follows that Inequality \eqref{15} always holds for $a\leq \omega-1$. So, we assume that $a\geq \omega$.\\
\indent Now, consider the polynomial
\begin{align*}
f(k)=k^2-(2\mu_1-1)k+(2c-1)\omega+2n-\omega^2,~~~\omega+1\leq k \leq \omega+c-1.
\end{align*}
Clearly, $f(k)$ is decreasing for all $k\in \left [\omega+1,\mu_1-\frac{1}{2}\right ]$ and increasing for all $k\in \left [\mu_1-\frac{1}{2},\omega+c-1\right ]$. Therefore, if $f\left (\mu_1-\frac{1}{2}\right )\geq 0$, then Inequality \eqref{15} always holds. We have  $f\left (\mu_1-\frac{1}{2}\right )\geq 0$ if
\begin{align*}
\left (\mu_1-\frac{1}{2}\right )^2-(2\mu_1-1)\left (\mu_1-\frac{1}{2}\right )+(2c-1)\omega+2n-\omega^2\geq 0,
\end{align*}
if and only if \[ \mu_1\leq \frac{1}{2}+\sqrt{2n-\omega^2+(2c-1)\omega}. \]
For this $\mu_1$, Inequality \eqref{15} always holds. Since, $\mu_1\leq a+1=\frac{n}{\omega}$, we have 
\[ \frac{n}{\omega}\leq \frac{1}{2}+\sqrt{2n-\omega^2+(2c-1)\omega}, \]
which implies that
\[ 2a+1\leq 2\sqrt{2a\omega-\omega^2+(2c+1)\omega}, \]
that is, \[ 4a^2-(8\omega-4)a+(4\omega^2-(8c+4)\omega+1)\leq 0, \]
which further gives 
\begin{align*}
\omega-\frac{1}{2}-\sqrt{(2c-1)\omega}\leq a \leq \omega-\frac{1}{2}+\sqrt{(2c-1)\omega}.
\end{align*}
This shows that Inequality \eqref{15} holds for all $a$, with \[a\in \left [\omega-\frac{1}{2}-\sqrt{(2c-1)\omega},\omega-\frac{1}{2}+\sqrt{(2c-1)\omega}\right ].\] Since this inequality  already holds for $a\leq \omega-1$, it follows that Inequality \eqref{15} holds for all $a$, with $a\in \left [1,\omega-\frac{1}{2}+\sqrt{(2c-1)\omega}\right ]$. Thus, it follows that Brouwer's conjecture holds for all $\omega+1\leq k\leq \omega+c-1$, provided $a\leq \omega-\frac{1}{2}+\sqrt{(2c-1)\omega}$.\qed 

\indent Now, we consider the graphs, having $K_{s,s}$, $s\geq 2$ as the
maximal complete bipartite subgraph. In this direction, we have the following.

\begin{theorem}\label{Theorem 3.10.}
Let $G$ be a connected graph of order $n\geq 4$ and size $m$ and
let $K_{s,s}$, $s\geq 2$ be the maximal complete bipartite subgraph of graph $G$.
If $H=G\setminus K_{s,s}$ is a graph having $r$ non-trivial components $C_1,C_2,\ldots, C_r$, each of which is a c-cyclic graph and $p\geq 0$ trivial components, then for $s\geq \frac{5+\sqrt{8r(c-1)+34}}{2}$, Brouwer's conjecture holds for all $k$; and for $s< \frac{5+\sqrt{8r(c-1)+34}}{2}$, Brouwer's conjecture holds for all $k\in [x_1,n]$ and for all $k\in [1,y_1]$
, where $x_1=\frac{2s+3+\sqrt{20s-4s^2+8r(c-1)+9}}{2}$ and
$y_1=\frac{2s+3-\sqrt{20s-4s^2+8r(c-1)+9}}{2}$.
\end{theorem}
\noindent{\bf Proof.}
Let $G$ be a connected graph of order $n$ and size $m$. If
$H=G\setminus K_{s,s}$  is a graph having $r$ non-trivial components $C_1,C_2,\ldots, C_r$, each of which is a $c$-cyclic graph and $p\geq 0$ trivial components, then $m=s^2+n-p+(c-1)r$. For $ 1\leq k\leq n$,
from Corollary \ref{Corollary 2.6.}, we have
\begin{align*}
S_k(G)&\leq s+k(s+2)+n-p+2r(c-1)-\alpha<s+k(s+2)+n-p+2r(c-1)\\&
\leq m+\frac{k(k+1)}{2}=s^2+n-p+r(c-1)+\frac{k(k+1)}{2},
\end{align*} if
\begin{align}\label{16}
k^2-(2s+3)k-(2s+2r(c-1)-2s^2)\geq 0.
 \end{align}
\indent Now, consider the polynomial
\begin{align*}
f(k)=k^2-(2s+3)k-(2s+2r(c-1)-2s^2).
\end{align*}
The discriminant of this polynomial is $d=20s-4s^2+8r(c-1)+9$. We have $d\leq 0$ if  $20s-4s^2+8r(c-1)+9\leq 0$, which gives $s\geq \frac{5+\sqrt{8r(c-1)+34}}{2}$. This shows that for $s\geq \frac{5+\sqrt{8r(c-1)+34}}{2}$ the Inequality (\ref{16}) and so the Brouwer's conjecture  always  holds. For $s< \frac{5+\sqrt{8r(c-1)+34}}{2}$ the roots of the polynomial $f(k)$ are $x_1=\frac{2s+3+\sqrt{20s-4s^2+8r(c-1)+9}}{2}$ and
$y_1=\frac{2s+3-\sqrt{20s-4s^2+8r(c-1)+9}}{2}$, which implies that
$f(k)\geq 0$ for all $k\geq x_1$ and $f(k)\geq 0$ for all $k\leq y_1$.  So, for $s< \frac{5+\sqrt{8r(c-1)+34}}{2}$, the Inequality \eqref{16} and Brouwer's conjecture  holds for all $k\in [x_1,n]$ and for all $k\in [1,y_1]$.\qed 

\indent The following result explores some interesting families of graphs given by Theorem \ref{Theorem 3.10.} for which Brouwer's conjecture holds. We will use the fact that Brouwer's conjecture is always true for all $ k\leq 2 $.

\begin{corollary}\label{corollary 3.11.}
Let $G$ be a connected graph of order $n\geq 4$ and size $m$ having clique number
$\omega\geq 2$. Let $H=G\setminus K_{s,s},~s\geq 2$ be a graph having $r$ non-trivial components $C_1,C_2,\ldots, C_r$, each of which is a $c$-cyclic graph.\\
{\bf (i)} If $c=0$, that is, each $C_i$ is a tree, then Brouwer's conjecture holds  for all $k$, if $s\geq 5$ and holds for all $k,~~k\notin [3,7]$, if $2\leq s\leq 4$.\\
{\bf (ii)} If $c=1$, that is, each  $C_i$ is a unicyclic graph, then Brouwer's conjecture holds for all $k$, if $s\geq 6$ and holds for all $k,~~k\notin [3,7]$, if $2\leq s\leq 5$.\\
{\bf (iii)} If $c=2$, that is, each  $C_i$ is a bicyclic graph, then  Brouwer's conjecture holds for all $k$, if $s\geq 10 $ and holds for all $k,~~k\notin [3,12]$, if $2\leq s\leq 9$.\\
{\bf (iv)} If $c=3$, that is, each  $C_i$ is a tricyclic graph, then Brouwer's conjecture holds for all $k$, if $s\geq 14 $ and holds for all $k,~~k\notin [3,17]$, if $2\leq s\leq 13$.\\
{\bf (v)} If $c=4$, that is, each  $C_i$ is a tetracyclic graph, then Brouwer's conjecture holds for all $k$, if $s\geq 18 $ and holds for all $k,~~k\notin [3,20]$, if $2\leq s\leq 17$.\\

\end{corollary}
\noindent{\bf Proof.}
{\bf (i)}. If $c=0$, then each of the $r$ non-trivial components of $H$ are trees and so  Brouwer's conjecture  holds for all $k$, if $s\geq 5$; and holds for all $k,~~k\ne 2,3,4,5,6,7$ if
$2\leq s\leq 4$ see Theorem 3.9 in  \cite{hpv}.\\
{\bf (ii)}. If $c=1$, then each of the $r$ non-trivial components of $H$ are unicyclic graphs and so
the discriminant of the polynomial given by the left hand side of (\ref{16}) becomes $d=20s-4s^2+9$. Clearly, for $s\geq 6$, the discriminant $d<0$ and therefore for such $s$, Brouwer's conjecture always  holds. For $s=5$, we have $x_1=\frac{13+\sqrt{9}}{2}=8$ and $y_1=\frac{13-\sqrt{9}}{2}=5,$ implying that Brouwer's conjecture holds for all $k$, $k\ne 6,7$. For $s=4$, we have $x_1=\frac{11+\sqrt{25}}{2}=8$ and $y_1=\frac{11-\sqrt{25}}{2}=3$, implying that Brouwer's conjecture holds for all $k$, $k\ne 4,5,6,7$. For $s=3$, we have
$x_1=\frac{9+\sqrt{33}}{2}=7.3722$ and $y_1=\frac{9-\sqrt{33}}{2}=1.6277$, implying that Brouwer's conjecture holds for all $k$, $k\ne 2,3,4,5,6,7$. For $s=2$, we have 
$x_1=\frac{9+\sqrt{33}}{2}=7.3722$ and $y_1=\frac{9-\sqrt{33}}{2}=1.6277$, implying that Brouwer's conjecture holds for all $k$, $k\ne 2,3,4,5,6,7$.\\
{\bf (iii),(iv),(v)}. These follow by proceeding similar to the above cases.\qed

We have the following observations for a connected graph $ G $ of order $ n\geq 4 $, size $ m $ and having $ K_{s,s} , s\geq 2$ as its maximal complete bipartite subgraph. Let $ H\setminus K_{s,s}, s\geq 2 $ be a graph having $ r $ non-trivial components $ C_{1}, C_{2},\dots, C_{r}, $ each of which is a $ c- $cyclic graph. If $ c\leq \frac{s}{2}, ~ r\leq \frac{s}{2} $, then Brouwer's conjecture holds for all $ k $, if $ s=7, s\geq 9 $. If $ s=8 $, then Brouwer's conjecture holds for all $ k $ for $ c\leq 4, r\leq 3 $; holds for all $ k, k\neq 9,10 $ for $ c=r=4 $. If $ s=6 $, then Brouwer's conjecture holds for all $ k $ for $ c\leq 2, r\leq 1 $: holds for all $ k $, for $ c=r=2 $; holds for all $ k, k\neq 7,8 $ for $ c=2, r=3 $; holds for all $ k, k\notin[6,9] $ for  $ c=3, r=2 $; and holds for all $ k, k\notin[5,10] $, for $ c=r=3 $. If $ s=5, $ then Brouwer's conjecture holds for all $ k, k\neq 6,7 $ for $ c\leq 1, r\leq 2 $; holds for all $ k, k\notin[5,8] $, for $ c=2, r\leq 2 $. If $ s=4 $, then Brouwer's conjecture holds for all $ k, k\notin[4,7] $, for $ c\leq 1, r\leq 2 $, holds for all $ k, k\notin[3,8] $, for $ c=2, r\leq 2 $. If $ s=3, $ then Brouwer's conjecture holds for all $ k, k\neq 3,4,5,6,7. $ If $ s=2 $, then Brouwer's conjecture holds for all $ k, k\notin[3,6] $.

\indent We require the following observation, which gives a sequence of odd perfect squares.

\begin{lemma}\label{corollary 3.12.}
If $b_t=b_{t-1}+8(t+1), t\geq 1$ with $b_0=0$, then $x_t=b_t+9$ is perfect square. In fact $x_t=(2t+3)^2$, for all $t\geq 1$. 
\end{lemma}
\noindent{\bf Proof.} To prove this, we use induction on $t$. We have $x_1=b_1+9=b_{0}+16+9=25=5^2=(2(1)+3)^2; x_2=b_2+9=b_{1}+8(2)+9=49=7^2=(2(2)+3)^2$. This shows that result is true for $t=1$ and $t=2$. Assume that the result is true for $t=u$. Then by induction hypothesis, we have $x_u=(2u+3)^2$. For $t=u+1$, we have $x_{u+1}=b_{u+1}+9=b_u+8(u+2)+9=b_u+9+8(u+2)=x_u+8(u+2)=(2u+3)^2+8(u+2)=(2(u+1)+3)^2$. This shows that the result is true for $t=u+1$. Thus  by  induction, the result is true for  all $t\geq 1$.\qed 

\indent The following observation shows the existence of  more families of graphs, given by Theorem \ref{Theorem 3.10.}, for which Brouwer's conjecture  holds.

\begin{corollary}\label{corollary 3.13.}
Let $G$ be a connected graph of order $n\geq 4$ and size $m$ having clique number
$\omega\geq 2$. Let $H=G\setminus K_{s,s},~s\geq 2$ be a graph having $r$ non-trivial components $C_1,C_2,\ldots, C_r$, each of which is a $c$-cyclic graph.\\
{\bf (i)} If $-8<20s-4s^2+8r(c-1)\leq 0$, then Brouwer's conjecture holds  for all but  $k\in [s+1,s+2]$.\\
{\bf (ii)} If $b_{t-1}<20s-4s^2+8r(c-1)\leq b_{t}$, then Brouwer's conjecture holds for all but $k\in [s-1,s+4]$, where $ b_{t}=b_{t-1}+8(t+1) $ with $ t\geq 1 $ and $ b_{0}=1. $\\

\end{corollary}
\noindent{\bf Proof.}
$ \textbf{(i)} $. If $-8<20s-4s^2+8r(c-1)=0$, then we have $s+2< x_{1}=\frac{2s+3+\sqrt{20s-4s^2+8r(c-1)+9}}{2}\leq s+3$ and $ s\leq y_{1}=\frac{2s+3-\sqrt{20s-4s^2+8r(c-1)+9}}{2}<s+1 $. Therefore the desired result follows from Theorem \ref{Theorem 3.10.}.\\
$ \textbf{(ii)} $. This follows from Lemma \ref{corollary 3.12.} and by using Theorem \ref{Theorem 3.10.}. \qed

We note that the inequality  $-8<20s-4s^2+8r(c-1)\leq 0$ in the unknowns $r, c$ have many solutions. Some of the solutions are $r\leq \frac{s}{2}, c\leq s-6$; $r\leq \frac{s-5}{2}, c\leq s+1$; $r\leq s-5, c\leq \frac{s+2}{2}$; $r\leq s, c\leq \frac{s-3}{2}$; $r\leq \frac{s}{4}, c\leq 2s-9$.

\begin{corollary}\label{corollary 3.14.}
Let $G$ be a connected graph of order $n\geq 4$ and size $m$ having clique number
$\omega\geq 2$. Let $H=G\setminus K_{s,s},~s\geq 2$ be a graph having $r$ non-trivial components $C_1,C_2,\ldots, C_r$, each of which is a $c$-cyclic graph. Let $d=20s-4s^2+8r(c-1)+9$ be the  discriminant of the polynomial given by the left hand side of (\ref{16}). \\
{\bf (i)} If $r\leq s$, then Brouwer's conjecture holds  for all but $k\in [s-c,s+c+3]$.\\
{\bf (ii)} If $c\leq s$, then Brouwer's conjecture holds  for all but $k\in [s-r,s+r+3]$.
\end{corollary}
\noindent{\bf Proof.} If $r\leq s$, then 
\begin{align*}
d=&20s-4s^2+8r(c-1)+9\leq 12s-4s^2+8sc+9\\
=&9+(2c+3)^2-(2s-2c-3)^2< (2c+4)^2.
\end{align*}
Therefore, we have $x_1=\frac{2s+3+\sqrt{d}}{2}<\frac{2s+3+2c+4}{2}=s+c+\frac{7}{2}$ and $y_1=\frac{2s+3-\sqrt{d}}{2}>\frac{2s+3-2c-4}{2}=s-c-\frac{1}{2}$. This shows that Brouwer's conjecture holds  for all  $k,~~k\notin [s-c,s+c+3]$. Similarly, we can prove the other part. \qed

\section{Brouwer's conjecture with restriction on the size of graph}\label{section chen}

Chen \cite{chen} verified that Brouwer's conjecture is true for graphs in which the size $ m $ is restricted.
\begin{theorem}\cite{chen} \label{Cor1}
Let $G$ be a connected graph with $n\geq 4$ vertices and $m$ edges having $p\geq 1$ pendent vertices.
\begin{itemize}
\item[{\bf (i)}] If $p<\frac{n}{2}$ and $m \geq \frac{(n-1)(3n-1)}{8}-\frac{(n-3)p}{2}$, then  Brouwer's conjecture holds for $k \in [1,\frac{n-1}{2}]$.
\item[{\bf (ii)}]  If $p>\frac{n}{2} $ and $m \geq \frac{(n-1)(3n-1)}{8}-\frac{(n-3)p}{2}$, then
Brouwer's conjecture holds for $k \in [\frac{n-1}{2},n]$.
\end{itemize}
\end{theorem}

Now, let the graph $G$ have $p$ vertices each of degree $r$. The following theorem verifies Brouwer's conjecture under certain restrictions on the size $ m $ of $ G $. In fact this is a generalization of Theorem \ref{Cor1}.

\begin{theorem} \label{Theorem 2}
Let $G$ be a connected graph with $n\geq 4$ vertices and $m$ edges having $p\geq 1$ vertices of degree $r\geq 1$.
\begin{itemize}
\item[{\bf (i)}] If $m\geq \frac{(2n-r-1)r}{2}$, then Brouwer's conjecture holds for $k \in [1,r]$.
\item[{\bf (ii)}] If $p<\frac{n}{2}$ and $m \geq \frac{(n-1)(3n-1)}{8}-\frac{(n-1-2r)p}{2}$, then Brouwer's conjecture holds for $k \in [r+1,\frac{n-1}{2}]$.
\item[{\bf (iii)}]  If $p>\frac{n}{2}$ and $m \geq \frac{(n-1)(3n-1)}{8}-\frac{(n-1-2r)p}{2}$, then
Brouwer's conjecture holds for $k \in [\frac{n-1}{2},n]$.
\end{itemize}
\end{theorem}
\noindent{\bf Proof.} Let $G$ be a connected graph with $n$ vertices and $m$ edges having $p\geq 1$ vertices of degree $r\geq 1$. By definition of conjugated degrees, we have $ d_1^*\leq n, d_2^*\leq n, \dots , d_r^*\leq n, d_i^*\leq n-p$, where $i=r+1,\dots,n$. We consider the cases $k\leq r$ and $k\geq r+1$ separately. For $1\leq k\leq r$, by  Grone-Merris-Bai Theorem, it follows that
\begin{align*}
S_k(G)\leq \sum\limits_{i=1}^{k}d_i^*(G)\leq kn\leq m + \frac{k(k+1)}{2},
\end{align*}
if $k^2-(2n-1)k+2m\geq 0$. This shows that Brouwer's conjecture holds for all \[ k\leq \frac{2n-1-\sqrt{(2n-1)^2-8m}}{2}~~~ \text{as}~~~  \frac{2n-1+\sqrt{(2n-1)^2-8m}}{2}\geq n-1. \] Thus, $\frac{2n-1-\sqrt{(2n-1)^2-8m}}{2}\geq r$ implies that $m\geq \frac{(2n-r-1)r}{2}$, completing the proof in this case.\\
\indent For $r+1\leq k\leq n-1$, by  Grone-Merris-Bai Theorem, it follows that
\begin{align*}
S_k(G)\leq\sum\limits_{i=1}^{k}d_i^*(G)\leq k(n-p)+pr\leq m + \frac{k(k+1)}{2},
\end{align*}
provided that
$k^2-k(2n-2p-1) +2(m-pr) \geq 0 $.
Consider the polynomial
$f(x)=k^2-(2n-2p-1)k +2(m-pr)$, for $k\in [r+1,n-1]$. The discriminant of this polynomial is $\bigtriangledown=(2n-2p-1)^2-8(m-pr)$. Since $G$ has $ p $ vertices of degree $ r $, we have $m\leq\binom{n-p}{2}+pr <\frac{(2n-2p-1)^2}{2}+pr$, which implies that $\bigtriangledown>0$. The roots of this polynomial are
\begin{align*}
x= \frac{(2n-2p-1)-\sqrt{\bigtriangledown}}{2},\quad
y= \frac{(2n-2p-1)+\sqrt{\bigtriangledown}}{2}.
\end{align*} This shows that
$f(x)\geq 0$ is true  for $k\in [r+1,x]$ and $k\in [y,n-1].$\\
\indent If $p<\frac{n}{2}$, then $m \geq \frac{(n-1)(3n-1)}{8}-\frac{(n-1-2r)p}{2}$ implies that $x=\frac{(2n-2p-1)-\sqrt{\bigtriangledown}}{2}\geq \frac{n-1}{2}$. Similarly, if $p>\frac{n}{2}$, then $m \geq \frac{(n-1)(3n-1)}{8}-\frac{(n-1-2r)p}{2}$ implies that $y=\frac{(2n-2p-1)+\sqrt{\bigtriangledown}}{2}\leq \frac{n-1}{2}$,
completing the proof.\qed 

Clearly, if we choose $ r=1 $ in Theorem \ref{Theorem 2}, we get Theorem \ref{Cor1}.

Since Brouwer's conjecture holds for all $k\leq 2$, therefore by choosing $ r=2 $ in Theorem \ref{Theorem 2}, we have the following observation.

\begin{corollary}\label{Cor2}
Let $G$ be a connected graph with $n\geq 4$ vertices, $m$ edges and having $p\geq 1$ vertices of degree $ 2 $.  \\
{\bf (i)} If $p<\frac{n}{2}$ and $m \geq \frac{(n-1)(3n-1)}{8}-\frac{(n-5)p}{2}$, then  Brouwer's conjecture holds for $k \in [3,\frac{n-1}{2}]$.\\
{\bf (ii)}  If $p>\frac{n}{2} $ and $m \geq \frac{(n-1)(3n-1)}{8}-\frac{(n-5)p}{2}$, then Brouwer's conjecture holds for $k \in [\frac{n-1}{2},n]$.
\end{corollary}
Further, we observe that Corollary  \ref{Cor2} is more general than Theorem \ref{Cor1}, in the sense that it holds for more families of graphs than those given by Theorem \ref{Cor1}. Also if a connected graph $G$ has no pendent vertices, then  Theorem \ref{Cor1} is not applicable, however Theorem \ref{Theorem 2} is always applicable.

If a connected graph $G$ has a vertex of degree $r=3$, we have the following observation from Theorem \ref{Theorem 2}
\begin{corollary} \label{cor5}
Let $G$ be a connected graph with $n\geq 4$ vertices and $m$ edges having $p\geq 1$ vertices of degree $3$.  \\
{\bf (i)} If $m\geq 3n-6$, then Brouwer's conjecture holds for $k \in [1,3]$. \\
{\bf (ii)} If $p<\frac{n}{2}$ and $m \geq \frac{(n-1)(3n-1)}{8}-\frac{(n-7)p}{2}$, then Brouwer's conjecture holds for $k \in [4,\frac{n-1}{2}]$.\\
{\bf (iii)}  If $p>\frac{n}{2}$ and $m \geq \frac{(n-1)(3n-1)}{8}-\frac{(n-7)p}{2}$, then Brouwer's conjecture holds for $k \in [\frac{n-1}{2},n]$.
\end{corollary}
From part {\textbf{(i)}} of Corollary \ref{cor5}, it is clear that if a graph $G$ has a vertex of degree $3$, then Brouwer's conjecture holds for $k=3$, provided that $m\ge 3n-6$.

Now, let $ G $ be a connected graph with order $ n\geq 4 $ and let $ G $ have $ p\geq 1 $ and $ q\geq 1 $ $ (p\neq q) $ vertices of degrees $ r $ and $ s~ (s>r\geq 1) $ respectively. The following theorem verifies Brouwer's conjecture for the graphs when the size $ m $ is restricted in terms of $ n,~p,$ and $ q. $

\begin{theorem}\label{Theorem 3}
Let $G$ be a connected graph with $n\geq 4$ vertices, $m$ edges and having  $p\geq 1$ and $q\geq 1(q\neq p)$ vertices of degrees $r$ and $s~(s> r\geq 1)$, respectively.\\
{\bf (i)} If $m\geq \frac{(2n-r-1)r}{2}$, then Brouwer's conjecture holds for $k \in [1,r]$.\\
{\bf(ii)} If $n>p+s+\frac{1}{2}$ and $m\geq \frac{s(2n-2p-s-1)}{2}+pr$; or $n<p+r+\frac{3}{2}$ and $m\geq \frac{(r+1)(2n-2p-r-2)}{2}+pr$, then Brouwer's conjecture holds for $k \in [r+1,s]$.\\
{\bf(iii)} If $p+q<\frac{n}{2} $ and $m \geq \frac{(n-1)(3n-1)}{8}-\frac{(n-2r-1)}{2}p-\frac{(n-2s-1)}{2}q$, then
Brouwer's conjecture holds for all $k$,  $s+1\leq k\leq(n-1)/2$.\\
{\bf (iv)} If $p+q>\frac{n}{2} $ and $m \geq \frac{(n-1)(3n-1)}{8}-\frac{(n-2r-1)}{2}p-\frac{(n-2s-1)}{2}q$, then
Brouwer's conjecture holds for all $k$, $(n-1)/2\leq k\leq n$.
\end{theorem}
\noindent{\bf Proof.} Consider a graph $G$ as given in the statement. Then by definition of conjugated degrees, we have $ d_1^*\leq n, d_2^*\leq n, \dots , d_r^*\leq n,d_{r+1}^*\leq n-p,\dots, d_s^*\leq n-p,d_i^*\leq n-p-q$, where $i=s+1,\dots,n$. We consider the following cases $k\leq r$, $r+1\leq k\leq s$ and $s+1\leq k\leq n-1$. For $1\leq k\leq r$, the result follows by part \textbf{(i)} of Theorem \ref{Theorem 2}.  For $r+1\leq k\leq s$, by  Grone-Merris-Bai Theorem, it follows that
\begin{align*}
S_k(G)\leq \sum\limits_{i=1}^{k}d_i^*(G)\leq nr+(k-r)(n-p)\leq m + \frac{k(k+1)}{2},
\end{align*}
if $k^2-(2n-2p-1)k+2m-2pr\geq 0$. Consider the polynomial $f(x)=k^2-(2n-2p-1)k +2(m-pr)$, for $k\in [r+1,s]$. The discriminant of this polynomial is $\bigtriangledown=(2n-2p-1)^2-8(m-pr)$. Since $G$ has $p$ vertices of degree $r$ and $q$ vertices of degree $s$, we have $m\leq\binom{n-p-q}{2} +pr+qs<\binom{n-p}{2} +pr <\frac{(2n-2p-1)^2}{8}+pr$, which implies that $\bigtriangledown>0$. The roots of this polynomial are
\begin{align*}
x= \frac{(2n-2p-1)-\sqrt{\bigtriangledown}}{2},\quad
y= \frac{(2n-2p-1)+\sqrt{\bigtriangledown}}{2}.
\end{align*} This shows that
$f(x)\geq 0$ is true  for $k\in [r+1,x]$ and $k\in [y,s].$ For $n>p+s+\frac{1}{2}$, we have $x= \frac{(2n-2p-1)-\sqrt{\bigtriangledown}}{2}\geq s$ implying that $m\geq \frac{s(2n-2p-s-1)}{2}+pr$. For $n<p+r+\frac{3}{2}$, we have $y= \frac{(2n-2p-1)+\sqrt{\bigtriangledown}}{2}\leq r+1$ implying that $m\geq \frac{(r+1)(2n-2p-r-2)}{2}+pr$. This completes the proof in this case. \\
\indent Lastly,  for $s+1\leq k\leq n-1$, by  Grone-Merris-Bai Theorem, it follows that
\begin{align*}
S_k(G)&\leq nr+ (s-r)(n-p)+(k-s)(n-p-q)\\&=k(n-p-q)+rp+sq
\leq m + \frac{k(k+1)}{2}
\end{align*}
 provided that
 \begin{equation}\label{2.6}
k^2-k(2n-2(p+q)-1) +2(m-rp-sq) \geq 0.
\end{equation}
To complete the proof, we consider the polynomial
\begin{align*}
f(k)= k^2-(2n-2(p+q)-1)k +2(m-rp-sq), \quad k\in [s+1,n-1].
\end{align*}  The discriminant of this polynomial is $\bigtriangledown=(2n-2(p+q)-1)^2-8(m-rp-sq)$. Since $G$ has respectively $p$ and $q$ vertices of degree $r$ and $s$, we have $m<\binom{n-(p+q)}{2}+rp+sq$. It is easy to see that  $\binom{n-(p+q)}{2}+rp+sq<\frac{(2n-2(p+q)-1)^2}{2}+rp+sq$, implying that $m<\frac{(2n-2(p+q)-1)^2}{2}+rp+sq$, which in turn implies that $\bigtriangledown > 0$. The roots of $f(x)=0$ are
\begin{align*}
x= \frac{2n-2(p+q)-1-\sqrt{\bigtriangledown}}{2},\quad y= \frac{2n-2(p+q)-1+\sqrt{\bigtriangledown}}{2}.
\end{align*}
This shows that \eqref{2.6} holds for $k\leq x$ and  $k\geq y.$
Consequently, if $p+q<\frac{n}{2}$, then \eqref{2.6} holds for $s+1\leq k\leq (n-1)/2$, when
\begin{align*}
\frac{2n-2(p+q)-1-\sqrt{(2n-2p-1)^2-8(m-rp-sq)}}{2}\geq \frac{n-1}{2},
\end{align*}  that is,
\begin{align*}
(n-2(p+q))^2\geq \sqrt{(2n-2p-1)^2-8(e(G)-rp-sq)},
\end{align*} which follows when $m \geq \frac{(n-1)(3n-1)}{8}-\frac{(n-2r-1)}{2}p-\frac{(n-2s-1)}{2}q.$\\
\indent If $p+q>\frac{n}{2}$, then \eqref{2.6} holds for $(n-1)/2\leq k\leq n$, when
\begin{align*}
\frac{2n-2(p+q)-1+\sqrt{(2n-2p-1)^2-8(m-rp-sq)}}{2}\leq \frac{n-1}{2},
\end{align*}
 that is,
 \begin{align*}
 2(p+q)-n\leq \sqrt{(2n-2p-1)^2-8(m-rp-sq)},
 \end{align*}
 which also follows when $m \geq \frac{(n-1)(3n-1)}{8}-\frac{(n-2r-1)}{2}p-\frac{(n-2s-1)}{2}q.$
This completes the proof.\qed 

In case $G$ has $p\geq 1$ pendent vertices and $q\geq 1$ vertices of degree two, we have the following consequence of Theorem \ref{Theorem 3}.

\begin{corollary}\label{cor6}
Let $G$ be a connected graph with $n\geq 4$ vertices and $m$ edges having  $p\geq 1$ pendent vertices and $q\geq 1$ vertices of degrees $2$.\\
{\bf(i)} If $p+q<\frac{n}{2} $ and $m \geq \frac{(n-1)(3n-1)}{8}-\frac{(n-3)}{2}p-\frac{(n-5)}{2}q$, then
Brouwer's conjecture holds for all $k$,  $3\leq k\leq(n-1)/2$.\\
{\bf (ii)} If $p+q>\frac{n}{2} $ and $m \geq \frac{(n-1)(3n-1)}{8}-\frac{(n-3)}{2}p-\frac{(n-5)}{2}q$, then
Brouwer's conjecture holds for all $k$, $(n-1)/2\leq k\leq n$.
\end{corollary}
\noindent{\bf Proof.} This follows from Theorem \ref{Theorem 3}, by choosing $r=1$ and $s=2$.\qed 

\begin{remark}
If $G$ has $p\geq 1$ pendent vertices and at least one vertex of degree two, then Corollary \ref{cor6} is an improvement of Theorem \ref{Cor1}.
\end{remark}

Similar to Corollary \ref{cor6}, the following observation is immediate from Theorem \ref{Theorem 3}.

\begin{corollary}\label{cor7}
Let $G$ be a connected graph with $n\geq 7$ vertices and $m$ edges having  $p\geq 1$ pendent vertices and $q\geq 1$ vertices of degrees $3$.\\
{\bf(i)} If $n>p+\frac{7}{2}$ and $m\geq 3n-2p-6$; or $n<p+\frac{3}{2}$ and $m\geq 2n-p-3$, then  Brouwer's conjecture holds for $k=3$.\\
{\bf(ii)} If $p+q<\frac{n}{2} $ and $m \geq \frac{(n-1)(3n-1)}{8}-\frac{(n-3)}{2}p-\frac{(n-7)}{2}q$, then
Brouwer's conjecture holds for all $k$,  $4\leq k\leq(n-1)/2$.\\
{\bf (iii)} If $p+q>\frac{n}{2} $ and $m \geq \frac{(n-1)(3n-1)}{8}-\frac{(n-3)}{2}p-\frac{(n-7)}{2}q$, then
Brouwer's conjecture holds for all $k$, $(n-1)/2\leq k\leq n$.
\end{corollary}
\noindent{\bf Proof.} This follows from Theorem \ref{Theorem 3}, by taking $r=1$ and $s=3$.\qed 

\indent We note that part {\bf (i)} of Corollary \ref{cor7} imposes conditions on the number of edges and the number of vertices in terms of the number of pendent vertices of graph $G$ for Brouwer's conjecture to hold for $k=3$. This information will be helpful for further investigations, as one can investigate the case $k=3$, just as the case $k=2$ has been  discussed for any graph $G$. In fact, part {\bf (i)} of Corollary \ref{cor7} guarantees that for a graph $G$ with $n\geq p+4$ vertices Brouwer's conjecture holds for $k=3$, provided that $m\geq 3n-2p-6$. For $n\le p+2$ vertices Brouwer's conjecture holds for $k=3$, provided that $m\geq 2n-p-3$. If in particular $p=\frac{n}{2}$, then  Brouwer's conjecture  holds for $k=3$, provided that $m\geq 2n-6$.

The following upper bound for $S_k(G)$, in terms of the order $n$, the size $m$, the maximum degree $\Delta$ and the number of pendent vertices $p$ can be found in \cite{rt}:
\begin{align}\label{tre}
S_k(G)\le 2m-n+3k-\Delta+p+1.
\end{align}
With $k=3$, upper bound (\ref{tre}) implies that $S_3(G)\le 2m-n+10-\Delta+p\le m+6$, provided that $m\le n+\Delta-p-4$. Thus we have the following observation.
\begin{corollary}
Let $G$ be a connected graph of order $n\ge 7$ having size $m$ and maximum degree $\Delta$. Let $p$ be the number of pendent vertices of $G$. If $3\le p\le n-4$, then Brouwer's conjecture holds for $k=3$, provided that
\begin{align*}
 m\ge n+\Delta-p-4\quad  or \quad m\ge 3n-2p-6.
\end{align*}
\end{corollary}

For regular graphs, it is well known that Brouwer's conjecture always holds \cite{Mayank}. For biregular graphs (a graph in which degree of vertices is either $r$ or $s$ is said to be an $(r,s)$-regular graph or a biregular graph) no such result can be found in the literature. However, we have the following observation for a connected $(r,s)$-regular graph $G$.

\begin{corollary}\label{cor4}
Let $G$ be a connected $(r,s)$-regular graph  with $n\geq 7$ vertices and $m$ edges having  $p\geq 1$ vertices of degree $r$ and $q\geq 1$ vertices of degree $s, ~s>r$, with $p+q=n$.\\
{\bf(i)} If $pr+2qs\leq qs+r(r+1)$, then  Brouwer's conjecture holds for $k\in [1,r]$.\\
{\bf(ii)} If $2qs+pr\leq s(s+1)$, then Brouwer's conjecture holds for all $k$,  $ k\in [r+1,s]$.\\
{\bf (iii)} If $(p+q)^{2}\geq 4(pr+qs)+1$, then Brouwer's conjecture holds for all $k\in [(n-1)/2,n]$.
\end{corollary}
\noindent{\bf Proof.}  This follows from Theorem \ref{Theorem 3} by using the fact that $n=p+q$ and $m=\frac{pr+qs}{2}$.\qed 

To see the strength of Corollary \ref{cor4}, we consider some examples. For a $(2,3)$-regular graph $G$ of order $n\geq 7$, part {\bf (iii)} of Corollary \ref{cor4} implies that $(p+q)^2\geq 4(2p+3s)+1$, that is,
\begin{align}\label{a}
p(p-6)+q^2-2p-1+2q(p-6)\geq 0.
\end{align}
It is easy to see that (\ref{a}) is  true for $p\geq 6$ and $q\geq p$. For $p=4,5$, it can be seen that (\ref{a}) holds for $q\ge 6$. For $p=3$, it can be seen that (\ref{a}) holds for $q\ge 8$.  Thus, it follows if $p=4,5$ and $q\ge 6$, Brouwer's conjecture holds for all $k\in [(n-1)/2,n]$.

Let $S_{\omega}(H_1,H_2,\dots,H_{\omega})$, where $H_i$ is a graph of order $a_i,~~0\leq a_i<\omega,~~1\leq i\leq \omega$, be the family of connected graphs of order $n=\sum\limits_{i=1}^{\omega}a_i$ and size $m$ obtained by identifying a vertex of the graph $H_i$  at the $i^{th}$ vertex of the clique $K_{\omega}$. If $H_i=K_{1,a}$,  for $i=1,2,\dots,\omega-1$, $H_{\omega}=H^*$, where $H^*$ is the graph obtained by identifying a vertex of cycle $C_t$ with the root vertex of $K_{a-2,1}$ and the vertex of $H_i$ to be identified at $i^{th}$ vertex of $K_{\omega}$ is the root vertex, then any graph $G$ in the family $ S_{\omega}(H_1,H_2,\dots,H_{\omega})$,  can be obtained from a split graph (a graph $G$ whose vertex set $V(G)$ can be partitioned into two parts $V_1$ and $V_2$, such that the subgraph induced by $V_1$ is empty graph and the subgraph induced by $V_2$ is a clique) by fusing a vertex of a cycle $C_t$ at a vertex of the clique of degree $\omega+a-3$. For this family of graphs, we have the following observation.

\begin{theorem}\label{Theorem4}
Let $G$ belong to the family $S_{\omega}(H_1,H_2,\dots,H_{\omega})$, $H_i=K_{1,a}$. Then  Brouwer's conjecture holds for $k \in [1,\omega-0.5-u]$ and $k\in [\omega-0.5+u,n]$, where $u=\sqrt{2t-3.75}$.
\end{theorem}
\noindent{\bf Proof.} Let $G$ be a connected graph as in the hypothesis. Then $n=(a+1)(\omega-1)+a-1+t-1=\omega(a+1)+t-3$ and $m=\frac{\omega(\omega-1)}{2}+a\omega+t-2$. By definition of conjugated degrees, we have $ d_1^*= \omega(a+1)+t-3, d_2^*\leq\omega+t-1, d_3^*\leq\omega,\dots,d^{*}_{\omega+a-1}\leq\omega$ and $d^{*}_i=0$, for $i=\omega+a,\dots,n$. Since Brouwer's conjecture is always true for $k\leq 2$, we assume that $k\geq 3$. For $k\geq 3$, by  Grone-Merris-Bai Theorem, it follows that
\begin{align*}
S_k(G)&\leq \sum\limits_{i=1}^{k}d_i^*(G)\leq a\omega +2t-4+k\omega\\&
\leq m + \frac{k(k+1)}{2}=\frac{\omega(\omega-1)}{2}+a\omega+t-2+\frac{k(k+1)}{2},
\end{align*}
 provided
 \begin{equation}\label{2.9}
k^2-(2\omega-1)k+\omega(\omega-1)-2t+4 \geq 0.
 \end{equation}
Consider the function $ f(k)=k^2-(2\omega-1)k+\omega(\omega-1)-2t+4$, for $k\in[3,n-1]$. Since $t\geq 3$, it follows that the discriminant $\bigtriangledown=8t-15$ of the polynomial $f(k)$ is always positive. The roots of this polynomial are
\begin{align*}
\alpha= \frac{2\omega-1-\sqrt{8t-15}}{2} \quad \text{and} \quad \beta= \frac{2\omega-1+\sqrt{8t-15}}{2}.
\end{align*}
This shows that \ref{2.9} holds for all $k\in [3,\alpha]$ and for all $k\in [\beta,n-1]$. This completes the proof.\qed 

In particular, if $t=3,4,5$, we have the following consequence of Theorem \ref{Theorem4}.

\begin{corollary}
Consider $G$ in the family $S_{\omega}(H_1,H_2,\dots,H_{\omega})$. If $t=3$, then  Brouwer's conjecture holds for all $k$, except $k=\omega-1$ or $\omega$. If $t=4,5$, then  Brouwer's conjecture holds for all $k$, except $k=\omega-2$ or $\omega-1$ or $\omega$ or $\omega+1$.
\end{corollary}

\begin{remark}
In \cite{rt}, it is shown that the disjoint union of graphs which satisfy Brouwer's conjecture also satisfy Brouwer's conjecture. Therefore, any disjoint union of graphs that satisfy the hypothesis of theorems given in Section \ref{section chen}  also satisfy Brouwer's conjecture. Further, if $\overline{G}$ is the complement of the graph $G$, then it is shown in \cite{hmt} that Brouwer's conjecture   holds for $n-k-1$ for the complement $\overline{G}$, with $1\le k\le n-2$, if  it  holds for $k$ for the graph $G$. In the light of these statement and the fact $m(\overline{G})=\frac{n(n-1)}{2}-m(G)$, it is clear that the lower bounds given in Theorems \ref{Theorem 2} and \ref{Theorem 3} for the number of edges $m(G)$ of $G$ will give the upper bounds for the number of edges $m(\overline{G})$ of the graph $\overline{G}$.
\end{remark}

\newpage

\chapter{Laplacian energy conjecture of trees}

In this chapter, we investigate the Laplacian energy of trees. We verify the truth of the Laplacian energy conjecture completely for trees with diameter $ 4 $. Further, we settle the Laplacian energy conjecture for all trees having at most $\frac{9n}{25}-2$ non-pendent vertices. We also obtain some sufficient conditions for the Laplacian energy conjecture to hold for all trees of order $n$.

\section{Introduction}
Let  $ L(G) $ be the Laplacian matrix of the graph $G$ with its Laplacian eigenvalues  $\Big \{\mu_{1}(G),\mu_{2}(G),\dots, \mu_{n-1},\mu_{n}(G)\Big \} $.  Gutman and Zhou \cite{gz} defined the Laplacian energy $LE(G)$ of a graph $G$ as
\begin{equation*}
LE(G)=\sum\limits_{i=1}^{n}\left |\mu_{i}-\overline{d}\right |,
\end{equation*}
where $ \overline{d}=\frac{2m}{n} $ is the average degree of $ G $.  Using the fact that $\sum\limits_{1=i}^{n-1}\mu_{i}=2m$, from \cite{fhrt_1}, we have
\begin{equation}\label{laplacian energy}
LE(G)=2\left (\sum\limits_{i=1}^{\sigma}\mu_{i}-\sigma \overline{d}\right )=2\max_{1\leq k\leq n}\left( \sum\limits_{1=i}^{k}\mu_{i}-k\overline{d} \right ),
\end{equation}
where $\sigma$ is the number of Laplacian eigenvalues greater than or equal to the average degree $ \overline{d}$. The parameter $\sigma$ is an active component of the present research and some work mostly on trees can be found in the literature \cite{dmt,jacobs1,sin,zhou}. In fact, it is shown in \cite{radenkovic} that the Laplacian energy has remarkable chemical applications beyond the molecular orbital theory of conjugated molecules. Laplacian graph energy is a broad measure of graph complexity. Song et al. \cite{s} have introduced component-wise Laplacian graph energy, as a
complexity measure useful to filter image description hierarchies. For some recent works on Laplacian energy and related results, we refer to \cite{hb,hbp,ph} and the references therein.

One of the interesting problems is to determine the extremal value for Laplacian energy $LE(G)$ and to characterize graphs which attain such extremal values. This problem has been considered for various families of graphs, like trees, unicyclic graphs and the graphs attaining the maximum/minimum values are completely determined. However, this seems to be a hard problem, which
to the best of our knowledge, is still open.  In this direction, Radenkovi\'{c} and Gutman \cite{radenkovic} studied the correlation between the energy and the Laplacian energy of trees and they computed the energy and Laplacian energy of all trees up to $14$ vertices. They formulated the following conjecture.
\begin{conjecture}[Laplacian energy conjecture]\label{conjecture}
If $ T $ is a  tree of order $ n $, then
\begin{equation*}
LE(P_{n})\leq LE(T) \leq LE(S_{n}).
\end{equation*}
\end{conjecture}
Trevisan et al. \cite{trevisan}, proved that Conjecture \ref{conjecture} is true for all trees of diameter $3$, and further by direct computations they showed that Conjecture \ref{conjecture} is true for all trees up to $18$ vertices. Fritscher et al. \cite{fhrt_1} proved that the right inequality of Conjecture \ref{conjecture} is true for all trees of order $n$. For the left inequality, Chang et al. \cite{chang} verified the conjecture \ref{conjecture} for trees of diameter $4$ and $5$ with perfect matching. Rahman et al. \cite{rahman} considered some families of trees of diameter $4$ (not chosen in \cite{chang}) and verified the truth of left inequality of Conjecture \ref{conjecture} for trees of such families. But in general the left hand inequality of Conjecture \ref{conjecture} is still open. From (\ref{laplacian energy}), it is clear that in comparing the Laplacian energies of two graphs with same number of vertices and edges, we need an effective information about the sum of their largest Laplacian eigenvalues and about the number of Laplacian eigenvalues that are larger than the average degree.\\
\indent In this chapter our aim is to verify the truth of left hand inequality of Conjecture \ref{conjecture} for trees of diameter $4$, completely.

\section{Laplacian energy of a tree}

To prove our main results, we make use of the following algorithm  \cite{jacobs}, which is applied to
estimate the Laplacian eigenvalues of a given tree.

\noindent\textbf{Algorithm (I).}

The algorithm associates each vertex $ v $, a rational function $ a(v)=\frac{r}{s} $. Here $ r $ and $ s $ are the members of the polynomial ring $ \mathbb{Q}[\mu] .$ These are computed bottom-up starting with the leaves which are assigned $ \mu-1 $ (the trees can be rooted in an arbitrary way). Once all the children of $ v $ have been processed, $ v $ is assigned the function \[ a(v)=\mu-d_{v}-\sum_{c\in C} \frac{1}{a(c)}, \] where $ C $ is the set of its children and $ d_{v} $ is the degree of $ v $. When all the vertices have been processed, we compute the characteristic polynomial by taking the product of all functions $ a(v) $: \[ p(\mu)=\prod_{v\in V} a(v). \]

\noindent\textbf{Algorithm (II).} \cite{jacobs1} \\
\texttt{ Input: tree $T$, scalar $\alpha$}\\
\texttt{Output: diagonal matrix D congruent to $L(T)+\alpha I$}\\
Algorithm Diagonalize ($T,\alpha$)

initialize $a(v):=d(v)+\alpha,$ for all vertices $ v $

order vertices bottom up

\textbf{for} $ k=1\ \text{to}\ n $

\qquad \textbf{if} $ v_{k} $ is a leaf then continue

\qquad \textbf{else if} $ a(c)\neq 0 $ for all children $ c $ of $ v_{k} $ then

\qquad \qquad $ a(v_{k}):= a(v_{k})-\sum \frac{1}{a(c)},$ summing over all children of $ v_{k} $

\qquad\textbf{ else}

\qquad \qquad select one child $ v_{j} $ of $ v_{k} $ for which $ a(v_{j})=0 $

\qquad \qquad $ a(v_{k}):= -\frac{1}{2} $

\qquad\qquad $ a(v_{j}) := 2$

\qquad\qquad \textbf{if} $ v_{k} $ has a parent $ v_{l} $, remove the edge $ v_{k}v_{} $

end loop

This algorithm is useful to calculate the number of eigenvalues of the Laplacian matrix of $T$ lying in a
given interval and hence helps us to estimate $\sigma$ for a tree. It is worth to mention that the diagonal elements of the output matrix correspond precisely to the values $a(v)$ on each node $v$ of the tree.

The following observation due to Jacobs and Trevisan \cite{jacobs1} is helpful throughout the chapter.
\begin{lemma}\cite{jacobs1}\label{number1}
Let $ T $ be a tree and $ D $ be the diagonal matrix produced by the algorithm Diagonalize $ (T,-\alpha) $. Then the following assertions hold.\\
$ (a)  $ The number of positive entries in $ D $ is the number of the Laplacian eigenvalues of $ T $ that are greater than $ \alpha. $\\
$ (b)  $ The number of negative entries in $ D $ is the number of the Laplacian eigenvalues of $ T $ that are smaller than $ \alpha. $\\
$ (c)  $ If there are $ j $ zeros in $ D $, then $ \alpha $ is the Laplacian eigenvalues of $ T $ with multiplicity $ j $.
\end{lemma}

Using direct computation and the fact that
\begin{align*}
\frac{\pi}{n}\sum\limits_{j=1}^{\lfloor\frac{n}{2}\rfloor}\cos\frac{\pi j}{n}\le \int\limits_{0}^{\frac{\pi}{2}}\sin x dx=1,
\end{align*}
the following upper bound for the Laplacian energy of a path $P_n$ was established in \cite{trevisan}.
\begin{lemma}\label{pathenergy}\cite{trevisan}
If $ P_{n}  $ is a path on $ n $ vertices, then
\[LE(P_{n})\leq 2+\frac{4n}{\pi}.
\]
\end{lemma}
For non-increasing real sequences $(x)=(x_1,x_2,\dots,x_n)$ and $(y)=(y_1,y_2,\dots,y_n)$ of length $n$, we say that $(x)$ is majorized by $(y)$ or $(y)$ majorizes $(x)$, denoted by $(x)\preceq (y)$ if
\begin{align*}
&\sum\limits_{i=1}^{n}x_i =\sum\limits_{i=1}^{n}y_i~~~\text{and}~~~\sum\limits_{i=1}^{k}x_i \leq\sum\limits_{i=1}^{k}y_i, ~~\text{for  all}~~k=1,2,\dots,n-1.
\end{align*}

The following observation can be found in \cite{grone}.

\begin{lemma}\label{lem12}
Let $L(G)$ be the Laplacian matrix of $ G $. Then $(d_1,d_2,\dots,d_n)\preceq(\mu_1,\mu_2,\dots,\mu_n)$, that is;
\begin{align*}
\sum\limits_{i=1}^{k}\mu_i \ge 1+\sum\limits_{i=1}^{k}d_i ~~\text{for  all}~~k=1,2,\dots,n-1.
\end{align*}
\end{lemma}

\begin{theorem}\label{thm1}
Let $T$ be a tree of order $n\ge 4$ and let $P_n$ be the path graph on $n$ vertices. If $T$ has $s$ internal (non-pendent) vertices, then
\begin{align*}
LE(T)\ge LE(P_n),
\end{align*} provided that $\Big(\frac{\pi-2}{\pi}\Big)n\ge s+2-\frac{2s}{n}$.
\end{theorem}
\noindent{\bf Proof.} Let $\mu_1(T)\ge \mu_2(T)\ge \cdots\ge \mu_{n-1}(T)\ge \mu_n(T)=0$ be the Laplacian eigenvalues of $T$.  If $s\ge 1$ is the number of internal vertices (that is, the vertices having degree greater than 1) of $T$, then the number of pendent vertices $p$ of $T$ is given as $p=n-s$. As $1\le s\le n-2$, so using definition of Laplacian energy (\ref{laplacian energy}), we have
\begin{align}\label{u}
 LE(T)&=2\max_{1\leq k\leq n}\left( \sum\limits_{1=i}^{k}\mu_{i}(T)-k\overline{d} \right)\nonumber\\&
 \ge 2\Big( \sum\limits_{1=i}^{s}\mu_{i}(T)-s\overline{d} \Big).
\end{align}
If $d_1(T)\ge d_2(T)\ge \dots\ge d_n(T)$ are the vertex degrees of $T$, then by Lemma \ref{lem12}, we have
\begin{align*}
\sum\limits_{1=i}^{k}\mu_{i}(T)\ge 1+\sum\limits_{1=i}^{k}d_i,~~~\text{for~~ all}~~~ k=1,2,\dots,n-1.
\end{align*}
Using this in (\ref{u}), we get
\begin{align}\label{u1}
LE(T)\ge 2\Big(1+\sum\limits_{1=i}^{s}d_i-s\overline{d} \Big).
\end{align}
Since $d_1+d_2+\cdots+d_n=2m=2n-2$ and $T$ has $p=n-s$ pendent vertices, it follows that $d_1+d_2+\cdots+d_s+n-s=2n-2$, which implies that $\sum\limits_{1=i}^{s}d_i=n+s-2$. Using this in (\ref{u1}), we obtain
\begin{align*}\label{u2}
LE(T)\ge 2n+2s-2-2s\overline{d}\ge 2+\frac{4n}{\pi},
\end{align*}
provided that $\Big(\frac{\pi-2}{\pi}\Big)n\ge s+2-\frac{2s}{n}$. Now, using Lemma \ref{pathenergy}, the result follows.\qed

From Theorem \ref{thm1}, we have the following observation.

\begin{corollary}\label{cor1}
For a tree $T$ on $n$ vertices having vertex degrees $d_1\ge d_2\ge\dots\ge d_n$, and $k$, $1\le k\le n-1$, is any positive integer, then
\begin{align*}
LE(T)\ge 2\Big(1+\sum\limits_{1=i}^{k}d_i-k\overline{d} \Big).
\end{align*}
\end{corollary}

From Corollary \ref{cor1}, it follows that any information about the degrees of a tree $T$ can be used to obtain a lower bound for Laplacian energy $LE(T)$ of $T$, which in turn can be helpful for verification of Conjecture \ref{conjecture}.

\indent The following observation is also immediate from Theorem \ref{thm1}.

\begin{corollary}\label{cor2}
 Let $L(T)$ be the Laplacian matrix of tree $ T $ on  $n$ vertices. If $k$, $1\le k\le n-1$, is any positive integer, then
 \begin{align*}
 LE(T)\ge 2\Big(\sum\limits_{1=i}^{k}\mu_i-k\overline{d} \Big).
 \end{align*}
 \end{corollary}

From  Corollary \ref{cor2}, it is clear that any lower bound for the sum of $k$ largest Laplacian eigenvalues  $S_k(T)=\sum\limits_{1=i}^{k}\mu_i$ of $T$  can be used to obtain a lower bound for Laplacian energy $LE(T)$, which in turn can be helpful to verify Conjecture \ref{conjecture}.

Another observation showing the importance of Theorem \ref{thm1} is as follows.
 \begin{corollary}\label{cor3}
Let $T$ be a tree of order $n\ge 4$ having $s$ internal (non-pendent) vertices. Then following holds,
\begin{itemize}
\item[{\bf (i)}] If $s=1$, then Conjecture \ref{conjecture} holds for all $n\ge 9$;
\item[{\bf (ii)}] If $s=2$, then Conjecture \ref{conjecture} holds for all $n\ge 12$;
\item[{\bf (iii)}] If $s=3$, then Conjecture \ref{conjecture} holds for all $n\ge 14$;
\item[{\bf (iv)}] If $s=4$, then Conjecture \ref{conjecture} holds for all $n\ge 17$;
\item[{\bf (v)}] If $s=5$, then Conjecture \ref{conjecture} holds for all $n\ge 20$;
\item[{\bf (vi)}] If $s=6$, then Conjecture \ref{conjecture} holds for all $n\ge 23$;
\item[{\bf (vii)}] If $s=7$, then Conjecture \ref{conjecture} holds for all $n\ge 25$;
\item[{\bf (viii)}] If $s\le \frac{9n}{25}-2$, then Conjecture \ref{conjecture} holds for all $n$.
\end{itemize}
\end{corollary}

The following lemma \cite{groneme} is the interlacing property of the Laplacian eigenvalues of a graph and its spanning subgraph.
\begin{lemma}\label{lem1}
If $G^{\prime}=G+e$ is the graph obtained from $G$ by adding a new edge $e$, then the Laplacian
eigenvalues of $G$ interlace the Laplacian eigenvalues of $G^{\prime}$, that is,
\begin{align*}
\mu_1(G^{\prime})\geq\mu_1(G)\geq\mu_2(G^{\prime})\geq \mu_2(G)\geq\cdots\geq \mu_n(G^{\prime})\geq\mu_n(G)=0.
\end{align*}
\end{lemma}

The next lemma can be seen in \cite{sin}.

\begin{lemma}\label{lem3}
 The number of Laplacian eigenvalues less than the average degree $2-\frac{2}{n}$ of a tree $T$ of order $n$ is at least $\lfloor\frac{n}{2}\rfloor$.
 \end{lemma}

Now, we obtain a lower bound for the Laplacian energy of a tree $T$ in terms of the sum of $k_i$ largest Laplacian eigenvalues of $T_i$, where $T_i$, for $i=1,2$, are the components of $T$ obtained by deleting any non-pendent edge.
\begin{theorem}\label{thm2} 
Let $T$ be a tree of order $n\ge 8$ and let $e$ be a non-pendent edge of $T$. Let $T-e=T_1\cup T_2$ and let $\sigma$ be the number of Laplacian eigenvalues of $T-e$ which are greater than or equal to the average degree $\overline{d}(T-e)$. Then
\begin{align*}
LE(T)\ge 2S_{k_{1}}(T_1)+2S_{k_{2}}(T_2)-4\sigma+\frac{4\sigma}{n},
\end{align*} where $k_1,k_2$ are respectively, the number of Laplacian eigenvalues of $T_1$, $T_2$ which are greater than or equal to  $\overline{d}(T-e)$ with $k_1+k_2=\sigma$ and $S_k(T)$ is the sum of $k$ largest Laplacian eigenvalues of $T$.
\end{theorem}
\noindent{\textbf{ Proof.}} Let $T$  be a tree  with Lalacian matrix $ L(T) $ and the average degree $\overline{d}(T)=2-\frac{2}{n}$. Let $e$ be a non-pendent edge in $T$ and let $T_1$, $T_2$ be the components of $T-e$. Let $|V(T_i)|=n_i$, for $i=1,2$ and $\overline{d}(T-e)$ be the average degree of $T-e$. Then $n=n_1+n_2$ and $\overline{d}(T-e)=2-\frac{4}{n}$. If $\sigma$ is the number of Laplacian eigenvalues of $T-e=T_1\cup T_2$, which are greater than or equal to the average degree $\overline{d}(T-e)$, then it clear that $1\le \sigma\le n-2$.
Since $1\le \sigma\le n-2$, therefore  we have
\begin{equation}\label{k}
 LE(T)=2\max_{1\leq k\leq n}\left( \sum\limits_{i=1}^{k}\mu_{i}(T)-k\overline{d}(T) \right)
 \ge 2\Big( \sum\limits_{i=1}^{\sigma}\mu_{i}(T)-\sigma(2-\frac{2}{n}) \Big).
\end{equation}
Let $\mu_1(T-e)\ge \mu_2(T-e)\ge \cdots\ge \mu_{n-2}(T-e)\ge\mu_{n-1}(T-e)= \mu_n(T-e)=0 $ be the Laplacian eigenvalues of $T-e$. By Lemma \ref{lem1}, we have $\mu_i(T)\ge \mu_i(T-e)$ implying that $\sum\limits_{i=1}^{\sigma}\mu_{i}(T) \ge \sum\limits_{i=1}^{\sigma}\mu_{i}(T-e)$. Now, using (\ref{k}), we get
\begin{align}\label{v1}
 LE(T)\ge 2\sum\limits_{i=1}^{\sigma}\mu_{i}(T-e)-4\sigma+\frac{4\sigma}{n}.
\end{align}
Since Laplacian spectrum of $T_1\cup T_2$ is the union of Laplacian spectrum of $T_1$ and Laplacian spectrum of $T_2$, therefore $k_1+k_2=\sigma$, where $k_i\ge 1$ be the number of Laplacian eigenvalues of $T_i$ which are greater than or equal to  $\overline{d}(T-e)$. Thus, from (\ref{v1}), we have
\begin{align*}
 LE(T)\ge 2\Big(\sum\limits_{i=1}^{k_1}\mu_{i}(T_1)+\sum\limits_{i=1}^{k_2}\mu_{i}(T_2)\Big)-4\sigma+\frac{4\sigma}{n}.
\end{align*}
The result now follows.\qed

Let $\sigma_i$ be the number of Laplacian eigenvalues of $T_i$ which are greater than or equal to the average degree $\overline{d}(T_i)=2-\frac{2}{n_i}$. If $n_1\ge n_2$, the it is easy to see that $\overline{d}(T_1)\ge \overline{d}(T-e) $ and $\overline{d}(T_2)\le \overline{d}(T-e)$. Therefore, it follows that $k_1\ge \sigma_1$ and $k_2\le \sigma_2$. If $k_1=\sigma_1$ and $k_2=\sigma_2$, then we have the following observation.
\begin{corollary}\label{coru}
Let $T$ be a tree of order $n\ge 8$ and let $e$ be a non-pendent edge of $T$. Let $T_1$  and  $T_2$ be the components of $T-e$ and  let $\sigma$ be the number of Laplacian eigenvalues of $T-e$ which are greater than or equal to the average degree $\overline{d}(T-e)$. Let  $k_i$ and $\sigma_i$ be  respectively, the number of Laplacian eigenvalues of $T_i$ which are greater than or equal to  $\overline{d}(T-e)$ with $k_1+k_2=\sigma$ and the number of Laplacian eigenvalues of $T_i$ which are greater than or equal to  $\overline{d}(T_i)$. If $\sigma_1=k_1$ and $\sigma_2=k_2$, then Conjecture \ref{conjecture} holds for $T$, provided that $LE(T_i)\ge 2+\frac{4n_i}{\pi}$, for $i=1,2$.
\end{corollary}
\noindent{\bf Proof.}
Let $\sigma_i$ be the number of Laplacian eigenvalues of $T_i$ which are greater than or equal to the average degree $\overline{d}(T_i)=2-\frac{2}{n_i}$. Then, we have
\begin{align}\label{v2}
 LE(T_i)=2\sum\limits_{j=1}^{\sigma_i}\mu_{j}(T_i)-2\sigma_i \overline{d}(T_i).
\end{align}
Suppose that $LE(T_i)\ge 2+\frac{4n_i}{\pi}$, for $i=1,2$. Then, from (\ref{v2}), it follows that
\begin{align*}
2S_{\sigma_i}(T_i)\ge 2+\frac{4n_i}{\pi}+2\sigma_i \left (2-\frac{2}{n_i}\right ).
\end{align*}
Now, if $\sigma_i=k_i$, for $i=1,2$, then from Theorem \ref{thm2}, it follows that
\begin{align*}
LE(T)\ge& 2+\frac{4n_1}{\pi}+2\sigma_1 \left (2-\frac{2}{n_1}\right )+2+\frac{4n_2}{\pi}+2\sigma_2 \left (2-\frac{2}{n_2}\right )-4\sigma+\frac{4\sigma}{n}\\&
=\frac{4n}{\pi}+4+\frac{4\sigma}{n}-\left (\frac{4\sigma_1}{n_1}+\frac{4\sigma_2}{n_2}\right )\\&
\ge \frac{4n}{\pi}+2,
\end{align*}
provided that $1+\frac{2\sigma}{n}\ge \frac{2\sigma_1}{n_1}+\frac{2\sigma_2}{n_2}$. Note that we have  used $\sigma_1+\sigma_2=k_1+k_2=\sigma$.  Since $n=n_1+n_2$ and $\sigma=\sigma_1+\sigma_2$, therefore $1+\frac{2\sigma}{n}\ge \frac{2\sigma_1}{n_1}+\frac{2\sigma_2}{n_2}$ implies that $n_1n_2(n_1+n_2)\ge 2\sigma_1n^{2}_2+2\sigma_2n^{2}_1$, which in turn implies that $n_1^{2}(n_2-2\sigma_2)+n_2^{2}(n_1-2\sigma_1)\ge 0$. Now, applying Lemma \ref{lem3} to tree $T_i$,~it follows that this last inequality always holds. This completes the proof.\qed

\section{Trees of diameter at most 4}

In this section, we prove that Conjecture \ref{conjecture} is true for trees of diameter at most $4$.

A \emph{double broom} of diameter four is a tree consisting of a path on five vertices whose central vertex has degree two, while the other two non-leaf vertices may have arbitrary degree. It is usually denoted by $T(a,b)$, where $a$ is the number of pendent vertices on one of the non-leaf vertex and $b$ is the number of pendent vertices on another  non-leaf vertex, such that $a+b=n-3$. It is shown in Figure $ 3.1 $.

Let $ \mathcal{T}_{n}(d)$ be the family of trees each  of diameter $d $ and order $ n\geq 3$. In particular, $ \mathcal{T}_{n}(4)$ is the family of trees with order $n$ and diameter $4$. In $\mathcal{T}_{n}(4)$, there are more than $19$ possible subfamilies which need to be investigated in order to verify the truth of Conjecture \ref{conjecture} for $\mathcal{T}_{n}(4)$. From these subfamilies of trees, Conjecture \ref{conjecture} has already been verified for $5$ subfamilies; one in  \cite{chang} and four in \cite{rahman}.  Here, we consider all the trees in the family $\mathcal{T}_{n}(4)$ and verify the truth of Conjecture \ref{conjecture}.

In \cite{trevisan}, it is proved that Conjecture \ref{conjecture} is true for trees of diameter at most $3$. We give a simple proof of this theorem.

\begin{theorem}\label{thm23}
 Conjecture \ref{conjecture} is true for trees of diameter at most $3$.
\end{theorem}
\noindent{\bf Proof.} Let $T$ be  a tree of order $n\ge 3$ having diameter $2\le d\le 3$. If $d=2$, then $T$ is the star $S_{n}$ and so it follows by direct calculation that $LE(S_{n})\ge LE(P_n)$. If $d=3$, then $T$ is a double broom $T(a,b)$ of diameter $3$. Since the number of non-pendent vertices in $T(a,b)$ is $s=2$, it follows by part {\bf (ii)} of Corollary \ref{cor3} that the  Conjecture \ref{conjecture} holds for all $n\ge 12$. For $3\le n\le 11$, the truth of Conjecture \ref{conjecture} is verified by direct calculations. \qed

The following lemma gives a lower bound for the $i^{th}$-largest Laplacian eigenvalue in terms of the $i^{th}$-largest degree of graph $G$.

\begin{lemma}\cite{brouwerheamers}\label{eigenvalue lower bound}
Let $L(G) $ be the Laplacian eigenvalue of  $ G $. Then
\[  \mu_{i}\geq d_{i}-i+2, \qquad i=1,2,\dots,n. \]
\end{lemma}

\indent A tree $T$ such that the deletion of any non-pendent edge $e$ of $T$ splits it into a star and a tree that is not a star is said to be a \emph{Star-NonStar} tree (or SNS-tree for short) \cite{fhrt_1}. As observed in \cite{fhrt_1}, one can view an  SNS-tree as a tree with a root vertex $v_0$ with which three different types of branches
may be incident: (i) a single vertex, also called pendant vertex, which we call a branch of type $0$;
(ii) a tree of height one, called a branch of type $1$;
(iii) a tree of height two whose root has degree one, called a branch of type $2$.
Moreover, any such  tree satisfy two additional properties, namely, the combined number of branches
of type $1$ and $2$ in such a tree is at least two; when there are exactly two branches of type $1$ and $2$ incident with $v_0$ and at least one of them has type $1$, then $v_0$ is adjacent to pendant vertices (otherwise it is a double broom whose diameter is four or five). We note that an SNS-tree of diameter $4$ has no branch of type $2$ and has at
least two branches of the type $1$ (actually, in the absence of pendant vertices, there are at least
three branches of type $1$.

Since Conjecture \ref{conjecture} is true for all trees with order up to $n\le 18$, therefore, in the rest of the paper, we consider trees with order $n\ge 19$. The following theorem shows that Conjecture \ref{conjecture} holds true for all trees of the family $\mathcal{T}_n(4)$.

\begin{theorem}\label{thm42}
Conjecture \ref{conjecture} is true for all trees of diameter 4, that is, for all trees of the family $\mathcal{T}_n(4)$.
\end{theorem}
\noindent{\bf Proof.} Let $T\in \mathcal{T}_n(4)$ be a tree with order $n\ge 19$ and diameter $4$. Let $\mu_1(T)\ge \mu_2(T)\ge\dots\ge\mu_{n-1}(T)\ge \mu_n=0$ be the Laplacian eigenvalues of $T$ and let $\overline{d}(T)=2-\frac{2}{n}$ be its average vertex degree. By Lemma \ref{pathenergy}, to show that Conjecture \ref{conjecture} holds for $T$,  it suffices to show that the inequality
\begin{align}\label{a1}
LE(T)\ge \frac{4n}{\pi}+2,
\end{align} holds for $T$, for all  $n\ge 19$.
Since diameter of $T$ is $4$, it follows that either {\bf (i)} $T$ is a double broom of diameter $4$ (shown in Figure $3.1$) or {\bf (ii)} $T$ is an SNS-tree of diameter $4$ (shown in Figure $3.3$).\\
\indent If $T$ is  a double broom of diameter $4$, from Figure $3.1$, it is clear that $T$ has $s=3$ internal vertices. Therefore, using part {\bf (iii)} of Corollary \ref{cor3}, it follows that inequality $LE(T)\ge \frac{4n}{\pi}+ 2$ always holds, implying that the result is true in this case. So, assume that $T$ is an SNS-tree of order $n\ge 19$ having root vertex $v_0$ with $p\ge 0$ pendent vertices of level $0$, $r\ge 2$ vertices  of level $1$ such that each $v_i$ has $s_i$ pendent vertices attached, where at least two $s_i$ are non-zero. Clearly, in this case, order of $T$ is $n=p+r+\sum\limits_{i=1}^{r}s_i$. If $p=0$ and  $s_i=1$, for all $i=1,2,\dots,r$, then $T$ is the tree $ T(4;2a,2b)$, with $a+b=r$, as shown in Figure $3.1$ . Applying Algorithm $ {\bf (I)} $ to  $ T(4;2a,2b) $, we find that  the characteristic polynomial of $ T(4;2a,2b) $ is given by $\phi(T(4;a,b),x)= x(x^2-3x+1)^{a+b-1}(x^2-x(a+b+3)+2a+2b+1).$
Applying Algorithm $ {\bf (II)} $ with $\alpha=-1$ to the tree $ T(4;2a,2b) $ and using Lemma \ref{number1}, we see that $ a+b+2 $ eigenvalues are greater than one, while as $ a+b+1 $ eigenvalues are less than one. Since zeros of $ x^2-3x+1 $ and $x^2-x(a+b+3)+2a+2b+1  $ are respectively $ \frac{3\pm\sqrt{5}}{2} $ and $ \frac{1}{2}\left (3+a+b\pm\sqrt{(a+b)^2-2(a+b)+5}\right ) $, it follows that $ \sigma= a+b$. Therefore, using $2(a+b)=n-1$, from Equation \eqref{laplacian energy}, we get
\begin{align*}
LE(T(4;a,b))&= 2\left (\sum\limits_{i=1}^{\sigma}\mu_{i}(T(4;a,b))-\sigma \overline{d}(T(4;a,b))\right )\\&
= 2\left (\sum\limits_{i=1}^{a+b}\mu_{i}(T(4;a,b))-(a+b) \overline{d}(T(4;a,b))\right )\\&
= \sqrt{5}(a+b-1)+\sqrt{(a+b)^2-2(a+b)+5}+(a+b)\frac{4}{n}\\&
=\frac{\sqrt{5}}{2}\Big(n-3\Big)+2+\frac{1}{2}\sqrt{n^2-6n+25}-\frac{2}{n}\\&
\ge \frac{\sqrt{5}}{2}\Big(n-3\Big)+2+\frac{1}{2}(n-5)-\frac{2}{n}
>2+\frac{4n}{\pi},
\end{align*}
provided that
\begin{align}\label{x}
\left (\frac{\sqrt{5}+1}{2}-\frac{4}{\pi}\right )n> \left (\frac{3\sqrt{5}+5}{2}+\frac{2}{n}\right ).
\end{align}
As $n\ge 19$ implies that $\frac{2}{n}<0.11$, it follows that Inequality (\ref{x}) holds for all $n\ge 18$. This shows that Inequality (\ref{a1}) holds for $ T(4;2a,2b)$. If $p=0$, $s_1\ge 2$ and $s_i=1$, for $i=2,3,\dots,r$, then $T$ is the tree $T^{'}$ shown in Figure $3.2$. Applying Algorithm $ {\bf (I)} $ to  $ T^{'} $, its characteristic polynomial is given by $
 \phi(T^{'},x)= x(x-1)^{s_1-1}(x^2-3x+1)^{r-2} p(x)$, where $p(x)=(x^4-(r+s_1+5)x^3+(s_1r+4r+3s_1+8)x^2-(2s_1r+5r+2s_1+4)x+s_1+2r).$  By Lemma \ref{eigenvalue lower bound}, we have $ \mu_{1}(T^{'})\geq \max\{r,s_1+1\}+1>\overline{d}(T^{'})$ and $ \mu_{2}(T^{'})\geq \min\{r,s_1+1\}>\overline{d}(T^{'})$. Let $x_3$ and $x_4$, $x_3\ge x_4$, be the smallest and the second smallest zeros of $p(x)$. We have $p(0)=2r+s_1>0$, $p(1)=-s_1(r-1)<0$ and $p(2)=s_1>0$. Therefore, by the Intermediate value theorem, $x_4\in (0,1)$ and $x_3\in (1,2)$. So, either $\sigma=r$ or $r+1$ for the tree $T^{'}$. If $\sigma=r$, we have
 \begin{align*}
 LE(T^{'})&= 2\left (\sum\limits_{i=1}^{\sigma}\mu_{i}(T^{'})-\sigma \overline{d}(T^{'})\right )
 = 2\left (\sum\limits_{i=1}^{r}\mu_{i}(T^{'})-r \overline{d}(T^{'})\right )\\&
 \ge  \left (\sqrt{5}+1\right )r+2s_1-\left (2+2\sqrt{5}\right )+\frac{4r}{n}
 >2+\frac{4n}{\pi},
 \end{align*}
 provided that
 \begin{align}\label{y}
 \left (\sqrt{5}+1-\frac{8}{\pi}\right )r+\left (2-\frac{4}{\pi}\right  )s_1-\left (4+2\sqrt{5}\right )+\frac{4r}{n}> 0.
 \end{align} It is easy to see that Inequality (\ref{y}) holds for all $r\ge 5$ and $s_1\ge 7$. As $n=2r+s_1$, it follows that Inequality (\ref{y}) holds for all $n\ge 17$.  Therefore, Inequality (\ref{a1}) holds for $ T^{\prime}$ in this case. For $\sigma=r+1$, proceeding similarly as above, we see that Inequality (\ref{a1}) holds for $ T^{\prime}$ in this case as well. If $p=0$, $s_1,s_2\ge 2$ and $s_i=1$, for $i=3,\dots,r$, then $T$ is the tree $T^{\prime\prime}$ shown in Figure $3.2$. Applying Algorithm $ {\bf (I)} $ to  $ T^{\prime\prime} $, we observe that its characteristic polynomial is given by $ \phi(T^{''},x)= x(x-1)^{s_1+s_2-2}(x^2-3x+1)^{r-3} g(x)$, where
\begin{equation*}
\begin{split}
&g(x)=(x^6-(r+s_1+s_2+7)x^5+\alpha_1 x^4-\alpha_2 x^3+\alpha_3 x^2+\alpha_4 x+  s_1+s_2+2r-1),\\&
\alpha_1=rs_1+rs_2+s_1s_2+5s_1+6r+5s_2+19,\\&
\alpha_2=rs_1s_2+4rs_1+3s_1s_2+4rs_2+9s_1+9s_2+14r+24, \\&
\alpha_3=2rs_1s_2+5rs_1+3s_1s_2+5rs_2+7s_1+7s_2+16r+13, \\&
\alpha_4=2rs_1+2s_1s_2+2rs_2+3s_1+3s_2+9r+1.  
\end{split}
\end{equation*}
By Lemma \ref{eigenvalue lower bound}, we have
\begin{align*}
  &\mu_{1}(T^{\prime\prime})\geq \max\{r,s_1+1,s_2+1\}+1>\overline{d}(T^{\prime\prime}),\\& \mu_{2}(T^{\prime\prime})\geq \max\{\{r,s_1+1,s_2+1\}\setminus\theta\}>\overline{d}(T^{\prime\prime}),\\&
  \mu_3(T^{\prime\prime})\ge \min\{r,s_1+1,s_2+1\}-1>\overline{d}(T^{\prime\prime}),
\end{align*}
where $\theta=\max\{r,s_1+1,s_2+1\}$. Let $x_4, x_5$ and $x_6$, with $x_4\ge x_5\ge x_6$, respectively be the smallest, the second smallest and the third smallest zeros of $g(x)$. It is easy to see that $x_5,x_6< \overline{d}(T^{\prime\prime})$, implying that $\sigma=r$ or $r+1$ for the tree $T^{\prime\prime}$. If $\sigma=r$, then we have
\begin{align*}
LE(T^{\prime\prime})&= 2\left (\sum\limits_{i=1}^{\sigma}\mu_{i}(T^{\prime\prime})-\sigma \overline{d}(T^{\prime\prime})\right )
   = 2\left (\sum\limits_{i=1}^{r}\mu_{i}(T^{\prime\prime})-r \overline{d}(T^{\prime\prime})\right )\\&
   \ge  \left (\sqrt{5}+1\right )r+2s_1+2s_2-\left ( 5+3\sqrt{5}\right )+\frac{4r}{n}
   \\& >~2+\frac{4n}{\pi},
\end{align*}
provided that
   \begin{align}\label{z}
   \left (\sqrt{5}+1-\frac{8}{\pi}\right )r+\left (2-\frac{4}{\pi}\right )s_1+\left (2-\frac{4}{\pi}\right  )s_2-\left (7+3\sqrt{5}-\frac{4}{\pi}\right )+\frac{4r}{n}> 0.
\end{align}
Clearly, Inequality (\ref{z}) holds for all $r\ge 2$ and $s_1\ge 8$. Since $n=2r+s_1+s_2-1$, it follows that Inequality (\ref{z}) holds for all $n\ge 19$.  This shows that Inequality (\ref{a1}) holds for $ T^{\prime\prime}$. For $\sigma=r+1$, by proceeding similarly as above, Inequality (\ref{a1}) holds for $ T^{\prime\prime}$. This proves the result in this case.   \\
\indent Now, let $T$ be an SNS-tree other than $T(4;2a,2b), T^{\prime}$ and $T^{\prime\prime}$. We prove the result in this case by using induction on $r$. If $r=2$, then $T$ is a tree as shown in Figure $3.3$. Clearly, $T$ has $s=3$ internal vertices,  and so  using part {\bf (iii)} of Corollary \ref{cor3}, the result holds in this case. Assume that the result holds for all trees with $r=k$ vertices of level $1$. We show the result also holds for trees  with $r=k+1$ vertices of level $1$. Let $T$ be a tree with $r=k+1$ vertices (say) $v_1,v_2,\dots,v_k, v_{k+1}$ of level $1$. Let $e=v_0v_{k+1}$ be the edge joining the root $v_0$ with the vertex $v_{k+1}$. Delete edge $e$, and let $T_1$ and $T_2$ be the components of $T_e$. Assume that the order of $T_i$ is $n_i$ and the average vertex degree of $T_i$ is  $\overline{d}(T_i)=2-\frac{2}{n_i}$,  for $i=1,2$, with $n_1\ge n_2$. As $ r\ge 3$, so $T_1$ is an SNS-tree of order $n_1$ having $r=k+1$ vertices, namely $v_1,v_2,\dots,v_k$, of level $1$, while as $T_2$ is a star $K_{1,s_{k+1}}$. By induction hypothesis, Inequality (\ref{a1}) holds for both $T_1$ and $T_2$. Let $\sigma_i$ be the number of Laplacian eigenvalues of $T_i$ which are greater than or equal to the average degree $\overline{d}(T_i)=2-\frac{2}{n_i}$. Since $T_2$ is a star with at least two vertices, and for a star all non-zero Laplacian eigenvalues except the  spectral radius are equal to $1$, it follows that $\sigma_2=1$. To compute $\sigma_1$, we use Algorithm $ {\textbf{(II)}} $ with $\alpha=-2+\frac{2}{n_1}$ to tree $ T_{1}$. Let $u$ correspond to the pendent vertices in $T_1$. We have
\begin{align*}
a(u)=1-2+\frac{2}{n_1}=-\frac{n_1-2}{n_1}<0,
\end{align*}
as $n_1\ge 4$. Therefore, the diagonal entries in the resulting diagonal matrix corresponding to each pendent vertex is negative. For the vertices $v_i$, $1\le i\le k$, of level $1$, we have
 \begin{align*}
 a(v_i)=s_i+1-2+\frac{2}{n_1}-\frac{s_i}{a(u)}=2s_1i-1+\frac{2s_i}{n_1-2}+\frac{2}{n_1}>0,
 \end{align*}
implying that the diagonal entries in the resulting diagonal matrix corresponding to each vertex $v_i$ of level $1$ is positive. From Lemma \ref{number1}, it is clear that $\sigma_1=k$ or $k+1$, depending on whether $a(v_0)<0$ or $a(v_0)>0$. For the root vertex $v_0$, we have
\begin{align}\label{b}
  a(v_0)&=k-2+\frac{2}{n_1}-p(-\frac{n_1-2}{n_1})-\sum\limits_{i=1}^{k}\frac{1}{a(v_i)}\\&
  =k-2+\frac{2}{n_1}+\frac{p(n_1-2)}{n_1}-\sum\limits_{i=1}^{k}\frac{1}{a(v_i)}.
\end{align}
If $p\ge 1$, then clearly $a(v_0)>0$. By Lemma \ref{number1}, we see that $\sigma_1=k+1$. So, let $p=0$ in $T_1$. If at least three $s_i$, say $s_1,s_2,s_3$, are greater or equal $2$, then again it can be seen that $a(v_0)>0$, which by Lemma \ref{number1} implies that $\sigma_1=k+1$ in this case as well. So, assume that at most two $s_i$'s are greater or equal to $2$. Then $T\in \{T(4;2a,2b), T^{\prime},T^{\prime\prime}\}$, which is not the case.  Thus, it follows that for the tree $T_1$, we have $\sigma_1=k+1$, the number of non-pendent vertices. Let $k_1$ be the number of Laplacian eigenvalues of $T_1$ greater or equal to $\overline{d}(T_1\cup T_2)=2-\frac{4}{n}$. Since $\overline{d}(T_1)\ge \overline{d}(T_1\cup T_2)$, it follows that $k_1\ge \sigma_1$. We claim that $k_1=\sigma_1$. Since $\sigma_1=k+1$, it follows that  $k_1\ge k+1$.  Applying Algorithm $ {\bf (II)} $ with $\alpha=-2+\frac{4}{n}$ to the tree $ T_{1}$, we get
\begin{align*}
a(u)=1-2+\frac{4}{n}=\frac{-(n-4)}{n}<0,
\end{align*} for all $u$ pendent vertices $u$ of $T_1$. This implies that $k_1\le k+1$. Thus, we must have $k_1=k+1=\sigma_1$, proving the claim in this case. Further for the tree $T_2$, let $k_2$ be the number of Laplacian eigenvalues greater than or equal to $\overline{d}(T_1\cup T_2)=2-\frac{4}{n}$. Since $\sigma_2=1$ and $T_2$ has at
least one edge, it follows that $k_2=\sigma_2=1$. Thus, for the components  $T_1$ and $T_2$ of $T-e$, we have shown that Inequality (\ref{a1}) holds and they satisfy the property that $k_1=\sigma_1$, $k_2=\sigma_2$. Now, applying  Corollary \ref{coru}, it follows that Inequality (\ref{a1}) holds for $T$, also. Hence, we conclude with the help of induction that the result is true for all $r\ge 2$. This completes the proof.\qed

\begin{figure}\label{double broom}
\centering

\begin{tikzpicture}
\draw[fill=black] (0,0) circle (2pt);

\draw[fill=black] (-1,1) circle (2pt);\draw[fill=black] (-1,0.5) circle (2pt);\draw[fill=black] (-1,-1) circle (2pt);

\draw[fill=black] (1,0) circle (2pt);

\draw[fill=black] (2,0) circle (2pt);
\draw[fill=black] (3,1) circle (2pt);\draw[fill=black] (3,.5) circle (2pt);\draw[fill=black] (3,-1) circle (2pt);

\draw[thin] (0,0)--(-1,1); \draw[thin] (0,0)--(-1,0.5); \draw[thin] (0,0)--(-1,-1);
\draw[thin] (0,0)--(1,0)--(2,0);\draw[thin] (2,0)--(3,1);\draw[thin] (2,0)--(3,-1);\draw[thin] (2,0)--(3,.5);

\node at (-1,-.2) {$ a $};
\node at (3,-.2) {$ b $};

\end{tikzpicture}\qquad
\begin{tikzpicture}
\draw[fill=black] (0,0) circle (2pt);
\draw[fill=black] (-1,1) circle (2pt);\draw[fill=black] (-2,1) circle (2pt);\draw[fill=black] (-1,.5) circle (2pt);\draw[fill=black] (-2,.5) circle (2pt);\draw[fill=black] (-1,-1) circle (2pt);\draw[fill=black] (-2,-1) circle (2pt);

\draw[fill=black] (1,1) circle (2pt);\draw[fill=black] (2,1) circle (2pt);\draw[fill=black] (1,.5) circle (2pt);\draw[fill=black] (2,.5) circle (2pt);\draw[fill=black] (1,-1) circle (2pt);\draw[fill=black] (2,-1) circle (2pt);

\draw[thin] (0,0)--(-1,1)--(-2,1);\draw[thin] (0,0)--(-1,.5)--(-2,.5);\draw[thin] (0,0)--(-1,-1)--(-2,-1);

\draw[thin] (0,0)--(1,1)--(2,1);\draw[thin] (0,0)--(1,.5)--(2,.5);\draw[thin] (0,0)--(1,-1)--(2,-1);

\node at (-2,-.2) {$ a $};

\node at (2,-.2) {$ b $};

\end{tikzpicture}
\caption{Double broom of diameter $4$ and the tree $T(4;2a,2b)$}
\end{figure}
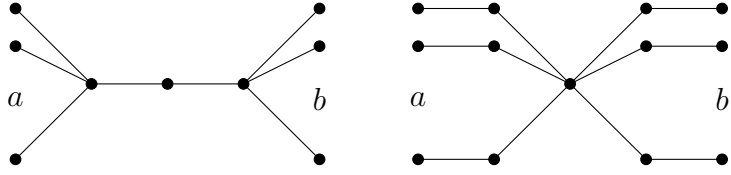

\begin{figure}\label{extended double broom}
\centering

\begin{tikzpicture}
\draw[fill=black] (0,0) circle (2pt);

\draw[fill=black] (-1,1) circle (2pt);\draw[fill=black] (-1,0.5) circle (2pt);\draw[fill=black] (-1,-1) circle (2pt);

\draw[fill=black] (1,0) circle (2pt);

\draw[fill=black] (2,1) circle (2pt);\draw[fill=black] (3,1) circle (2pt);\draw[fill=black] (2,.3) circle (2pt);\draw[fill=black] (3,.3) circle (2pt);\draw[fill=black] (2,-1) circle (2pt);\draw[fill=black] (3,-1) circle (2pt);

\draw[thin] (0,0)--(-1,1); \draw[thin] (0,0)--(-1,0.5); \draw[thin] (0,0)--(-1,-1);

\draw[thin] (0,0)--(1,0);

\draw[thin] (1,0)--(2,1)--(3,1);\draw[thin] (1,0)--(2,.3)--(3,.3);\draw[thin] (1,0)--(2,-1)--(3,-1);

\node at (-1,-.2) {$ s_{1} $};\node at (0,0.3) {$ v_{1} $};\node at (1,0.3) {$ v_{0} $};
\node at (2,1.3) {$ v_{2} $}; \node at (2,.6) {$ v_{3} $};\node at (2,-.7) {$ v_{r} $};


\end{tikzpicture}\qquad
\begin{tikzpicture}
\draw[fill=black] (0,0) circle (2pt);

\draw[fill=black] (-1,1) circle (2pt);\draw[fill=black] (-1,0.5) circle (2pt);\draw[fill=black] (-1,-1) circle (2pt);

\draw[fill=black] (1,0) circle (2pt);

\draw[fill=black] (0.5,-1) circle (2pt);\draw[fill=black] (-.5,-1.8) circle (2pt);\draw[fill=black] (-.1,-2) circle (2pt);\draw[fill=black] (1,-2) circle (2pt);

\draw[fill=black] (2,1) circle (2pt);\draw[fill=black] (3,1) circle (2pt);\draw[fill=black] (2,.3) circle (2pt);\draw[fill=black] (3,.3) circle (2pt);\draw[fill=black] (2,-1) circle (2pt);\draw[fill=black] (3,-1) circle (2pt);

\draw[thin] (0,0)--(-1,1); \draw[thin] (0,0)--(-1,0.5); \draw[thin] (0,0)--(-1,-1);

\draw[thin] (0,0)--(1,0);
\draw[thin] (1,0)--(0.5,-1)--(-.5,-1.8);\draw[thin](0.5,-1)--(-.1,-2);\draw[thin](0.5,-1)--(1,-2);

\draw[thin] (1,0)--(2,1)--(3,1);\draw[thin] (1,0)--(2,.3)--(3,.3);\draw[thin] (1,0)--(2,-1)--(3,-1);

\node at (-1,-.2) {$ s_{1} $};\node at (.5,-2) {$ s_{2} $};

\node at (0,0.3) {$ v_{1} $};\node at (1,0.3) {$ v_{0} $};
\node at (.8,-1) {$ v_{2} $};\node at (2,1.3) {$ v_{3} $};\node at (2,.6) {$ v_{4} $};\node at (2,-.7) {$ v_{r} $};


\end{tikzpicture}

\caption{Trees $ T^{'} $ and $ T^{''} $}
\end{figure}
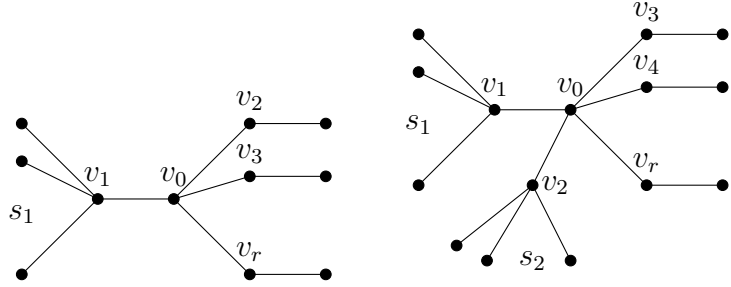
\begin{figure}\label{star non star}
\centering
\begin{tikzpicture}
\draw[fill=black] (-1,1) circle (2pt);\draw[fill=black] (-.3,1) circle (2pt);\draw[fill=black] (1,1) circle (2pt);\node at (0.4,1) {$...\ p\ ..$};

\draw[fill=black] (0,0) circle (2pt); \node at (0.4,0) {$v_{0}$};

\draw[fill=black] (-3,-2) circle (2pt);   \node at (-3.4,-2) {$v_{1}$};      \draw[fill=black] (-4,-3) circle (2pt);\draw[fill=black] (-3.5,-3) circle (2pt);\draw[fill=black] (-2,-3) circle (2pt);\node at (-2.8,-3) {$... s_{1} ...$};

\draw[fill=black] (0,-2) circle (2pt);   \node at (-0.4,-2) {$v_{2}$};     \draw[fill=black] (-1,-3) circle (2pt);\draw[fill=black] (-0.5,-3) circle (2pt);\draw[fill=black] (1,-3) circle (2pt);\node at (0.3,-3) {$... s_{2} ...$};

\node at (1.4,-2) {$\dots$};

\draw[fill=black] (3,-2) circle (2pt);   \node at (3.4,-2) {$v_{1}$};   \draw[fill=black] (2,-3) circle (2pt);\draw[fill=black] (2.5,-3) circle (2pt);\draw[fill=black] (4,-3) circle (2pt);\node at (3.3,-3) {$... s_{r} ...$};

\draw[thin] (0,0)--(-1,1);\draw[thin] (0,0)--(-.3,1);\draw[thin] (0,0)--(1,1);

\draw[thin] (0,0)--(-3,-2);\draw[thin] (-3,-2)--(-4,-3);\draw[thin] (-3,-2)--(-3.5,-3);\draw[thin] (-3,-2)--(-2,-3);

\draw[thin] (0,0)--(0,-2);\draw[thin] (0,-2)--(-1,-3);\draw[thin] (0,-2)--(-.5,-3);\draw[thin] (0,-2)--(1,-3);

\draw[thin] (0,0)--(3,-2);\draw[thin] (3,-2)--(2,-3);\draw[thin] (3,-2)--(2.5,-3);\draw[thin] (3,-2)--(4,-3);

\end{tikzpicture}
\caption{An SNS-tree of diameter $ 4 $}
\end{figure}
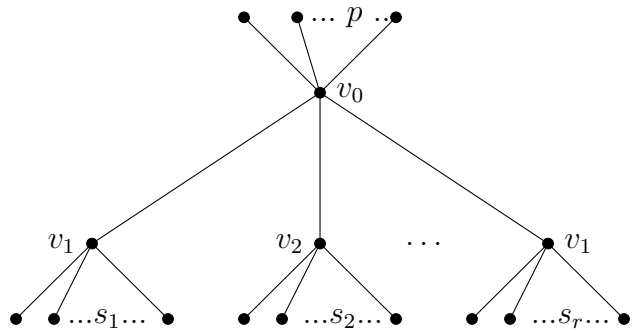

\section{Trees of any diameter}

Assume Conjecture \ref{conjecture} holds for two components obtained by deleting a non-pendent edge of a tree $T$. Also, for both the components, let the number of Laplacian eigenvalues greater than or equal to the average degree be equal to their number of non-pendent vertices. Then Conjecture \ref{conjecture} holds for $T$, as can be seen as below.

\begin{theorem}\label{thm51}
Let $T_1$ be  a tree of order $n_1$  having $r_1$ non-pendent vertices and let $T_2$ be a tree of diameter at most $3$ having order $n_2$ with $n_1\ge n_2\ge 6$. Let $\sigma_1$ be the number of Laplacian eigenvalues of $T_1$ greater than or equal to the average vertex degree $\overline{d}(T_1)=2-\frac{2}{n_1}$. Let $T=T_{1}\cup T_{2}\cup \{u,v\}$, where $ u\in T_{1} $ and $ v\in T_{2} $. If $\sigma_1=r_1$, then Inequality (\ref{a1}) (and hence Conjecture \ref{conjecture}) holds for $T$, provided that $LE(T_1)\ge 2+\frac{4n_1}{\pi}$.
\end{theorem}
\noindent{\bf Proof.} Let $\sigma$ be the number of Laplacian eigenvalues of $T_1\cup T_2$ which are greater than or equal to the average vertex degree $\overline{d}(T_1\cup T_2)=2-\frac{4}{n}$ and let $n_1+n_2=n$. Let $k_i$ be the number of Laplacian eigenvalues of $T_i$ which are greater than or equal to the average vertex degree $\overline{d}(T_1\cup T_2)=2-\frac{4}{n}$. Since $n_1\ge n_2$ implies that $\overline{d}(T_1)\ge \overline{d}(T_1\cup T_2)$ and $\overline{d}(T_2)\le \overline{d}(T_1\cup T_2)$, it follows that $k_1\ge \sigma_1$. Let $r_1$ be the number of non-pendent vertices in $T_1$. Assume that $\sigma_1=r_1$. Then $k_1\ge r_1$. We claim that $k_1=\sigma_1=r_1$.  Applying Algorithm $ {\bf (II)} $ with $\alpha=-2+\frac{4}{n}$ to the tree $ T_{1}$, we get
\begin{align*}
a(u)=1-2+\frac{4}{n}=\frac{-(n-4)}{n}<0,
\end{align*}
for all pendent vertices $u$ of $T_1$. This implies that $k_1\le r_1$. Combining, $k_1\ge r_1$ and  $k_1\le r_1$, we get $k_1=r_1=\sigma_1$, proving the claim.  Let $\sigma_2$ be the number of Laplacian eigenvalues of $T_2$ which are greater than or equal to average vertex degree $\overline{d}(T_2)=2-\frac{2}{n_2}$. Given that $T_2$ is a tree of diameter at most $3$, so therefore $T_2=S_{n_{2}}$ or $T_2=T(a,b)$, a double broom of diameter $3$. If $T_2=S_{n_{2}}$, then $\sigma_2=1$, the number of non-pendent vertices. Since  $T_2$ has at
least two edges and $\sigma_2=1$, it follows that $k_2=\sigma_2=1$. On the other hand, if $T_2=T(a,b)$, a double broom of diameter $3$, then as shown in \cite{trevisan}, we have $\sigma_2=2$, the number of non-pendent vertices. As  $T_2$ is of order $n_2\ge 6$ and $\sigma_2=2$, the number of non-pendent vertices, we must have $k_2=\sigma_2=2$. Thus, for trees  $T_1$ and $T_2$, we have  $k_1=\sigma_1$ and  $k_2=\sigma_2$. Since Inequality (\ref{a1}) always holds for $T_2$, by Theorem \ref{thm23}, and by using  Corollary \ref{coru}, it follows that Inequality (\ref{a1}) (and hence Conjecture \ref{conjecture})  holds for $T$, provided that $LE(T_1)\ge 2+\frac{4n_1}{\pi}$. This completes the proof.\qed

Assume $T_2$ is an SNS-tree of diameter $4$ other than the trees $T(4;2a,2b)$, or$ T^{\prime},$ or $T^{\prime\prime}$ defined in Section $3$. Then using Theorem \ref{thm42}, and proceeding similarly as in Theorem \ref{thm51}, we arrive at the following result.

\begin{theorem}\label{thm52}
Let $T_1$ be  a tree of order $n_1$  having $r_1$ non-pendent vertices. Let $T_2$ be an SNS-tree of order $n_2$ ($n_1\ge n_2\ge 6$)  and diameter $4$ other than the trees $T(4;2a,2b),$ or $ T^{\prime},$ or $T^{\prime\prime}$. Let $\sigma_1$ be the number of Laplacian eigenvalues of $T_1$ greater than or equal to average vertex degree $\overline{d}(T_1)=2-\frac{2}{n_1}$ and let $T=T_{1}\cup T_{2}\cup \{u,v\}$, where $ u\in T_{1} $ and $ v\in T_{2} $. If $\sigma_1=r_1$, then Inequality (\ref{a1}) (and hence Conjecture \ref{conjecture}) holds for $T$, provided that $LE(T_1)\ge 2+\frac{4n_1}{\pi}$.
\end{theorem}

Assume that Conjecture \ref{conjecture} holds for two components obtained by deleting a non-pendent edge of a tree $T$. Further, let one of the components $T_1$ satisfies $\mu_{\sigma_1+1}(T_1)- \overline{d}(T_1)<\frac{-2}{n}$ and the other component $T_2$ has the property that $\sigma_2$ is  equal to the number of non-pendent vertices. Then Conjecture \ref{conjecture} holds for $T$, as shown in the following theorem. Note that $\sigma_i$ is the  number of Laplacian eigenvalues of $T_i$ greater than or equal to average vertex degree $\overline{d}(T_i)$.

\begin{theorem}
Let $T_1$ be  a tree of order $n_1$  having Laplacian eigenvalues $\mu_1(T_1)\ge\dots\ge \mu_{n-1}(T_1),\mu_n(T_1)=0$. Let $T_2$ be a tree of order $n_2$ and diameter at most $3$ or an SNS-tree of order $n_2$ and diameter $4$ other than the trees $T(4;2a,2b),$ or $ T^{\prime},$ or $T^{\prime\prime}$  ($n_1\ge n_2\ge 6$). Let $\sigma_1$ be the number of Laplacian eigenvalues of $T_1$ greater than or equal to average vertex degree $\overline{d}(T_1)=2-\frac{2}{n_1}$. Let $T=T_{1}\cup T_{2}\cup \{u,v\}$, where $ u\in T_{1} $ and $ v\in T_{2} $, be the tree of order $n=n_1+n_2$. If $\mu_{\sigma_1+1}- \overline{d}(T_1)<\frac{-2}{n}$, then Inequality (\ref{a1}) (and hence Conjecture \ref{conjecture}) holds for $T$, provided that $LE(T_1)\ge 2+\frac{4n_1}{\pi}$.
\end{theorem}
\noindent{\bf Proof.} Let $\sigma$ be the number of Laplacian eigenvalues of $T_1\cup T_2$ which are greater than or equal to the average vertex degree $\overline{d}(T_1\cup T_2)=2-\frac{4}{n}$, and let $n_1+n_2=n$. Let $k_i$ be the number of Laplacian eigenvalues of $T_i$ which are greater than or equal to average vertex degree $\overline{d}(T_1\cup T_2)=2-\frac{4}{n}$. Since $n_1\ge n_2$ implies that $\overline{d}(T_1)\ge \overline{d}(T_1\cup T_2)$ and $\overline{d}(T_2)\le \overline{d}(T_1\cup T_2)$, it follows that $k_1\ge \sigma_1$. Since $\mu_{\sigma_1+1}- \overline{d}(T_1)<\frac{-2}{n}$, therefore
\begin{equation*}
\begin{split}
\mu_{\sigma_1+1}- \overline{d}(T_1\cup T_2)&=\mu_{\sigma_1+1}-2+\frac{4}{n}\\&=\mu_{\sigma_1+1}-2+\frac{2}{n_1}+\frac{4}{n}-\frac{2}{n_1}\\&
<\frac{-2}{n}+\frac{4}{n}-\frac{2}{n_1}\\&<0,
\end{split}
\end{equation*}
as $n>n_1$. Thus, it follows that $k_1\le \sigma_1$ and so we must have $\sigma_1=k_1$. Now, proceeding similar to Theorem \ref{thm51} the result follows.
\section{Conclusion}
In order to prove the conjecture in general, one must be aware of the distribution of Laplacian eigenvalues of trees around the average degree. The graph invariant $ \sigma $ investigated in \cite{dasmojallal2019}, plays a fundamental role in finding the lower bounds of  $ S_{k}(G) $, which in turn may help in proving the Laplacian energy conjecture. Thus to prove the Laplacian energy conjecture, we must study $ \sigma $ of the trees and use the  gained information in verifying the Laplacian energy conjecture. Some tree transformation as in \cite{sin} can also help in verifying the Laplacian energy conjecture.

\newpage

\chapter{Normalized Laplacian eigenvalues and eigenvalues of power graphs of finite cyclic groups}
In this chapter, we discuss the normalized Laplacian eigenvalues of the joined union of regular graphs in terms of the adjacency eigenvalues and the eigenvalues of the quotient matrix associated with graph $ G $. As a consequence of the joined union of graphs, we investigate the normalized Laplacian eigenvalues of the power graphs of finite cyclic group $ \mathbb{Z}_{n}. $
\section{Introduction}
The normalized Laplacian matrix of a graph $ G $, denoted by $ \mathcal{L}(G) $, whose rows and columns are indexed by the vertices of $ G $, is defined as 
\begin{equation*}
\mathcal{L}(G)=\begin{cases}
1 & \text{if} ~ v_{i}=v_{j} ~\text{and} ~ d_{v_{i}}\neq 0,\\
\dfrac{-1}{\sqrt{d_{v_{i}}d_{v_{j}}}} & \text{if}~ v_{i}\sim v_{j},\\
0 & \text{otherwise}.

\end{cases}
\end{equation*}
This matrix was introduced by Chung \cite{chung} to study the random walks of $ G $. This matrix is equivalently defined as 
\begin{align*}
\mathcal{L}(G)=& D(G)^{-\frac{1}{2}}L(G)D(G)^{-\frac{1}{2}}\\
=& D(G)^{-\frac{1}{2}}(D(G)-A(G))D(G)^{-\frac{1}{2}}\\
=& I-D(G)^{-\frac{1}{2}}A(G)D(G)^{-\frac{1}{2}},
\end{align*}
where $ D(G)^{-\frac{1}{2}} $ is the diagonal matrix whose $ i$-th diagonal entry is $ \frac{1}{\sqrt{d_{i}}} $. Clearly, $ \mathcal{L}(G) $ is real symmetric and positive semi-definite matrix so its eigenvalues are real, ordered as $ 0=\rho_{n}(\mathcal{L}) \leq \rho_{n-1}(\mathcal{L})\leq \dots \leq \rho_{1}(\mathcal{L})$ and are known as normalized Laplacian eigenvalues of $ G $. We see that $ D^{-\frac{1}{2}} \textbf{J}_{n\times 1} $ is an eigenvector of $ \mathcal{L}(G) $ with corresponding eigenvalue $ 0 $, where $ \textbf{J} $ is a vector with each entry $ 1 $. In \cite{chung}, it is shown that $ \sum\limits_{i=1}^{n} \rho_{i}(\mathcal{L})=n $, for each $ 1\leq i\leq n $, $~ \rho_{i}(\mathcal{L})\in [0,2] $ and $  \rho_{n}(\mathcal{L})= 2 $ if and only if a connected component of $ G $ is bipartite. Cavers et. al. \cite{caverslaa2010} studied the normalized Laplacian energy and characterized the extremal graphs. Two non-isomorphic graphs $ G_{1} $ and $ G_{2} $ are normalized Laplacian \emph{cospectral} if $ G_{1} $ and $ G_{2} $ have same normalized Laplacian eigenvalues. Butler \cite{butler2015} constructed normalized Laplacian cospectral graphs using twins and scaling technique. More work on $\mathcal{L}(G)$ can be seen in \cite{cavers, db, db1, sun} and references therein.
\section{Normalized Laplacian eigenvalues of the joined union of graphs}\label{section 4.2}

Consider the matrix
\begin{equation*}
M= \begin{pmatrix}
m_{1,1} & m_{1,2} & \cdots & m_{1,s} \\
m_{2,1} & m_{2,2} & \cdots & m_{2,s} \\
\vdots & \vdots & \ddots & \vdots \\
m_{s,1} & m_{s,2} & \cdots & m_{s,s} \\
\end{pmatrix}_{n\times n}, \end{equation*}
whose rows and columns are partitioned according to a partition $P=\{ P_{1}, P_{2},\dots , P_{m}\} $ of the set $X= \{1,2,\dots,n\}. $ The quotient matrix $ \mathcal{Q} $ of the matrix $ M $ is the $s \times s$ matrix whose entries are the average row sums of the blocks $ m_{i,j} $. The partition $ P $ is said to be \emph{ equitable} if each block $ m_{i,j} $ of $ M $ has constant row (and column) sum and in this case the matrix $ \mathcal{Q} $ is called as \emph{ equitable quotient matrix}. In general, the eigenvalues of $ \mathcal{Q}  $ interlace the eigenvalues of $ M $. In case the partition is equitable, we have following lemma.

\begin{lemma}\cite{BH1, cds}
If the partition $ P $ of $ X $ of matrix $ M $ is equitable, then each eigenvalue of $ \mathcal{Q} $ is an eigenvalue of $ M. $
\end{lemma}

If each of the $ G_{i} $ is $ r_{i} $ regular, then the following result gives the normalized Laplacian eigenvalues of $G[G_1,G_2,\dots,G_n]$ in terms of the adjacency eigenvalues of the graphs $G_1,G_2,\dots,G_n$ and the eigenvalues of the quotient matrix.

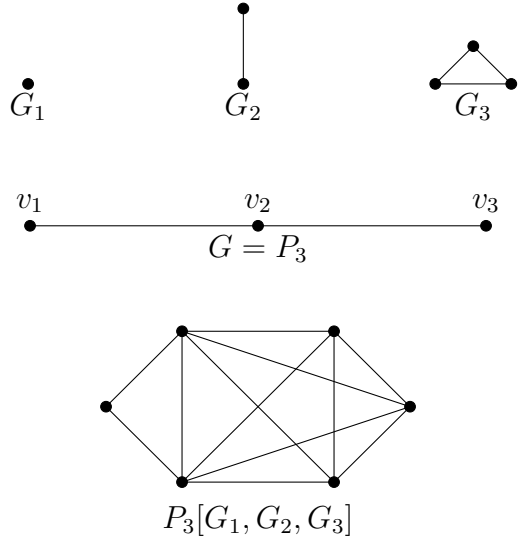
\begin{figure}\label{Joined union of graphs}
\centering

\begin{tikzpicture}
\draw[fill=black] (0,0) circle (2pt);

\node at (0,-.3) {$ G_{1} $};
\end{tikzpicture}\qquad\qquad\quad
\begin{tikzpicture}
\draw[fill=black] (0,0) circle (2pt);\draw[fill=black] (0,1) circle (2pt);
\draw[thin] (0,0)--(0,1);

\node at (0,-.3) {$ G_{2} $};
\end{tikzpicture}\qquad\qquad\quad
\begin{tikzpicture}
\draw[fill=black] (0,0) circle (2pt);\draw[fill=black] (1,0) circle (2pt);\draw[fill=black] (0.5,0.5) circle (2pt);

\draw[thin] (0,0)--(1,0); \draw[thin] (0,0)--(0.5,0.5);\draw[thin] (1,0)--(0.5,0.5);

\node at (.5,-.3) {$ G_{3} $};
\end{tikzpicture}
\\~\\
\begin{tikzpicture}
\draw[fill=black] (-1,0) circle (2pt);\draw[fill=black] (2,0) circle (2pt);\draw[fill=black] (5,0) circle (2pt);

\draw[thin] (-1,0)--(2,0); \draw[thin] (2,0)--(5,0);

\node at (-1,.3) {$ v_{1} $};\node at (2,.3) {$ v_{2} $};\node at (5,.3) {$ v_{3} $};

\node at (2,-.3) {$ G=P_{3} $};
\end{tikzpicture}
\\~\\
\begin{tikzpicture}
\draw[fill=black] (0,0) circle (2pt);\draw[fill=black] (1,1) circle (2pt); \draw[fill=black] (1,-1) circle (2pt);
\draw[fill=black] (4,0) circle (2pt);\draw[fill=black] (3,1) circle (2pt); \draw[fill=black] (3,-1) circle (2pt);

\draw[thin] (0,0)--(1,1); \draw[thin] (0,0)--(1,-1);\draw[thin] (1,1)--(1,-1);
\draw[thin] (4,0)--(3,1); \draw[thin] (4,0)--(3,-1);\draw[thin] (3,1)--(3,-1);

\draw[thin] (1,1)--(4,0);\draw[thin] (1,1)--(3,1);\draw[thin] (1,1)--(3,-1);
\draw[thin] (1,-1)--(4,0);\draw[thin] (1,-1)--(3,1);\draw[thin] (1,-1)--(3,-1);

\node at (2,-1.5) {$ P_{3}[G_{1},G_{2},G_{3}] $};
\end{tikzpicture}

\caption{Joined union of three graphs }
\end{figure}

\begin{theorem} \label{joined union}
Let $G$ be a graph of order $n$ and size $m$. Let $ G_{i}$ be $r_{i}$ regular graphs of order $ n_{i} $ having adjacency eigenvalues $\lambda_{i1}=r_{i}\geq \lambda_{i2}\geq\ldots\geq \lambda_{in_{i}}, $ where $ i=1,2, \ldots, n$.  Then the normalized Laplacian eigenvalues of the graph $ G[G_{1},\ldots, G_{n}] $ consists of the eigenvalues $ 1-\dfrac{1}{r_{i}+\alpha_{i}}\lambda_{ik}(G_{i})$, for $ i=1,\ldots,n $ and $ k=2,3,\ldots, n_{i} $, where $\alpha_i=\sum\limits_{v_j\in N_{G}(v_{i})}n_{i}$ is the sum of the orders of the graphs $G_j, j\ne i$ which correspond to the neighbours of vertex $ v_{i}\in G $. The remaining $n$ eigenvalues are the eigenvalues of matrix \eqref{Qmat of joined union}.

\end{theorem}
\noindent\textbf{ Proof.}
Let $V(G)=\{ v_{1}, \ldots, v_{n}\}$ be the vertex set of $ G $ and let $V(G_i)=\{ v_{i1}, \ldots, v_{in_i}\}$ be the vertex set of the graph $G_i$, for $i=1,2,\dots,n$. Let  $H=G[G_{1}, \ldots, G_{n}]$ be the joined union  of $ r_{i} $ regular graphs $G_i,$ for $ i=1,2,\dots,n $. It is clear that the order of $ H $ is $ N=\sum\limits_{i=1}^{n}n_{i} $. Since degree of each vertex $v_{ij}\in V(H)$ is degree inside $ G_{i} $ and  the sum of orders of $ G_{j}\text{'s}, j\ne i $, which corresponds to the neighbours of the vertex $ v_{i} $ in $ G $, where $1\le i\leq n$ and $1\leq j\leq n_i$, therefore, for  each $v_{ij}\in V(G_{i})$, we have
\begin{equation}\label{eq23}
d_{H}(v_{ij})=r_{i}+\sum\limits_{v_j\in N_{G}(v_i)}n_{j}=r_i+\alpha_i,
\end{equation} where $\displaystyle\alpha_i=\sum\limits_{v_j\in N_{G}(v_i)}n_{j}$. Under suitable labelling of the vertices in $ H $, the normalized Laplacian matrix of $H$ can be written as
\begin{equation*}
\mathcal{L}(H)=\begin{pmatrix}
g_1 & a(v_1,v_2)& \ldots &  a(v_1,v_n) \\
a(v_2,v_1) & g_2 & \ldots &  a(v_2,v_n)\\
\vdots &\vdots &\ddots &\vdots\\
 a(v_n,v_1) & a(v_n,v_2) &\ldots & g_n
\end{pmatrix},
\end{equation*}
where,  for $i=1,2,\ldots,n,$
\begin{equation*}
g_{i}=I_{n_{i}}-\dfrac{1}{r_{i}+\alpha_{i}}A(G_{i})
\end{equation*} and
\begin{equation*} a(v_{i}, v_{j}) = \begin{cases}
\dfrac{1}{\sqrt{(\alpha_{i}+r_{i})(\alpha_{j}+r_{j})}}J_{n_{i}\times n_{j}} ,& \text{if}~  v_{i}\sim v_{j}~  \text{in}~ G\\
\textbf{0}_{n_{i}\times n_{j}}, & \text{otherwise.}
\end{cases}
\end{equation*}
$ A(G_i)$ is the adjacency matrix of $ G_{i}$, $J_{n_i\times n_{j}}$ is the matrix having all entries $ 1 $, $ \textbf{0}_{n_{i}\times n_{j}} $ is the zero matrix of order $ n_{i}\times n_{j} $ and $ I_{n_i}$ is the identity matrix of order $n_i$. \\
\indent As each $ G_{i} $ is an $ r_{i} $ regular graph, so the all one vector $e_{n_i}=(\underbrace{1,1,\dots,1}_{n_{i}})^T$  is the eigenvector of the adjacency matrix $ A(G_{i}) $ corresponding to the eigenvalue $r_{i} $ and all other eigenvectors are orthogonal to $e_{n_i}.$ Let $ \lambda_{ik}$, $2\leq k\leq n_i$, be any eigenvalue of $ A(G_{i})$ with the corresponding eigenvector $X=(x_{i1},x_{i2},\dots,x_{in_i})^T$ satisfying $e_{n_i}^TX=0.$ Clearly, the column vector $X$ can be regarded as a function defined on $ V(G_i) $ assigning the vertex $ v_{ij} $ to $ x_{ij} $, that is, $ X(v_{ij})=x_{ij} $ for $ i=1,2,\ldots,n $ and $j=1,2,\dots,n_i$. Now, consider the vector $Y=(y_{1},y_{2},\dots,y_{n})^T$, where
\begin{equation*}
y_{j}=\left\{
\begin{array}{rl}
x_{ij}&~\text{if}~ v_{ij}\in V(G_i)\\
0&~\text{otherwise.}\\
\end{array}\right.
\end{equation*}
Since, $ e_{n_{i}}^{T}X=0 $ and coordinates of the  vector $Y$ corresponding to vertices in $\cup_{j\ne i}V_j$ of $H$ are zeros,  we have
\begin{align*}
\mathcal{L}(H)Y&=\begin{pmatrix}
0\\
\vdots\\
0\\
X-\dfrac{1}{r_{i}+\alpha_{i}}\lambda_{ik} X\\
0\\
\vdots\\
0
\end{pmatrix}\\
&=\left (1-\dfrac{1}{r_{i}+\alpha_{i}}\lambda_{ik}\right )Y.
\end{align*}
This shows that $Y$ is an eigenvector of $\mathcal{L}(H) $ corresponding to the eigenvalue 
\[ 1-\dfrac{1}{r_{i}+\alpha_{i}}\lambda_{ik}, \] for every eigenvalue $\lambda_{ik}$, $2\leq k\leq n_i$, of $A(G_i)$. In this way, we have obtained \[  \sum\limits_{i=1}^{n}n_{i}-n=N-n  \] eigenvalues. The remaining $ n $ normalized Laplacian eigenvalues of $ H $ are the eigenvalues of the equitable quotient matrix
\begin{equation}\label{Qmat of joined union}
\begin{pmatrix}
\dfrac{\alpha_{1}}{\alpha_1+r_1}& \dfrac{-n_2a_{12}}{\sqrt{(r_{1}+\alpha_{1})(r_{2}+\alpha_{2})}}&\ldots & \dfrac{-n_na_{1n}}{\sqrt{(r_{1}+\alpha_{1})(r_{n}+\alpha_{n})}}\\
\dfrac{-n_1a_{21}}{\sqrt{(r_{2}+\alpha_{2})(r_{1}+\alpha_{1})}}&\dfrac{\alpha_{2}}{\alpha_2+r_2}& \ldots & \dfrac{-n_na_{2n}}{\sqrt{(r_{2}+\alpha_{2})(r_{n}+\alpha_{n})}}\\
\vdots &\vdots  &\ddots &\vdots\\
\dfrac{-n_1a_{n1}}{\sqrt{(r_{n}+\alpha_{n})(r_{1}+\alpha_{1})}}& \dfrac{-n_2a_{n2}}{\sqrt{(r_{n}+\alpha_{n})(r_{2}+\alpha_{2})}}&\ldots &\dfrac{\alpha_{n}}{\alpha_n+r_n}
\end{pmatrix},
\end{equation} 
where, for $i\ne j$,
\begin{equation*}
a_{ij}=\begin{cases}
1, & v_i\sim v_j\\
0, & \text{otherwise}.

\end{cases}
\end{equation*} \qed

The next observation is a consequence of Theorem \ref{joined union} and gives the normalized Laplacian eigenvalues of the complete multipartite graph.

\begin{corollary}\label{p-partite}
The normalized Laplacian eigenvalues of the complete $p$-partite graph $K_{n_1,n_2,\dots,n_p}=K_{p}[\overline{K}_{n_1},\overline{K}_{n_2},\dots,\overline{K}_{n_p}]$ with $N=\displaystyle\sum_{i=1}^{p}n_{i}$ consists of the eigenvalue $1$ with multiplicity $N-p$ and the remaining $ p $ eigenvalues are given by the matrix
\begin{equation*}
\begin{pmatrix}
1&\dfrac{-n_{2}}{\sqrt{\alpha_{1}\alpha_{2}}}&\ldots & \dfrac{-n_{p}}{\sqrt{\alpha_{1}\alpha_{p}}}\\
\dfrac{-n_{1}}{\sqrt{\alpha_{2}\alpha_{1}}}&1& \ldots & \dfrac{-n_{p}}{\sqrt{\alpha_{2}\alpha_{p}}}\\
\vdots &\vdots  &\ddots &\vdots\\
\dfrac{-n_{1}}{\sqrt{\alpha_{p}\alpha_{2}}}& \dfrac{-n_{2}}{\sqrt{\alpha_{p}\alpha_{2}}}&\ldots &1
\end{pmatrix}.
\end{equation*}
\end{corollary}
\noindent\textbf{Proof.} This follows from Theorem \ref{joined union}, by taking $ G_{i}=\overline{K}_{i} $ and $ \lambda_{ik}(G_{i})=0 $ for each $ i $ and each $ k $. \qed

In particular, if partite sets are of equal size, say $ n_{1}=n_{2}=\dots= n_{p}=t $, then we have the following observation.
\begin{corollary}
Let $ G= K_{t,t,\dots,t} $ be a complete $ p$-partite graph of order $ N=pt .$ Then the normalized Laplacian eigenvalues of $ G $ consists of the eigenvalue $1 $ with multiplicity $ N-p $, the eigenvalue $\dfrac{p}{p-1}$ with multiplicity $p-1$ and the eigenvalue $ 0 $.
\end{corollary}
\noindent\textbf{Proof.} By Theorem \ref{joined union}, we have $ \alpha_{i}=t(p-1) $, for $ i=1,2,\dots,p $. Also, by Corollary \ref{p-partite}, we see that $ 1 $ is an eigenvalue with multiplicity $ pt-p $ and other eigenvalues are given by
\begin{equation*}
M_{p}=\begin{pmatrix}
1&\dfrac{-1}{p-1}&\ldots & \dfrac{-1}{p-1}\\
\dfrac{-1}{p-1}&1& \ldots & \dfrac{-1}{p-1}\\
\vdots &\vdots  &\ddots &\vdots\\
\dfrac{-1}{p-1}& \dfrac{-1}{p-1}&\ldots &1
\end{pmatrix}.
\end{equation*}
By simple calculations, we see that the normalized Laplacian eigenvalues of matrix $ M_{p} $ consists of the eigenvalue $ \dfrac{p}{p-1} $ with multiplicity $ p-1 $ and simple eigenvalue $ 0. $ \qed

Another consequence of Theorem  \ref{joined union}, gives the normalized Laplacian eigenvalues of the join of two regular graphs.
\begin{corollary}\label{join of two graphs}
Let $ G_{i}$ be an $ r_{i} $ regular graph of order $ n_{i} $ for $i=1,2$. Let $ \lambda_{ik}, 2\leq k\leq n_{i}, i=1,2 $ be the adjacency eigenvalues of $ G_{i} $. Then the normalized Laplacian eigenvalues of $ G=G_{1}\triangledown G_{2} $ consists of the eigenvalue $ 1-\dfrac{1}{r_{1}+n_{2}}\lambda_{1k}A(G_{1}) $, $k=2,\dots, n_{1} $, the eigenvalues $  1-\dfrac{1}{r_{2}+n_{1}}\lambda_{2k}A(G_{1}), k=2,\dots, n_{2} $ and the remaining two eigenvalues are given by the quotient matrix
\begin{equation} \label{Qmat of join of two graphs}
\begin{pmatrix}
 \dfrac{n_{2}}{r_{1}+n_{2}}& \dfrac{-n_{2}}{\sqrt{(r_{1}+n_{2})(r_{2}+n_{1})}}\\
\dfrac{-n_{1}}{\sqrt{(r_{1}+n_{2})(r_{2}+n_{1})}} & \dfrac{n_{1}}{r_{2}+n_{1}}
\end{pmatrix}.
\end{equation}
\end{corollary}

Since $ G_{1} $ and $ G_{2} $ are regular graphs, we observe that the two eigenvalues of matrix \eqref{Qmat of join of two graphs} are the largest and the smallest normalized Laplacian eigenvalue of $ G=G_{1}\triangledown G_{2} $.
\begin{proposition}
The largest and the smallest normalized Laplacian eigenvalues of $ G_{1}\triangledown G_{2} $ are the eigenvalues of the matrix \eqref{Qmat of join of two graphs}.

\end{proposition}
\begin{proposition}\label{spec of join of two graphs}
\begin{itemize}
\item[ \bf{(i)} ] The normalized Laplacian eigenvalues of the complete bipartite graph $ K_{a,b}=K_{a}\triangledown K_{b} $ are 
\[ 
\left  \{ 0,1^{[a+b-2]},2 \right \}.
 \]
\item[ \bf{(ii)} ] The normalized Laplacian eigenvalues of the complete split graph $CS_{\omega,n-\omega}=K_{\omega}\triangledown \overline{K}_{n-\omega}$, with clique number $ \omega $ and independence number $ n-\omega $ are
\[ 
\left  \{ 0,\left (\dfrac{n}{n-1}\right)^{[\omega-1]},\dfrac{2n-\omega+1}{n-1} \right \}.
 \]
\item[ \bf{(iii)}] The normalized Laplacian eigenvalues of the cone graph $ C_{a,b}=C_{a}\triangledown \overline{K}_{b} $ consists of the eigenvalues $ 1-\dfrac{1}{2+b}2\cos\left( \dfrac{2\pi k}{n}\right) $, where $ k=2,\dots,n-1 $, the simple eigenvalues $ 0 $ and $ \dfrac{2b+2}{b+2} $.
\item[ \bf{(iv)} ] The normalized Laplacian eigenvalues of the wheel graph $ W_{n}=C_{n-1}\triangledown K_{1} $ consists of the eigenvalues $ 1-\dfrac{2}{3}\cos\left( \dfrac{2\pi k}{m}\right) $, where $ k=2,\dots,n-2 $, and the eigenvalues $ \left \{0,\dfrac{4}{3}\right \}. $
\end{itemize}
\end{proposition}
\noindent\textbf{Proof.} $ \textbf{(i)} $. This follows from Corollary \ref{join of two graphs},  by taking $ n_{1}=a,n_{2}=b, r_{1}=r_{2}=0 $ and $ \lambda_{1k}=0,$ for $ k=2,\dots,a$ and $ \lambda_{2k} =0$ for each $ k=2,\dots,b$.\\
$ \textbf{(ii)} $. We recall that the adjacency spectrum of $ K_{\omega} $ is $ \{ \omega-1,-1^{(\omega-1)} \} $. Now, the result follows from Corollary \ref{join of two graphs} by taking $ n_{1}=\omega, n_{2}=n-\omega, r_{1}=\omega-1, r_{2}=0, \lambda_{1k}=-1,$ for $ k=2,\dots,\omega $ and $ \lambda_{2k}=0$ for $k=2,3,\dots,n-\omega .$ \\
$ \textbf{(iii)} $. Since adjacency spectrum of $ C_{n} $ is $ \left\lbrace 2\cos\left( \frac{2\pi k}{n}\right) : k=1,2,\dots,n\right\rbrace  $, by taking $ n_{1}=a,n_{2}=b, r_{1}=2, r_{2}=0 $ and $ \lambda_{1k} =2\cos\left( \frac{2\pi i}{m}\right )$ for $ k=2,3,\dots,a-1 $ and $ \lambda_{2,k} $ for $ k=2,\dots,b-1 $ in Corollary \ref{join of two graphs}, we get the required eigenvalues.\\
$ \textbf{(iv)} $. This is a special case of part \textbf{(iii)} with $ a=n-1 $ and $ b=1 .$ \qed

A \emph{friendship} graph $ F_{n} $ is a graph of order $ 2n+1 $, obtained by joining $ K_{1} $ with $ n $ copies of $ K_{2} $, that is, $ F_{n}=K_{1}\triangledown (nK_{2}) =K_{1,n}[K_{1},\underbrace{K_{2},K_{2},\dots,K_{2}}_{n}]$, where $K_1$ corresponds to the root vertex (vertex of degree greater than one) in $K_{1,n}$.  In particular, replacing some of $ K_{2} $'s by $ K_{1} $'s in $F_n$ we get \emph{firefly} type graph denoted by $ F_{p,n-p}$ and written as 
\begin{equation*}
F_{p,n-p}=K_{1,n}[K_{1},\underbrace{K_{1},K_{1},\dots,K_{1}}_{p}\underbrace{K_{2},K_{2},\dots,K_{2}}_{n-p}].
\end{equation*}
A \emph{generalized or multi-step wheel network} $ W_{a,b} $ is a graph derived from $ a $ copies of $ C_{b} $ and $ K_{1} $, in a such a way that all vertices of each $ C_{b} $ are adjacent to $ K_{1} $. Its order is $ ab+1 $ and can be written as $ W_{a,b}=K_{1}\triangledown (aC_{b})= K_{1,a}[K_{1},\underbrace{C_{b},\dots,C_{b}}_{a} ]$.

The normalized Laplacian eigenvalues of the \emph{friendship} graph $ F_{n} $, the \emph{firefly} type graph $ F_{p,n-p}$ and $ W_{a,b} $ are given by the following result.
\begin{proposition}\begin{itemize}
\item[ \bf{(i)}] The normalized Laplacian eigenvalues of $ F_{n} $ are 
\[ \bigg\{0, \left ( \dfrac{1}{2} \right )^{n-1}, \left ( \dfrac{3}{2} \right )^{n+1} \bigg\}. \]

\item[ \bf{(ii)}] The normalized Laplacian eigenvalues of $ F_{p,n-p} $  are 
\[ 
\bigg\{0, \left ( \dfrac{1}{2} \right )^{[n-p-1]}, 1^{[p-1]}, \left ( \dfrac{3}{2} \right )^{[n-p]}, \dfrac{5\sqrt{2n-p}\pm\sqrt{2n+7p}}{4\sqrt{2n-p}} \bigg\}.
 \]

\item[ \bf{(iii)}] The normalized Laplacian eigenvalues of $ W_{a,b} $ consists of the eigenvalues $ 0$, the eigenvalue $ \dfrac{4}{3} $ and the eigenvalues $ 1-\dfrac{2}{3}\cos\left (\dfrac{2\pi k}{b}\right ),$ for  $ k=2,\dots,b. $

\end{itemize}
\end{proposition}
\noindent\textbf{Proof.} \textbf{(i)}. By Theorem \ref{joined union} and the definition of $ F_{n} $, we have 
\[ \alpha_{1}=2n, \alpha_{2}=\dots=\alpha_{n+1}=1~~ \text{and}~~ r_{1}=0, r_{2}=\dots =r_{n+1}=1. \]
So, by Theorem \ref{joined union}, we see that $ \dfrac{3}{2} $ is the normalized Laplacian eigenvalues of $ F_{n} $ with multiplicity $ n. $ The remaining eigenvalues are given by the block matrix
 \begin{equation}\label{Qmat of Fn}
\left (\begin{array}{c|c c c c}
1 & \dfrac{-1}{\sqrt{n}}& \dots & \dfrac{-1}{\sqrt{n}}&\dfrac{-1}{\sqrt{n}}\\
\hline
\dfrac{-1}{2\sqrt{n}} & \dfrac{1}{2} & \dots & 0& 0\\
\vdots & \vdots &\ddots& \vdots & \vdots\\
\dfrac{-1}{2\sqrt{n}}& 0 & \dots &  \dfrac{1}{2} & 0\\
\dfrac{-1}{2\sqrt{n}}& 0 & \dots & 0 & \dfrac{1}{2}
\end{array}\right ).
\end{equation}
Clearly, $ \dfrac{1}{2} $ is the normalized laplacian eigenvalue of \eqref{Qmat of Fn} with multiplicity $ n-1 $ and the remaining two eigenvalues of block matrix \eqref{Qmat of Fn} are given by the quotient matrix
\[ \begin{pmatrix}
1 & \frac{-n}{\sqrt{n}}\\
\dfrac{-1}{2\sqrt{n}} & \dfrac{1}{2}

\end{pmatrix}.
 \]
First we will prove \textbf{(iii)} and then \textbf{(ii)}.\\
\textbf{(iii)}. As in part \textbf{(iii)} of Proposition \ref{spec of join of two graphs}, we see that $ 1-\dfrac{2}{3}\cos\left (\dfrac{2\pi k}{b}\right ),$ for  $ k=2,\dots,b. $ are the normalized Laplacian eigenvalues of $ W_{a,b} $. The other eigenvalues are given by the block matrix
\begin{equation*}
\left( \begin{array}{c| c c c c}
1 & \dfrac{-b}{\sqrt{3ab}}& \dots & \dfrac{-b}{\sqrt{3ab}}&\dfrac{-b}{\sqrt{3ab}}\\
\hline
\dfrac{-1}{\sqrt{3ab}} & \dfrac{1}{3} & \dots & 0 &  0\\
\vdots & \vdots &\ddots& \vdots & \vdots\\
\dfrac{-1}{\sqrt{3ab}} & 0 & \dots& \dfrac{1}{3} & 0\\
\dfrac{-1}{\sqrt{3ab}}& 0 & \dots & 0 & \dfrac{1}{3}
\end{array}\right).
\end{equation*}
Now, as in part \textbf{(i)}, $ \bigg\{0,\left (\dfrac{1}{3}\right )^{a-1}, \dfrac{4}{3}\bigg\} $ are the remaining normalized Laplacian eigenvalues of $ W_{a,b} $. \\
\textbf{(ii)}. Since $ \alpha_{1}=p+2(n-p)=2n-p $ and $ \alpha_{2}=\dots=\alpha_{2n+1-p}=1 $, so by Theorem \ref{joined union}, with $ r_{1}=\dots=r_{p+1}=0, ~ r_{p+2}=\dots=r_{2n+1-p}=1,$ we see that $ \dfrac{3}{2} $ is the normalized Laplacian eigenvalue of $ F_{p,n-p} $ with multiplicity $ n-p $. The other normalized Laplacian eigenvalues of $ F_{p,n-p} $ are given by the block matrix
\begin{equation}\label{Qmat Fp,n-p}
\left(
\begin{array}{c| c c c| c c c}
1 & \dfrac{-1}{\sqrt{2n-p}} & \dots & \dfrac{-1}{\sqrt{2n-p}}& \dfrac{-2}{\sqrt{2(2n-p)}}& \dots & \dfrac{-2}{\sqrt{2(2n-p)}}\\
\hline
\dfrac{-1}{\sqrt{2n-p}} & 1 & \dots & 0 & 0 &\dots & 0\\
\vdots & \vdots & \ddots & \vdots & \vdots &  \dots & \vdots\\
\dfrac{-1}{\sqrt{2n-p}} & 0 & \dots & 1 & 0 & \dots & 0\\

\hline
\dfrac{-2}{\sqrt{2(2n-p)}} & 0 & \dots & 0 & \dfrac{1}{2} & \dots & 0\\
\vdots & \vdots & \dots & \vdots & \vdots &  \ddots & \vdots\\
\dfrac{-2}{\sqrt{2(2n-p)}} & 0 & \dots & 0 & 0 & \dots & \dfrac{1}{2}\\
\end{array}
\right).
\end{equation}
By simple calculations, $ 1 $ and $ \dfrac{1}{2} $ are the normalized Laplacian eigenvalues of \eqref{Qmat Fp,n-p} and the remaining eigenvalues of block matrix \eqref{Qmat Fp,n-p} are given by the quotient matrix
\begin{equation}\label{qqmat}
\begin{pmatrix}
1 & \dfrac{-p}{\sqrt{2n-p}} & \dfrac{-2(n-p)}{\sqrt{2(2n-p)}}\\
\dfrac{-1}{\sqrt{2n-p}} & 1 & 0\\
\dfrac{-1}{\sqrt{2(2n-p)}} & 0 & \dfrac{1}{2}
\end{pmatrix}.
\end{equation}
Now, it is easy to see that $ 0 $ and $ \dfrac{5\sqrt{2n-p}\pm\sqrt{2n+7p}}{4\sqrt{2n-p}} $ are the normalized Laplacian eigenvalues of quotient matrix \eqref{qqmat}. \qed

\section{Normalized Laplacian eigenvalues of the power graphs of cyclic group $ \mathbb{Z}_{n} $ }
In this section, we consider the power graphs of finite cyclic group $ \mathbb{Z}_{n} $. As an application to Theorem \ref{joined union} and its consequences obtained in Section \ref{section 4.2}, we determine the normalized Laplacian eigenvalues of power graph of $ \mathbb{Z}_{n}. $

All groups are assumed to be finite and every cyclic group of order $ n $ is taken as isomorphic copy of integral additive modulo group $ \mathbb{Z}_{n} $ with identity denoted by $ 0 $. Let $ \mathcal{G} $ be a finite group of order $ n $ with identity  $ e $. The power graph of group $ \mathcal{G} $, denoted by $ \mathcal{P}(\mathcal{G}) $, is the simple graph with vertex set as the elements of group $ \mathcal{G} $ and two distinct vertices $ x,y\in \mathcal{G} $ are adjacent if and only if one is the positive power of the other, that is, $ x\sim y $ if and only if $ x^{i}=y $ or $ y^{j}=x $, for positive integers $i,j$ with  $ 2\leq i,j\leq n $. These graphs were introduced in \cite{kelarev}, see also \cite{sen}. Such graphs have valuable applications and are related to automata theory \cite{kelarev9}, besides being useful in characterizing finite groups.  We let $ U_{n}^{*}=\{x\in\mathbb{Z}_{n}: (x,n)=1 \}\cup\{0\} $, where $ (x,n) $ denotes greatest common divisor of $ x $ and $ n $. Our other group theory notations are standard and can be taken from \cite{roman}. More work on power graphs can be seen in \cite{cameron1, sen, mehreen, survey,tamiza} and the references therein. \\

The adjacency spectrum, the Laplacian and the signless Laplacian spectrum of power graphs of finite cyclic and dihedral groups have been investigated in  \cite{banerjee,sriparna,asma, mehreen1, panda}. The normalized Laplacian eigenvalues of power graphs of certain finite groups were studied in \cite{asma1}.\\

Let $ n $ be a positive integer. Then $ \tau(n) $ denotes the number of positive factors of $ n $, that is \[ \tau(n)=\sum\limits_{d|n} 1, \] where $ d|n $ denotes $ d $ divides $ n $.

The \emph{Euler's totient function} $ \phi(n) $ denotes the number of positive integers less or equal to $ n $ and relatively prime to $ n $.

We say $ n $ is in \emph{canonical decomposition} if $n=p_{1}^{n_{1}}p_{2}^{n_{2}}\dots p_{r}^{n_{r}} $, where $ r,n_{1},n_{2},\dots,n_{r} $ are positive integers and $ p_{1},p_{2},\dots,p_{r} $ are distinct primes.

The following result counts the values of $ \tau(n) $.
\begin{theorem}\cite{thomas}\label{tau function}
Let $ n $ be a positive integer with canonical decomposition $ n=p_{1}^{n_{1}}p_{2}^{n_{2}}\dots p_{r}^{n_{r}} $. Then 
\[ \tau(n)= (n_{1}+1) (n_{2}+1)\dots(n_{r}+1) \]
\end{theorem}
The following result gives some properties of Euler's totient function.
\begin{theorem}\cite{thomas}\label{phi function}
Let $ \phi $ be the Euler's totient function. Then following hold.
\begin{itemize}
\item[\bf(i)] $ \phi $ is multiplicative, that is 
$ \phi(st)=\phi(s)\phi(t), $ whenever $ s $ and $ t $ are relatively prime.
\item[\bf(ii)] Let $ n $ be a positive integer. Then  
$ \sum\limits_{d|n}\phi(d)=n. $
\item[\bf(iii)] Let $ p $ be a prime. Then  
$ \sum\limits_{i=1}^{l}\phi(p^{l})=p^{l}-1. $
\end{itemize}
\end{theorem}

Let $ n $ be a positive integer and $ d $ divides $ n $. The divisor $ d $ is the proper divisor of $ n $, if $ 1<d<n . $ Let $\mathbb{G}_{n}$ be a simple graph with vertex set as the proper divisor set $ \{d_{i}: 1,n\neq d_{i}|n,  ~ 1 \leq i \leq t\} $ and edge set $\{ d_{i}d_{j}: d_{i}|d_{j}, ~ 1 \leq i< j\leq t\} $, for $ 1\leq i<j\leq t $. If the \emph{canonical decomposition} of $n$ is $ n=p_{1}^{n_{1}}p_{2}^{n_{2}}\dots p_{r}^{n_{r}} $, then by Theorem \ref{tau function}, the order of the graph  $\mathbb{G}_{n}$ is  $ |V( \mathbb{G}_{n})|=\prod\limits_{i=1}^{r}(n_{i}+1)-2. $ Also, $\mathbb{G}_{n}$ is a connected graph \cite{bilal}, provided $ n $ is neither a prime power nor the product of two distinct primes. In \cite{mehreen}, $\mathbb{G}_{n}$ is used as the underlying graph for studying the power graph of finite cyclic group $ \mathbb{Z}_{n} $ and it has been shown that for each proper divisor $ d_{i} $ of $ n $, $ \mathcal{P}(\mathbb{Z}_{n})$ has a complete subgraph of order $ \phi(d_{i}) $.  \\

The following theorem shows that $ \dfrac{n}{n-1} $ is always the normalized Laplacian eigenvalue of the power graph $ \mathcal{P}(\mathbb{Z}_{n})$.

\begin{theorem}\label{mul of alpha n/n-1 of Zn}
Let $ \mathbb{Z}_{n} $ be a finite cyclic group of order $ n\geq 3 $. Then $ \dfrac{n}{n-1} $ is normalized Laplacian eigenvalue of $ \mathcal{P}(\mathbb{Z}_{n}) $ with multiplicity at least $ \phi(n)$.
\end{theorem}
\noindent\textbf{Proof.} Let $ \mathbb{Z}_{n} $ be the cyclic group of order $ n\geq 3 $. Then the identity $ 0 $ and invertible elements of  the group $ \mathbb{Z}_{n} $ in the power graph $ \mathcal{P}(\mathbb{Z}_{n}) $ are adjacent to every other vertex in $ \mathcal{P}(\mathbb{Z}_{n}) $. Since it is well known that the number of invertible elements of $ \mathbb{Z}_{n} $ are $ \phi(n) $ in number, so the induced power graph $ \mathcal{P}(U_{n}^{*}) $ is the complete graph $ K_{\phi(n)+1} $. Thus, by Theorem \ref{mehreen}, we see that $ \mathcal{P}(\mathbb{Z}_{n})=K_{\phi(n)+1}\triangledown \mathcal{P}(\mathbb{Z}_{n}\setminus U_{n}^{*})$. By applying Corollary \ref{join of two graphs}, we get \begin{equation*}
1-\dfrac{1}{r_{1}+\alpha_{1}}(-1)=1+\dfrac{1}{\phi(n)+n-\phi(n)-1}=\dfrac{n}{n-1}
\end{equation*}
as the normalized Laplacian eigenvalue with multiplicity at least $ \phi(n) $, since $ \dfrac{n}{n-1} $ can also be the normalized Laplacian eigenvalue of quotient matrix \eqref{Qmat of join of two graphs}.\qed

If $ n=p^{z} $, where $ p $ is prime and $ z $ is a positive integer, then we have following observation.

\begin{corollary}\label{NL spectra of p^z}
If $ n=p^{z} $, where $ p $ is prime and $ z $ is a positive integer, then the normalized Laplacian eigenvalues of  $\mathcal{P}(\mathbb{Z}_n) $ are 
\[  \left \{ 0, \left (\dfrac{n}{n-1}\right )^{[(n-1)]} \right \}. \]
\end{corollary}
\noindent\textbf{Proof.} If $ n=p^z $, where $ p $ is prime and $ z $ is a positive integer, then as shown in \cite{sen}, $ \mathcal{P}(\mathbb{Z}_n)$ is isomorphic to the complete graph $ K_{n}$ and hence the result follows.\qed

The next observation gives the normalized Laplacian eigenvalues of $ \mathcal{P}(\mathbb{Z}_n) $, when $ n $ is product of two primes.
\begin{corollary}\label{NL spectra of pq}
Let $ n=pq$ be product of two distinct primes. Then the normalized Laplacian eigenvalues of $ \mathcal{P}(\mathbb{Z}_n) $ are 
\[ \left \{ 0,\left (\dfrac{n}{n-1}\right )^{[\phi(n)]}, \left (1+\dfrac{1}{q\phi(p)}\right )^{[\phi(p)-1]}, \left (1+\dfrac{1}{p\phi(q)}\right )^{[\phi(q)-1]} \right  \} \]  and the zeros of polynomial
\begin{align*}
p(x)=x\Bigg(x^2&-x\left (\dfrac{\phi(n)+1}{q\phi(p)}+\dfrac{\phi(p)+\phi(q)}{q\phi(p)+\phi(q)}+\dfrac{\phi(n)+1}{p\phi(q)}\right )+\dfrac{(\phi(n)+1)\phi(p)}{p\phi(q)(q\phi(p)+\phi(q))} \\
&+\dfrac{(\phi(n)+1)^{2}}{n\phi(n)}+\dfrac{(\phi(n)+1)\phi(q)}{q\phi(p)(q\phi(p)+\phi(q))} \Bigg).
\end{align*}

\end{corollary}
\noindent\textbf{Proof.} If $ n=pq $, where $p$ and $q$, ($ p<q) $) are primes, then $ \mathcal{P}(\mathbb{Z}_{n}) $ \cite{tamiza} can be written as
\[  \mathcal{P}(\mathbb{Z}_n)=(K_{\phi(p)}\cup K_{\phi(q)})\triangledown K_{\phi(n)+1}= P_{3}[K_{\phi(p)},K_{\phi(pq)+1},K_{\phi(q)} ] . \]
By Theorem \ref{mul of alpha n/n-1 of Zn}, $ \dfrac{n}{n-1} $ is the normalized Laplacian eigenvalue with multiplicity $ \phi(n) $. Again, by Theorems \ref{joined union} and \ref{zn}, we see that $ \dfrac{1}{q\phi(p)} $ and $ \dfrac{1}{q\phi(p)} $ are the normalized Laplacian eigenvalues of $ \mathcal{P}(\mathbb{Z}_n) $ with multiplicity $ \phi(p)-1 $ and $ \phi(q)-1 $, respectively. The remaining three normalized Laplacian eigenvalues are given by the following matrix
\begin{equation*} \label{Qmat pq}
\begin{pmatrix}
\dfrac{\phi(n)+1}{q\phi(p)} & \dfrac{-(\phi(n)+1)}{\sqrt{q\phi(p)(q\phi(p)+\phi(q))}} & 0\\
\dfrac{-\phi(p)}{\sqrt{q\phi(p)(q\phi(p)+\phi(q))}} &\dfrac{\phi(p)+\phi(q)}{q\phi(p)+\phi(q)} & \dfrac{-\phi(q)}{\sqrt{p\phi(q)(q\phi(p)+\phi(q))}}\\
0 & \dfrac{-(\phi(n)+1)}{\sqrt{p\phi(q)(q\phi(p)+\phi(q))}} & \dfrac{\phi(n)+1}{p\phi(q)}
\end{pmatrix}.
\end{equation*}\qed

By Corollaries \ref{NL spectra of p^z} and \ref{NL spectra of pq}, we have the following proposition.
\begin{proposition}
Equality holds in Theorem \eqref{mul of alpha n/n-1 of Zn}, if n is some prime or product of two primes.
\end{proposition}

The following theorem \cite{mehreen} shows that the power graph of a finite cyclic group $ \mathbb{Z}_{n} $ can be written as the joined union, each of whose components are cliques.

\begin{theorem}\label{mehreen}
If $ \mathbb{Z}_{n} $ is a finite cyclic group of order $n\ge 3$, then the power graph $\mathcal{P}(\mathbb{Z}_{n})$ is given by  
\begin{equation*}
\mathcal{P}(\mathbb{Z}_{n}) =K_{\phi(n)+1}\triangledown \mathbb{G}_{n}[K_{\phi(d_{1})},K_{\phi(d_{2})},\dots,K_{\phi(d_{t})}],
\end{equation*} 
where $\mathbb{G}_{n}$ is the graph of order $t$ defined above. 
\end{theorem}

Using Theorem \ref{joined union} and its consequences, we can compute the normalized Laplacian eigenvalues of $ \mathcal{P}(\mathbb{Z}_{n}) $ in terms of the adjacency spectrum of $ K_{\omega} $ and zeros of the characteristic polynomial of the quotient matrix.

We form a connected graph $ H=K_{1}\triangledown \mathbb{G}_{n} $ which is of diameter at most two if $ \mathbb{G}_{n} $ is not complete, otherwise its diameter is $ 1 $. In the following result, we compute the normalized Laplacian eigenvalues of the power graph of $ \mathbb{Z}_{n} $ by using Theorems \ref{joined union} and \ref{mehreen}.

\begin{theorem}\label{zn}
The normalized Laplacian eigenvalues of $ \mathcal{P}(\mathbb{Z}_n) $ are
\begin{equation*}
\left\lbrace \left (\dfrac{n}{n-1}\right )^{(\phi(n))}, \left (\dfrac{\phi(d_{1})+\alpha_{2}}{\phi(d_{1})+\alpha_{2}-1}\right )^{[\phi(d_{1})-1]},\dots, \left (\dfrac{\phi(d_{t})+\alpha_{r+1}}{\phi(d_{t})+\alpha_{t+1}-1}\right )^{[\phi(d_{t})-1]} \right\rbrace
\end{equation*}
and the $ t+1 $ eigenvalues of  matrix \eqref{quotient matrix of Z_n}
\begin{equation}\label{quotient matrix of Z_n}
\begin{pmatrix}
\dfrac{n-1-\phi(n)}{n-1}& \dfrac{-\phi(d_{1})a_{12}}{\sqrt{(\phi(n)+\alpha_{1})(r_{2}+\alpha_{2})}}&\ldots & \dfrac{-\phi(d_{t})a_{1(t+1)}}{\sqrt{b}}\\
\dfrac{-\phi(d_{1})a_{21}}{\sqrt{(r_{2}+\alpha_{2})(\phi(n)+\alpha_{1})}}&\dfrac{\alpha_{2}}{\alpha_2+r_2}& \ldots & \dfrac{-\phi(d_{t})a_{2(t+1)}}{\sqrt{(r_{2}+\alpha_{2})(r_{t+1}+\alpha_{t+1})}}\\
\vdots &\vdots  &\ddots &\vdots\\
\dfrac{-\phi(d_{1})a_{(t+1)1}}{\sqrt{(r_{n}+\alpha_{t+1})(\phi(n)+\alpha_{1})}}& \dfrac{-\phi(d_{2})a_{(t+1)2}}{\sqrt{(r_{n}+\alpha_{t+1})(r_{2}+\alpha_{2})}}&\ldots &\dfrac{\alpha_{t+1}}{\alpha_{t+1}+r_{t+1}}
\end{pmatrix},
\end{equation} 
where, for $i\ne j$,
\begin{equation*}
a_{ij}=\begin{cases}
1, & v_i\sim v_j\\
0, & v_i\nsim v_j

\end{cases},
\end{equation*}
$ b=(\phi(n)+\alpha_{1})(r_{t+1}+\alpha_{t+1}) $
and $ r_{i}=\phi(d_{i})-1, ~ \text{for} ~i=2,\dots, t+1. $
\end{theorem}
\noindent\textbf{Proof.} Let $ \mathbb{Z}_{n} $ be a finite cyclic group of order $ n $. Since the identity element $ 0 $ and the $ \phi(n) $ generators of the group $ \mathbb{Z}_{n} $ are adjacent to every other vertex of $ \mathcal{P}(\mathbb{Z}_{n}) $, therefore, by Theorem \ref{mehreen}, we have
\begin{equation*}
\mathcal{P}(\mathbb{Z}_{n}) =K_{\phi(n)+1}\triangledown \mathbb{G}_{n}[K_{\phi(d_{1})},K_{\phi(d_{2})},\dots,K_{\phi(d_{t})}]= H[K_{\phi(n)+1},K_{\phi(d_{1})},K_{\phi(d_{2})},\dots,K_{\phi(d_{t})}],
\end{equation*}
where $ H=K_{1}\triangledown \mathbb{G}_{n} $ is the graph with vertex set $ \{v_{1},\dots, v_{t+1}\} $. Taking $ G_{1}=K_{\phi(n)+1}$ and $ G_{i}=K_{\phi(d_{i-1})} $, for $ i=2,\dots,t+1 $, in  Theorem \ref{joined union} and using the fact that the adjacency spectrum of $K_{\omega}$ consists of the eigenvalue $\omega-1$ with multiplicity $1$ and the eigenvalue $-1$ with multiplicity $\omega-1$, it follows that  
\begin{equation*}
1-\dfrac{1}{r_{1}+\alpha_{1}}\lambda_{1k}A(G_{1}) =1-\dfrac{1}{r_{1}+\alpha_{1}}(-1)= 1+\dfrac{1}{\phi(n)+n-\phi(n)-1}=\dfrac{n}{n-1} 
\end{equation*}
is a normalized eigenvalue of $ \mathcal{P}(\mathbb{Z}_n) $ with multiplicity $ \phi(n)$. Note that we have used the fact that vertex $ v_{1} $ of graph $ H $ is adjacent to every other vertex of $ H $ and $ \alpha_1=\sum\limits_{ d|n, d\ne 1,n}\phi(d)=n-1-\phi(n)$, as $ \sum\limits_{ d|s}\phi(d) =s.$ Similarly, we can show that $ \dfrac{\phi(d_{1})+\alpha_{2}}{\phi(d_{1})+\alpha_{2}-1},\dots, \dfrac{\phi(d_{t})+\alpha_{t+1}}{\phi(d_{t})+\alpha_{t+1}-1} $ are the normalized Laplacian eigenvalues of $ \mathcal{P}(\mathbb{Z}_n) $ with multiplicities $ \phi(d_{1})-1, \dots, \phi(d_{t})-1 $, respectively.  The remaining normalized Laplacian eigenvalues are the eigenvalues of quotient matrix  \eqref{quotient matrix of Z_n}.\qed\\

From Theorem \ref{zn}, it is clear that all the normalized Laplacian eigenvalues of the power graph $ \mathcal{P}(\mathbb{Z}_{n})$ are completely determined except the $t+1$ eigenvalues, which are the eigenvalues of matrix \eqref{quotient matrix of Z_n}. Further, it is also clear that matrix \eqref{quotient matrix of Z_n}  depends upon the structure of the graph $\mathbb{G}_{n} $, which is not known in general. However, if we give some particular value to $n$, then it may be possible to know the structure of graph $\mathbb{G}_{n}$ and hence about matrix \eqref{quotient matrix of Z_n}. This information may be helpful to determine the $t+1$ remaining normalized Laplacian eigenvalues of the power graph $ \mathcal{P}(\mathbb{Z}_{n})$.\\

We discuss some particular cases of Theorem \ref{zn}.\\

\begin{figure}
\centering

\begin{tikzpicture}
\draw[fill=black] (0,0) circle (2pt);\draw[fill=black] (2,0) circle (2pt);\draw[fill=black] (4,0) circle (2pt);
\draw[fill=black] (-1,-1) circle (2pt);\draw[fill=black] (1,-1) circle (2pt);\draw[fill=black] (3,-1) circle (2pt);

\draw[thin] (0,0)--(-1,-1); \draw[thin] (0,0)--(1,-1);\draw[thin] (2,0)--(-1,-1);
\draw[thin] (2,0)--(3,-1); \draw[thin] (4,0)--(1,-1);\draw[thin] (4,0)--(3,-1);

\node at (0,.3) {$ p $};\node at (2,.3) {$ q $};\node at (4,.3) {$ r $};
\node at (-1,-1.3) {$ pq $};\node at (1,-1.3) {$ pr $};\node at (3,-1.3) {$ qr $};
\end{tikzpicture}\quad
\begin{tikzpicture}
\draw[fill=black] (0,0) circle (2pt);\draw[fill=black] (2,0) circle (2pt);\draw[fill=black] (4,0) circle (2pt);
\draw[fill=black] (-1,-1) circle (2pt);\draw[fill=black] (1,-1) circle (2pt);\draw[fill=black] (3,-1) circle (2pt);
\draw[fill=black] (2.9,-2.3) circle (2pt);

\draw[thin] (0,0)--(2.9,-2.3);\draw[thin] (2,0)--(2.9,-2.3);\draw[thin] (3,-1)--(2.9,-2.3);
\draw[thin] (2,0)--(2.9,-2.3); \draw[thin] (4,0)--(2.9,-2.3);\draw[thin] (-1,-1)--(2.9,-2.3);

\draw[thin] (0,0)--(-1,-1); \draw[thin] (0,0)--(1,-1);\draw[thin] (2,0)--(-1,-1);
\draw[thin] (2,0)--(3,-1); \draw[thin] (4,0)--(1,-1);\draw[thin] (4,0)--(3,-1);

\node at (0,.3) {$ p $};\node at (2,.3) {$ q $};\node at (4,.3) {$ r $};
\node at (-1,-1.3) {$ pq $};\node at (1,-1.3) {$ pr $};\node at (3,-0.7) {$ qr $}; \node at (2.9,-2.7) {$ K_{1} $};
\end{tikzpicture}

\caption{Divisor graph $ \mathbb{G}_{pqr} $ and $ H=K_{1}\triangledown \mathbb{G}_{pqr} $. }
\end{figure}

Now, let $ n=pqr $, where $p,~ q,~ r$  with $ p<q<r$ are primes. From the definition of $\mathbb{G}_{n}$, the  vertex set and edge set of $\mathbb{G}_{n}$ are $\{ p, q, r, pq, pr, qr \}$ and $\{(p,pq),(p,pr),(q,pq),(q,qr),(r,pr),(r,qr)\} $, respectively, and is shown in Figure $ (4.1) $. Let $ H=K_{1}\triangledown \mathbb{G}_{n}.$ Then 
\begin{equation*}
 \mathcal{P}(\mathbb{Z}_n)= H[K_{\phi(n)+1}, K_{\phi(p)},K_{\phi(q)},K_{\phi(r)},K_{\phi(pq)},K_{\phi(pr)}, K_{\phi(qr)}].
\end{equation*}
By Theorem \ref{joined union}, we have \begin{align*}
\Big (\alpha_{1},\alpha_{2}, \alpha_{3},\alpha_{4},\alpha_{5},\alpha_{6},\alpha_{7}\Big)=&\Big(n-\phi(n)-1,~\phi(n)+1+\phi(pq)+\phi(pr),~\phi(n)+1\\
&+\phi(pq)+\phi(qr)),~\phi(n)+1+\phi(pr)+\phi(qr),~\phi(n)+1\\
&+\phi(p)+\phi(q),\phi(n)+1+\phi(p)+\phi(r),~\phi(n)+1+\phi(q)\\
&+\phi(r)\Big ).
\end{align*}
Also, for $ i=1,2,\dots,7$,~  $ r_{i}^{'}=\alpha_{i}+r_{i} $, we have
 \begin{align*}
\Big (r_{1}^{'},r_{2}^{'}, r_{3}^{'}, r_{4}^{'} ,r_{5}^{'},r_{6}^{'},r_{7}^{'}\Big)=& \Big(n-1,~\phi(n)+\phi(p)+\phi(pq)+\phi(pr),~\phi(n)+\phi(q)+\phi(pq)\\
&+\phi(qr),\phi(n)+\phi(r)+\phi(pr)+\phi(qr),~\phi(n)+\phi(pq)+\phi(p)\\
&+\phi(q),~\phi(n)+\phi(pr)+\phi(p)+\phi(r),~\phi(n)+\phi(qr)+\phi(q)\\
&+\phi(r)\Big ). 
\end{align*}
Now, by Theorem \ref{mul of alpha n/n-1 of Zn}, $ \dfrac{n}{n-1} $ is the normalized Laplacian eigenvalue with multiplicity $ \phi(n) $.  Using the above information and Theorem \ref{zn}, second distinct normalized Laplacian eigenvalue is $ 1+\dfrac{1}{r_{2}+\alpha_{2}}=1+\dfrac{1}{\phi(n)+\phi(p)+\phi(pq)+\phi(pr)} $ with multiplicity $ \phi(p)-1 $. In a similar way, we see that other eigenvalues are 
\begin{align*}
&1+\dfrac{1}{\phi(n)+\phi(q)+\phi(pq)+\phi(qr)},1+\dfrac{1}{\phi(n)+\phi(r)+\phi(pr)+\phi(qr)},\\
&1+\dfrac{1}{\phi(n)+\phi(pq)+\phi(q)+\phi(q)},1+\dfrac{1}{\phi(n)+\phi(pr)+\phi(p)+\phi(r)},\\
&1+\dfrac{1}{\phi(n)+\phi(qr)+\phi(q)+\phi(r)}
\end{align*}
with multiplicities $ \phi(q)-1,~ \phi(r)-1,~ \phi(pq)-1,~ \phi(pr)-1,~ \phi(qr)-1$, respectively. The remaining $7$ eigenvalues are given by the following matrix
\begin{equation*}\label{pqr}
\begin{pmatrix}
z_{1}& -\phi(p)c_{12} & -\phi(q)c_{13} & -\phi(r)c_{14}& -\phi(pq)c_{15}& -\phi(pr)c_{16}& -\phi(qr)c_{17}\\
(\phi(n)+1)c_{21} & z_{2} & 0 & 0& -\phi(pq)c_{25} & -\phi(pr)c_{26}& 0 \\
(\phi(n)+1)c_{31} & 0 & z_{3} & 0& -\phi(pq)c_{35} & 0 & -\phi(qr)c_{37} \\
(\phi(n)+1)c_{41} & 0 & 0 & z_{4} & 0&  -\phi(pr)c_{46} & -\phi(qr)c_{47} \\
(\phi(n)+1)c_{51} & -\phi(p)c_{25} & -\phi(q)c_{35} & 0 & z_{5} & 0  & 0  \\
(\phi(n)+1)c_{61} & -\phi(p)c_{26} & 0 & -\phi(r)c_{64} & 0 & z_{6} & 0  \\
(\phi(n)+1)c_{71} & 0 & -\phi(q)c_{75} & -\phi(r)c_{74} & 0 & 0 &  z_{7} 
\end{pmatrix},
\end{equation*}
where, \begin{align*}
z_{1}=&\frac{n-\phi(n)-1}{n-1}, ~ z_{2}= \frac{\phi(n)+1+\phi(pq)+\phi(pr)}{\phi(n)+\phi(p)+\phi(pq)+\phi(pr)},\\
&~ z_{3}= \frac{\phi(n)+1+\phi(pq)+\phi(qr)}{\phi(n)+\phi(q)+\phi(pq)+\phi(qr)},~z_{4}=\dfrac{\phi(n)+1+\phi(pr)+\phi(qr)}{\phi(n)+\phi(r)+\phi(pr)+\phi(qr)},\\
&z_{5}= \dfrac{\phi(n)+1+\phi(p)+\phi(q)}{\phi(n)+\phi(pq)+\phi(p)+\phi(q)},~ z_{6}= \dfrac{\phi(n)+1+\phi(p)+\phi(r)}{\phi(n)+\phi(pr)+\phi(p)+\phi(r)},\\ &z_{7}=\dfrac{\phi(n)+1+\phi(q)+\phi(r)}{\phi(n)+\phi(qr)+\phi(q)+\phi(r)}\\
\text{and}~c_{ij}&=c_{ji}=\dfrac{1}{\sqrt{(r_{i}+\alpha_{i})(r_{j}+\alpha_{j})}}.
\end{align*}\\

Next, we discuss the normalized Laplacian eigenvalues of the finite cyclic group $\mathbb{Z}_{n}$, with $n=p^{n_{1}}q^{n_{2}},$ where $ p<q $ are primes and $ n_{1}\leq  n_{2}$ are positive integers. We consider the case when both $ n_{1} $ and $ n_{2} $ are even, and the case when they are odd can be discussed similarly.
\begin{theorem} 
Let $ \mathcal{P}(\mathbb{Z}_{p^{n_{1}}q^{n_{2}}}) $ be the power graph of the finite cyclic group $\mathbb{Z}_{p^{n_{1}}q^{n_{2}}}$ of order $ n={p^{n_{1}}q^{n_{2}}} $, where $p<q$ are primes and $n_{1}=2m_{1}\leq n_{2}=2m_{2}$ are even positive integers. Then the normalized Laplacian eigenvalues of $ \mathcal{P}(\mathbb{Z}_{p^{n_{1}}q^{n_{2}}}) $ consists of eigenvalue set 
\begin{align*}
&\left (\dfrac{n}{n-1}\right )^{[\phi(n)]}, \left (\dfrac{n-q^{n_{2}}+1}{n-q^{n_{2}}}\right )^{[\phi(p)-1]},\\
&\qquad \vdots\\
&\left (\dfrac{p^{m_{1}-1}+q^{n_{2}}(p^{n_{1}}-p^{m_{1}-1})}{p^{m_{1}-1}+q^{n_{2}}(p^{n_{1}}-p^{m_{1}-1})}-1\right )^{[\phi(p^{m_{1}})-1]},\\
&\qquad\vdots\\
&\left (\dfrac{p^{n_{1}-1}+q^{n_{2}}\phi(p^{n_{1}})}{p^{n_{1}-1}+q^{n_{2}}\phi(p^{n_{1}})-1}\right )^{\left [\phi\left (p^{n_{1}}\right )-1\right ]},\left (\dfrac{n-p^{n_{1}}+1}{n-p^{n_{1}}}\right )^{\left [\phi\left (q\right )-1\right ]},\\
&\qquad\vdots\\
&\left (\dfrac{q^{m_{2}-1}+p^{n_{1}}(q^{n_{2}}-q^{m_{2}-1})}{q^{m_{2}-1}+p^{n_{1}}(q^{n_{2}}-q^{m_{2}-1})-1}\right )^{\left [\phi\left (q^{m_{2}}\right )-1\right ]},\\
&\qquad\vdots\\
&\left (\dfrac{q^{n_{2}-1}+p^{n_{1}}\phi(q^{n_{2}})}{q^{n_{2}-1}+p^{n_{1}}\phi(q^{n_{2}})-1}\right )^{\left [\phi\left (q^{n_{2}}\right )-1\right ]},\left (\dfrac{\phi(p)+\phi(q)+(q^{n_{2}}-1)(p^{n_{1}}-1)+1}{\phi(p)+\phi(q)+(q^{n_{2}}-1)(p^{n_{1}}-1)}\right )^{\left [\phi\left (pq\right )-1\right ]},\\
&\qquad\vdots\\
&\left (\dfrac{q^{n_{2}}(p^{n_{1}}-1)+q^{m_{2}}-q^{m_{2}-1}(p^{n_{1}}-p)}{q^{n_{2}}(p^{n_{1}}-1)+q^{m_{2}}-q^{m_{2}-1}(p^{n_{1}}-p)-1}\right )^{\left [\phi\left (pq^{m_{2}}\right )-1\right ]},\\
&\qquad\vdots\\
& \left (\dfrac{pq^{n_{2}}+\phi(p^{n_{1}})(q^{n_{2}}-q)}{pq^{n_{2}}+\phi(p^{n_{1}})(q^{n_{2}}-q)-1} \right )^{\left [\phi\left (pq^{n_{2}}\right )-1\right ]},\\
&\qquad\vdots\\
&\left (\dfrac{p^{m_{1}}+p^{n_{1}}(q^{n_{2}}-1)-p^{m_{1}-1}(q^{n_{2}}-q)}{p^{m_{1}}+p^{n_{1}}(q^{n_{2}}-1)-p^{m_{1}-1}(q^{n_{2}}-q)-1} \right )^{\left [\phi\left (p^{m_{1}}q\right )-1\right ]},\\
&\qquad\vdots\\
&\left (\dfrac{n+p^{m_{1}}q^{m_{2}}+p^{m_{1}-1}q^{m_{2}-1}-\phi(p^{m_{1}}q^{m_{2}})-p^{n_{1}}q^{m_{2}-1}-p^{m_{1}-1}q^{n_{2}}}{n+p^{m_{1}}q^{m_{2}}+p^{m_{1}-1}q^{m_{2}-1}-\phi(p^{m_{1}}q^{m_{2}})-p^{n_{1}}q^{m_{2}-1}-p^{m_{1}-1}q^{n_{2}}-1} \right )^{\left [\phi\left (p^{m_{1}}q^{m_{2}}\right )-1\right ]},\\
&\qquad\vdots\\
&\left (\dfrac{p^{m_{1}}q^{n_{2}}+\phi(q^{n_{2}})(p^{n_{1}}-p^{m_{1}})}{p^{m_{1}}q^{n_{2}}+\phi(q^{n_{2}})(p^{n_{1}}-p^{m_{1}})-1} \right )^{\left [\phi\left (p^{m_{1}}q^{n_{2}}\right )-1\right ]},\\
&\qquad\vdots\\
&\left (\dfrac{p^{n_{1}}q+\phi(p^{n_{1}})(q^{n_{2}}-q)}{p^{n_{1}}q+\phi(p^{n_{1}})(q^{n_{2}}-q)-1} \right )^{\left [\phi\left (p^{n_{1}}q\right )-1\right ]},\\
&\qquad\vdots\\
& \left ( \dfrac{p^{n_{2}}q^{m_{1}}+\phi(p^{n_{1}})(q^{n_{2}}-q^{m_{1}})}{p^{n_{2}}q^{m_{1}}+\phi(p^{n_{1}})(q^{n_{2}}-q^{m_{1}})-1} \right )^{\left [\phi\left (p^{n_{1}}q^{m_{2}}\right )-1\right ]},\\
&\qquad\vdots\\
&\left ( \dfrac{p^{n_{1}}q^{n_{2}-1}+\phi(n)}{p^{n_{1}}q^{n_{2}-1}+\phi(n)-1} \right )^{\left [\phi\left (p^{n_{1}}q^{n_{2}-1}\right )-1\right ]}\end{align*}
and the remaining eigenvalues are given by matrix \eqref{quotient matrix of Z_n}.
\end{theorem}
\noindent\textbf{Proof.} Suppose that $n=p^{n_{1}}q^{n_{2}}$, where $ n_{1}=2m_{1} $ and $ n_{2}=2m_{2} $ are even with $ n_{1}\leq n_{2} $ and $ m_{1} $ and $ m_{2} $ are positive integers. Since the total number of divisors of $ n $ are $ (n_{1}+1)(n_{2}+1) $, so the order of $ \mathbb{G}_{p^{n_{1}}q^{n_{2}}} $ is $ (n_{1}+1)(n_{2}+1)-2 $. The proper divisor set of $n$ is 
\begin{align*}
D(n)=&\bigg\{p,p^{2},\dots,p^{m_{1}},\dots,p^{n_{1}}, q,q^2,\dots,q^{m_{2}},\dots,q^{n_{2}}, pq, pq^2,\dots,pq^{m_{2}},\dots,pq^{n_{2}}, \dots,\\
&p^{m_{1}}q, p^{m_{1}}q^2,\dots, p^{m_{1}}q^{m_{2}},\dots,p^{m_{1}}q^{n_{2}},\dots,p^{n_{1}}q, p^{n_{1}}q^2,\dots,p^{n_{1}}q^{m_{2}},\dots,p^{n_{1}}q^{n_{2}-1}\bigg\}. 
\end{align*}
By the definition of graph $ \mathbb{G}_{n} $, we see that $ p$ is not adjacent to $ p, q,q^{2},\cdots,q^{m_{2}},\cdots,q^{n_{2}} $. So we write adjacency of vertices in terms of iterations and avoid divisors outside the  set $ D(n) $. Thus, we observe that
\begin{align*}
p&\sim p^{i},p^{j}q^{k},~ \text{for}~ i=2,3,\dots,n_{1},~  j=1,2,\dots,n_{1},~ k=1,2,\dots,n_{2},  \\ 
&\vdots\\
p^{m_{1}}&\sim p^{i},p^{j}q^{k},~ \text{for}~ i=1,2,\dots,n_{1}, ~i\neq m_{1},~ j=m_{1},\dots,n_{1},~ k=1,2,\dots,n_{2},\\
&\vdots\\
p^{n_{1}}&\sim p^{i},p^{n_{1}}q^{j},~ \text{for}~ i=1,2,\dots,n_{1}-1,~  j=1,2,\dots,n_{2}-1,\\ 
q&\sim q^{i},p^{j}q^{k},~ \text{for}~ i=2,3,\dots,n_{2}, ~ j=1,2,\dots,n_{1},~ k=1,2,\dots,n_{2},\\
& \vdots\\
q^{m_{2}}&\sim q^{i},p^{j}q^{k},~ \text{for}~ i=1,2,3,\dots,n_{1}, ~i\neq m_{2}, ~j=1,2,3,\dots,n_{1},~ k=m_{2},\dots,n_{2}\\
&\vdots\\
q^{n_{2}}&\sim q^{i},p^{j}q^{n_{2}},~ \text{for}~ i=1,2,3,\dots,n_{2}-1,~ j=1,2,3,\dots,n_{1}-1,\\ 
pq&\sim p, q, p^{i}q^{j},~ \text{for } ~ i=1,2,3,\dots,n_{1},~ j=1,2,3,\dots,n_{2},\\
& \vdots\\
pq^{m_{2}}&\sim  p, q^{i},pq^{j},p^{k}q^{k} ~ \text{for}~ i=1,2,\dots,m_{2}, ~ j=1,2,3,\dots,n_{2},~ j\neq m_{2},~ k=2,3,\\
&\dots,n_{1},~l=m_{2},\dots,n_{2},\\
&\vdots\\
pq^{n_{2}}&\sim  p, q^{i}, pq^{j},p^{k}q^{n_{2}}, ~\text{for}~ i=1,2,\dots,n_{2},~ j=1,2,\dots,n_{2}-1,k=2,3,\dots,n_{1}-1,\\ 
& \vdots\\
p^{m_{1}}q &\sim  p^{i}, q,p^{m_{1}}q^{j},p^{k}q,p^{l}q^{m} ~\text{for}~ i=1,\dots,m_{1}, ~ j=2,3,\dots,n_{2},~k=1,\dots,m_{1}-1,\\
&l=m_{1}+1,\dots,n_{1},~ m=1,2,\dots,n_{2},\\
& \vdots\\
p^{m_{1}}q^{m_{2}} &\sim  p^{i}, q^{j}, p^{k}q^{l}~ \text{for}~ i=1,2,\dots,m_{1},j=1,2,\dots,m_{2},k=1,2,\dots,n_{1},l=1,\dots,n_{2},\\
&\vdots\\
p^{m_{1}}q^{n_{2}} &\sim  p^{i}, q^{j},p^{k}q^{n_{2}}, p^{i}q^{j}~ \text{for}~ i=1,2,\dots,m_{1},j=1,2,\dots,n_{2},k=m_{1}+1,\dots,n_{1}-1,\\
&\vdots\\
p^{n_{1}}q &\sim  p^{i}, q ,p^{j}q, p^{n}q^{k}~ \text{for}~ i=1,2,\dots,n_{1},~j=1,2,\dots,n_{1}-1,~k=2,3\dots,n_{2}-1,\\
&\vdots\\
p^{n_{1}}q^{m_{2}} &\sim  p^{i}, q^{j} ,p^{n_{1}}q^{k}, p^{i}q^{j}~ \text{for}~ i=1,2,\dots,n_{1},j=1,2,\dots,m_{2},k=m_{2}+1,\dots,n_{2}-1,\\
&\vdots\\
p^{n_{1}}q^{n_{2}-1} &\sim  p^{i}, q^{j} , p^{i}q^{j}~ \text{for}~ i=1,2,\dots,n_{1},~j=1,2,\dots,n_{2}-1.\\
\end{align*}
Therefore, by Theorem \ref{mehreen}, we have \begin{align*}
\mathcal{P}(\mathbb{Z}_{n})= &K_{\phi(n)+1}\triangledown \mathbb{G}_{n}[K_{\phi(p)},\dots,K_{\phi(p^{m_{1}})},\dots,K_{\phi(p^{n_{1}})},K_{\phi(q)},\dots,K_{\phi(q^{m_{2}})},\dots,K_{\phi(q^{n_{2}})},\\
&K_{\phi(pq)},\dots,K_{\phi(pq^{m_{2}})},\dots,K_{\phi(pq^{n_{2}})},\dots,K_{\phi(p^{m_{1}}q)},\dots,K_{\phi(p^{m_{1}}q^{m_{2}})},\dots,K_{\phi(p^{m_{1}}q^{n_{2}})},\\
&\dots,K_{\phi(p^{n_{1}}q)},\dots,K_{\phi(p^{n_{1}}q^{m_{2}})},\dots,K_{\phi(p^{n_{1}}q^{n_{2}-1})}].
\end{align*}
Now, by using Theorem \ref{joined union}, we calculate the values of $ \alpha_{i} $'s and $ r_{i}+\alpha_{i}=r_{i}^{'} $'s. We recall that \cite{thomas} $\phi(xy)=\phi(x)\phi(y),$ provided that $(x,y)=1$, $\sum\limits_{i=1}^{k}\phi(p^i)=p^k-1$ and $ \sum\limits_{d|s}\phi(d)=s $. Using this information and definition of $ \alpha_{i} $'s, we have \[
\alpha_{1}=\sum\limits_{1,n\neq d|n}\phi(d)=n-1-\phi(n). \]
Similarly, \begin{align*}
\alpha_{2}=& \phi(p^{2})+\dots+\phi(p^{m_{1}})+\dots+\phi(p^{n_{1}})+\phi(pq)+\dots+\phi(pq^{m_{2}})+\dots\phi(pq^{n_{2}})\\
&+\phi(p^{m_{1}}q)+\dots+\phi(p^{m_{1}}q^{m_{2}})+\dots\phi(p^{m_{1}}q^{n_{2}})+\dots+\phi(p^{n_{1}}q)+\dots+\phi(p^{n_{1}}q^{m_{2}})\\
&+\dots+\phi(p^{n_{1}}q^{n_{2}-1})+\phi(n)+1\\
=&\sum\limits_{1,p,n\neq d|n}\phi(d)-[\phi(q)+\dots+\phi(q^{m_{1}})+\dots+\phi(q^{n_{1}})]+\phi(n)+1\\
=&n-1-\phi(p)-\phi(n)-[q^{n_{2}}-1]+\phi(n)+1=n-\phi(p)-q^{n_{2}}+1.
\end{align*} Proceeding in the same way as above, other $ \alpha_{i} $'s are 
\begin{align*}
\alpha_{3}=&q^{n_{2}}(p^{n_{1}}-p)+p-\phi(p^{2}),\\
& \vdots\\
\alpha_{m_{1}+1}=& p^{m_{1}-1}+q^{n_{2}}(p^{n_{1}}-p^{m_{1}-1}) -\phi(p^{m_{1}}),\\
\vdots\\
\alpha_{n_{1}+1}=&p^{n_{1}-1}+\phi(p^{n_{1}})(q^{n_{2}-1}-1),~\alpha_{n_{1}+2} =n-\phi(q)-p^{n_{1}}+1,\\
&\vdots\\
\alpha_{n_{1}+m_{1}+1}=&q^{m_{2}-1}+p^{n_{1}}(q^{n_{2}}-q^{m_{2}-1})-\phi(q^{m_{2}}),\\
\vdots\\
\alpha_{n_{1}+n_{2}+1} =&q^{n_{2}-1}+\phi(q^{n_{2}})(p^{n_{1}}-1),\\
\alpha_{n_{1}+n_{2}+2}=&\phi(p)+\phi(q)+1-\phi(pq)+(q^{n_{2}}-1)(p^{n_{1}}-1),\\
&\vdots\\
\alpha_{n_{1}+n_{2}+m_{1}+1}=&q^{n_{2}}(p^{n_{1}}-1)+q^{m_{2}}-q^{m_{2}-1}(p^{n_{1}}-p)-\phi(pq^{m_{2}}),\\
&\vdots\\
\alpha_{2n_{1}+n_{2}+1}=&pq^{n_{2}}-\phi(pq^{n_{2}})+\phi(p^{n_{1}})(q^{n_{2}}-q),\\
&\vdots\\
\alpha_{m_{1}n_{2}+n_{1}+2}=&p^{m_{1}}-\phi(p^{m_{1}}q)+p^{n_{1}}(q^{n_{2}}-1) -p^{m_{1}-1} (q^{v_{1}}-q),\\
&\vdots\\
 \alpha_{m_{1}n_{2}+n_{1}+m_{1}+1}=&n+p^{m_{1}-1}q^{m_{2}-1}+p^{m_{1}}q^{m_{2}}-2\phi(p^{m_{1}}q^{m_{2}})-q^{n_{2}}p^{m_{1}-1}-p^{n_{2}}q^{m_{2}-1},\\
&\vdots\\
\alpha_{(m_{1}+1)n_{2}+n_{1}+1}=&p^{m_{1}}q^{n_{2}}-\phi(p^{m_{1}}q^{n_{2}})+\phi(q^{n_{2}})(p^{n_{1}}-p^{m_{1}}),\\
&\vdots\\
\alpha_{n_{1}n_{2}+n_{1}+2}=&p^{n_{1}}q+\phi(p^{n_{1}})(q^{n_{2}}-q)-\phi(p^{n_{1}}q),\\
&\vdots\\
\alpha_{n_{1}n_{2}+n_{1}+m_{1}+1}=&p^{n_{1}}q^{m_{2}}+\phi(p^{n_{1}})(q^{n_{2}}-q^{m_{2}})-\phi(p^{n_{1}}q^{m_{2}}),\\
&\vdots\\
\alpha_{(n_{1}+1)n_{2}+n_{1}}=&p^{n_{1}}q^{n_{2}-1}+\phi(n)-\phi(p^{n_{1}}q^{n_{2}-1}).
\end{align*}
Also, value of $r_{i}+\alpha_{i}= r_{i}^{'} $'s are given by 
\begin{align*}
r_{1}^{'}=&n-1,~ r_{2}^{'}=n-q^{n_{2}},\\
&\vdots\\
r_{m_{1}+1}^{'}=&p^{m_{1}-1}+q^{n_{2}}(p^{n_{1}}-p^{m_{1}-1})-1,\\
&\vdots\\
r_{n_{1}+1}^{'}=&p^{n_{1}-1}+\phi(p^{n_{1}})q^{n_{2}}-1, ~r_{n_{1}+2}^{'}=n-p^{n_{1}},\\
&\vdots\\
r_{n_{1}+m_{1}+1}^{'}=&q^{m_{2}-1}+p^{n_{1}}(q^{n_{2}}-q^{m_{2}-1})-1,\\
&\vdots\\
r_{n_{1}+n_{2}+1}^{'}=&q^{n_{2}-1}+\phi(q^{n_{2}})p^{n_{1}}-1,~r_{n_{1}+n_{2}+2}^{'}=\phi(p)+\phi(q)+(q^{n_{2}}-1)(p^{n_{1}}-1),\\
&\vdots\\
r_{n_{1}+n_{2}+m_{1}+1}^{'}=&q^{n_{2}}(p^{n_{1}}-1)+q^{m_{2}}-q^{m_{2}-1}(p^{n_{1}}-p)-1,\\
&\vdots\\
r_{2n_{1}+n_{2}+1}^{'}=&pq^{n_{2}}+\phi(p^{n_{1}})(q^{n_{2}}-q)-1,\\
&\vdots\\
r_{m_{1}n_{2}+n_{1}+2}^{'}=&p^{m_{1}}+p^{n_{1}}(q^{n_{2}}-1) -p^{m_{1}-1} (q^{n_{1}}-q)-1,\\
&\vdots\\
r_{m_{1}N_{2}+N_{1}+m_{1}+1}^{'} =&n+p^{m_{1}-1}q^{m_{2}-1}+p^{m_{1}}q^{m_{2}}-\phi(p^{m_{1}}q^{m_{2}})-q^{n_{2}}p^{m_{1}-1}-p^{n_{2}}q^{m_{2}-1}-1,\\
&\vdots\\
r_{(m_{1}+1)n_{2}+n_{1}+1}^{'}=&p^{m_{1}}q^{n_{2}}+\phi(q^{n_{2}})(p^{n_{1}}-p^{m_{1}})-1,\\
&\vdots\\
r_{n_{1}n_{2}+n_{1}+2}^{'}=&p^{n_{1}}q+\phi(p^{n_{1}})(q^{n_{2}}-q)-1,\\
&\vdots\\
r_{n_{1}n_{2}+n_{1}+m_{1}+1}^{'}
=&p^{n_{1}}q^{m_{2}}+\phi(p^{n_{1}})(q^{n_{2}}-q^{m_{2}})-1,\\
&\vdots\\
r_{n_{1}n_{2}+n_{1}+m_{1}+1}^{'} =&p^{n_{1}}q^{n_{2}-1}+\phi(n)-1.
\end{align*}
By  Theorem \ref{mul of alpha n/n-1 of Zn}, $ \dfrac{n}{n-1} $ is the normalized Laplacian eigenvalue of $ \mathcal{P}(\mathbb{Z}_{n}) $ with multiplicity $ \phi(n). $  We note that each of $ G_{i}=K_{i} $ and by Theorem \ref{zn}, we see that 
\[ \dfrac{\phi(d_{1})+\alpha_{2}}{\phi(d_{1})+\alpha_{2}-1}=\dfrac{\phi(p)+n-\phi(p)-q^{n_{2}}+1}{\phi(p)+n-\phi(p)-q^{n_{2}}+1-1}=\dfrac{n-q^{n_{2}}+1}{n-q^{n_{2}}}, \] is the normalized Laplacian eigenvalue of $ \mathcal{P}(\mathbb{Z}_{n}) $ with multiplicity $ \phi(p)-1 $. Similarly, other normalized Laplacian eigenvalues of $ \mathcal{P}(\mathbb{Z}_{n}) $ can be found. Also, substituting the values of $ \alpha_{i} $'s, $ r_{i}^{'} $'s and using the adjacency relations, the remaining normalized Laplacian eigenvalues of $ \mathcal{P}(\mathbb{Z}_{n}) $ are the eigenvalues of matrix \eqref{quotient matrix of Z_n}.\qed

\chapter{Distance signless Laplacian eigenvalues  and generalized distance energy of graphs}

In this chapter, we obtain the distance signless Laplacian eigenvalues of the joined union of regular graphs. We also find the distance signless Laplacian eigenvalues of the zero-divisor graphs associated to ring $ \mathbb{Z}_{n} $. Also, we define the generalized distance energy and obtain the bounds for it. We characterize the graphs attaining these bounds. We show that the complete bipartite graph and the star graph have the minimum generalized distance energy among the class of bipartite graphs and trees respectively.

\section{ Introduction}
\indent In a graph $G$, the \textit{distance} between any two vertices $u,v\in V(G),$ denoted by $d(u,v)$, is defined as the length of a shortest path between $u$ and $v$. The \textit{distance matrix} of $G$, denoted by $\mathcal{D}(G)$, is defined as $\mathcal{D}(G)=(d_{uv})$, where $u,v\in V(G)$. More about distance matrix can be found in the survey paper \cite{aouchiche2014}. The \textit{transmission degree} $Tr_{G}(v)$ of a vertex $v$ is defined to be the sum of the distances from $v$ to all other vertices in $G$, that is, $Tr_{G}(v)=\sum\limits_{u\in V(G)}d(u,v)$. A graph $ G $ is said to be $k$-\textit{transmission regular} if $ Tr_{G}(v)=k,$ for each $ v\in V(G). $  The \textit{transmission} number or Wiener index of a graph $ G, $ denoted by $W(G), $ is the sum of distances between all unordered pairs of vertices in $G$.
It is clear that $W(G)=\frac{1}{2}\displaystyle\sum_{v\in V(G)}Tr_{G}(v)$. If $Tr_G(v_i)$ (or simply $Tr_{i}$) is the transmission degree of the vertex $v_i\in V(G)$, the sequence $\{Tr_{1},Tr_{2},\ldots,Tr_{n}\}$ is called the \textit{transmission degree sequence} of  $G$.\\
\indent Let $Tr(G)=diag (Tr_1,Tr_2,\ldots,Tr_n) $ be the diagonal matrix of vertex transmissions of $G$. Aouchiche and Hansen \cite{AH2} introduced the Laplacian and the signless Laplacian for the distance matrix of a connected graph. The matrix $ \mathcal{D}^L(G)=Tr(G)-\mathcal{D}(G) $ is called the \textit{distance Laplacian matrix} of $G$, while the matrix  $\mathcal{D}^{Q}(G)=Tr(G)+\mathcal{D}(G)$ is called the \textit{distance signless Laplacian matrix of $G$.} The matrices $ \mathcal{D}(G) $ and $\mathcal{D}^{Q}(G)$ are real symmetric and positive definite for $ n\geq 3 $, so their eigenvalues can be arranged as $ \rho^{\mathcal{D}}_n(G)\leq \ldots\leq\rho^{\mathcal{D}}_2(G)\leq \rho^{\mathcal{D}}_1(G) $ and  $\rho^{Q}_n(G)\leq \ldots\leq\rho^{Q}_2(G)\leq \rho^{Q}_1(G)$, respectively.  The eigenvalues $ \rho^{\mathcal{D}}_1(G) $ and $ \rho^{Q}_1(G) $ are called the distance  and distance signless Laplacian spectral radii of $ G. $ Since the matrix $\mathcal{D}^{Q}(G)$ is non-negative and irreducible, so by \emph{Perron-Frobenius Theorem}, $ \rho^{Q}_1(G) $ is positive and simple zero of the characteristic polynomial of $\mathcal{D}^{Q}(G)$. There corresponds a unique eigenvector with positive entries having unit length, which is called the distance signless Laplacian \emph{Perron vector} of graph $ G. $ More work on distance Laplacian and distance signless Laplacian matrix can be found in \cite{AH3,AH4,dah,hilal, HBRP,lwm,ld,bilal1,bilalkorea} and references therein.
\section{Distance signless Laplacian eigenvalues of the joined union of regular graphs}

We compute the distance signless Laplacian spectrum of the joined union of graphs $G_1,G_2,\dots,G_n$ in terms of the adjacency spectrum of the graphs $G_1,G_2,\dots,G_n$ and the eigenvalues of quotient matrix.

\begin{theorem} \label{thm 22}
Let $G$ be a graph of order $n$ having vertex set $V(G)=\{ v_{1}, \ldots, v_{n}\}.$ Let $ G_{i}$ be $r_{i}$-regular graphs of order $ n_{i} $ having adjacency eigenvalues $\lambda_{i1}=r_{i}\geq \lambda_{i2}\geq\ldots\geq \lambda_{in_{i}}, $ where $ i=1,2, \ldots, n$.  The distance signless Laplacian spectrum of the joined union graph $ G[G_{1},\ldots, G_{n}] $ of order $ \sum\limits_{i=1}^{n}n_{i} $ consists of the eigenvalues $ 2n_i+n_{i}^{\prime}-r_{i}-\lambda_{ik}-4$ for $ i=1,\ldots,n $ and $ k=2,3,\ldots, n_{i} $, where $n_{i}^{\prime}=\sum\limits_{k=1,k\ne i}^{n}n_kd_{G}(v_i,v_k)$. The remaining $n$ eigenvalues are given by the equitable quotient matrix
\begin{equation*}
\mathcal{Q}=\begin{pmatrix}
4n_{1}+n_{1}^{'}-2r_{1}-4& n_{2}d_{G}(v_{1}, v_{2})&\ldots & n_{n}d_{G}(v_{1},v_{n})\\
n_{1}d_{G}(v_{2}, v_{1})&4n_{2}+n_{2}^{'}-2r_{2}-4& \ldots & n_{n}d_{G}(v_{2}, v_{n})\\
\vdots &\vdots  &\ddots &\vdots\\
n_{1}d_{G}(v_{n}, v_{1})& n_{2}d_{G}(v_{n}, v_{2})&\ldots &4n_{n}+n_{n}^{'}-2r_{n}-4
\end{pmatrix}.
\end{equation*}
\end{theorem}
\noindent{\bf Proof.}
Let $G$ be a connected  graph of order $n$ having vertex set $V(G)=\{ v_{1}, \ldots, v_{n}\}$ and let $V(G_i)=\{ v_{i1}, \ldots, v_{in_i}\}$ be the vertex set of the graph $G_i$, for $i=1,2,\dots,n$. Let  $H=G[G_{1}, \ldots, G_{n}]$ be the joined union  of the graphs $G_1,G_2,\dots, G_n$. The transmission degree of each vertex $v_{ij}\in V(H)$, for $1\le i\leq n$ and $1\leq j\leq n_i$, is given by
\begin{equation}\label{transmission}
Tr(v_{ij})=2n_i+n_{i}^{'}-2-d(v_{ij})
\end{equation}
where $n_{i}^{'}=\sum\limits_{k=1}^{n}n_kd_{G}(v_i,v_k)$ is same for every vertex belonging to $ V(G_i)$.
We note that
$ Tr(v_{i,j})$ is equal to the distance of $j$th vertex in $G_{i} $ to all other vertices of $ G$, which is further equal to the distance within vertices of $ G_{i} $ and the distance of vertices outside of  $ G_{i} $. So
\begin{align*}
Tr(v_{i,j})&= d(v_{i,j})+2(n_{i}-1-d(v_{i,j})) + n_{1}d(v_{i},v_{1})+n_{2}d(v_{i},v_{2})+\dots+n_{n}d(v_{i},v_{n})\\&=2n_{i}-2-d(v_{i,j})+n_{i}^{'},
\end{align*}
and thus Equation (\ref{transmission}) is true for any diameter.
Under the suitably labelling of the vertices of $H$, the distance signless Laplacian matrix $D^{Q}(H)$ can be written as
\begin{equation*}
\mathcal{D}^{Q}(H)=\begin{pmatrix}
h_1&d_{G}(v_{1}, v_{2})J_{n_1\times n_2}& \ldots & d_{G}(v_{1}, v_{n})J_{n_1\times n_n} \\
d_{G}(v_{2}, v_{1})J_{n_2\times n_1}&h_2& \ldots & d_{G}(v_{2}, v_{n})J_{n_2\times n_n}\\
\vdots &\vdots &\ddots &\vdots\\
d_{G}(v_{n}, v_{1})J_{n_n\times n_1}& d_{G}(v_{n}, v_{2})J_{n_n\times n_{n-1}}&\ldots &h_n
\end{pmatrix},
\end{equation*}
where,  for $i=1,2,\ldots,n,$
\begin{equation*}
h_{i}=(2n_{i}+n_{i}^{'}-2-d(v_{ij}))I_{n_{i}}+2(J_{n_{i}}-I_{n_{i}})-A(G_{i}).
\end{equation*}
in which, for $ i=1,2,\ldots,n,$
 $A(G_i)$ is the adjacency matrix of graph $ G_{i}$, and $J_{n_i}$ is the matrix having all entries $ 1 $ and $ I_{n_i}$ is the identity matrix of order $n_i$. \\
\indent As $ G_{i} $ is an $ r_{i} $ regular graph, so $e_{n_i}=(1,1,\dots,1)^T$ with $ n_{i} $ entries is the eigenvector of the adjacency matrix $ A(G_{i}) $ corresponding to the eigenvalue $r_{i} $ and all other eigenvectors are orthogonal to $e_{n_i}.$ Let $ \lambda_{ik}$, $2\leq k\leq n_i$, be any eigenvalue of $ A(G_{i})$ with the corresponding eigenvector $X=(x_{i1},x_{i2},\dots,x_{in_i})^T$ satisfying $e_{n_i}^TX=0.$ It is a well known fact that the column vector $X$ can be regarded as a function defined on $ V(G_i) $ relating vertex $ v_{ij} $ to $ x_{ij} $, that is, $ X(v_{ij})=x_{ij} $ for $ i=1,2,\ldots,n $ and $j=1,2,\dots,n_i$. Now, consider the vector $Y=(y_{1},y_{2},\dots,y_{n})^T$, where
\begin{equation*}
y_{j}=\left\{
\begin{array}{rl}
x_{ij}&~\text{if}~ v_{ij}\in V(G_i)\\
0&~\text{otherwise.}\\
\end{array}\right.
\end{equation*}

Noting that $ e_{n_{i}}^{T}X=0 $ and coordinates of the  vector $Y$ corresponding to vertices in $\cup_{j\ne i}V_j$ of $H$ are zeros, therefore, we have
\begin{align*}
D^{Q}(H)Y=\begin{pmatrix}
0\\
\vdots\\
0\\
(2n_i+n_{i}^{'}-d(v_{ij})-4)X-\lambda_{ik} X\\
0\\
\vdots\\
0
\end{pmatrix}=(2n_i+n_{i}^{'}-d(v_{ij})-\lambda_{ik}-4)Y.
\end{align*}

This shows that the vector $Y$ is an eigenvector of $\mathcal{D}^{Q}(H) $ corresponding to the eigenvalue $2n_i+n_{i}^{\prime}-d(v_{ij})-\lambda_{ik}-4$, for every eigenvalue $\lambda_{ik}$, $2\leq k\leq n_i$, of $A(G_i)$. It follows that $2n_i+n_{i}^{\prime}-r_{i}-\lambda_{ik}-4$, for $1\leq i\leq n$ and $2\leq k\leq n_i$, is an eigenvalue of $ \mathcal{D}^{Q}(H) $ and in this way we have obtained $ \sum\limits_{i=1}^{n}n_{i}-n$  eigenvalues. The remaining $ n $ eigenvalues are the zeros of the characteristic polynomial of the following equitable quotient matrix.
\begin{align*}
\begin{pmatrix}
4n_{1}+n_{1}^{'}-2r_{1}-4& n_{2}d_{G}(v_{1}, v_{2})&\ldots & n_{n}d_{G}(v_{1},v_{n})\\
n_{1}d_{G}(v_{2}, v_{1})&4n_{2}+n_{2}^{'}-2r_{2}-4& \ldots & n_{n}d_{G}(v_{2}, v_{n})\\
\vdots &\vdots  &\ddots &\vdots\\
n_{1}d_{G}(v_{n}, v_{1})& n_{2}d_{G}(v_{n}, v_{2})&\ldots &4n_{n}+n_{n}^{'}-2r_{n}-4
\end{pmatrix}.
\end{align*} This completes the proof.\qed

The following consequence of Theorem \ref{thm 22} gives the distance signless Laplacian spectrum of the join of two regular graphs, in which one of the graphs is the union of two regular graphs having distinct vertex sets.

\begin{corollary}\label{th2}
For $ i=1,2,3$, let $ G_{i} $ be $ r_{i} $ regular graphs of orders $ n_{i} $ having adjacency eigenvalues $ \lambda_{i,1}=r_{i}\geq \lambda_{i,2}\geq \dots\geq\lambda_{i,n_{i}} $.  Then the distance signless Laplacian spectrum of $ G_{1}\triangledown (G_{2}\cup G_{3}) $ is
\begin{align*}
\Big\{ (n+n_{1}-r_{1}-\lambda_{1,k}-4)^{[n_{1}-1]}, (2n-n_{1}-&r_{2}-\lambda_{2,k}-4)^{[n_{2}-1]},\\ &(2n-n_{1}-r_{3}-\lambda_{3,k}-4)^{[n_{3}-1]} \Big\},
\end{align*}
where $ k=2,3,\dots,n_{i} $, for $ i=1,2,3 $ and $ n=n_{1}+n_{2}+n_{3} .$ The remaining three eigenvalues are given by the equitable quotient matrix
\begin{equation*}
\begin{pmatrix}
n+3n_{1}-2r_{1}-4 &n_{2} &n_{3}\\
n_{1} & 2n+2n_{2}-n_{1}-2r_{2}-4 &2n_{3}\\
n_{1} &2n_{2} & 2n+2n_{3}-n_{1}-2r_{3}-4
\end{pmatrix}.
\end{equation*}

\end{corollary}

The next consequence gives the distance signless Laplacian spectrum of the join of two regular graphs.

\begin{corollary}\label{th21}
Let $ G_{1} $ and $ G_{2} $ be $ r_{1} $ and $ r_{2} $ regular graphs of order $ n_{1} $ and $ n_{2} $, respectively. Let $ \lambda_{1}=r_{1},\lambda_{2}, \dots, \lambda_{n_{1}} $ and $ \lambda^{'}_{1}=r_{2},\lambda^{'}_{2},\dots,\lambda^{'}_{n_2} $ be the adjacency eigenvalues of $ G_{1} $ and $ G_{2} $, respectively. Then the distance signless Laplacian spectrum of $ G_{1}\triangledown G_{2} $ consists of eigenvalues $ 2n-n_{2}-r_{1}-4-\lambda_i $ and $ 2n-n_{1}-r_{2}-4-\lambda^{'}_j $ with multiplicities $ n_{1}-1 $ and $ n_{2}-1 $ respectively, where $n=n_1+n_2,~2\leq i\leq n_1$ and $2\leq j\leq n_2$. The remaining two eigenvalues are given by the eigenvalues of the equitable quotient matrix
\begin{equation} \label{eq21}
\begin{pmatrix}
n+3n_{1}-2r_{1}-4& n_{2}\\
n_{1}& n+3n_{2}-2r_{2}-4

\end{pmatrix}.
\end{equation}
\end{corollary}

The following are immediate consequences of Corollary \ref{th2} and give the distance signless Laplacian eigenvalues of a complete bipartite and a complete split graph.
\begin{corollary}\label{kab}
The distance signless Laplacian spectrum of $ K_{a,b} $ consists of the eigenvalue $ 2a+b-4 $ with multiplicity $ a-1$, the eigenvalue $2b+a-4$ with multiplicity $ b-1 $ and the eigenvalues
\begin{equation*}
\dfrac{1}{2}\left ( 5(a+b)-8\pm \sqrt{9(a-b)^2+4ab} \right ).
\end{equation*}
\end{corollary}
\begin{corollary}\label{cp}
The distance signless Laplacian eigenvalues of $CS_{\omega,n-\omega}$ are given by \[ \left\lbrace  (n-2)^{[\omega-1]}, (2n-\omega-4)^{[n-\omega-1]}, \frac{1}{2}\left(5n-2\omega-6\pm \sqrt{D}\right)\right\rbrace,  \]
where $ (3(2\omega-n)-2(\omega-1))^{2}+4\omega(n-\omega). $
\end{corollary}
\section{Distance signless Laplacian eigenvalues of zero-divisor graphs associate to the ring $ \mathbb{Z}_{n} $ }

For a commutative ring $R$ with identity denoted by $1$, let the set $Z(R)$ denote the set of zero-divisors and let $Z^{*}(R)=Z(R)\setminus \{0\}$ be the set of non-zero zero-divisors of $R$.  The zero-divisor graph of $R$, denoted by $\Gamma(R)$, is a simple graph whose vertex set is $Z^{*}(R)$ and two vertices $u, v \in Z^*(R)$ are adjacent if and only if $uv=vu=0$. 
The zero-divisor graphs of commutative rings were first introduced by Beck \cite{ib}. In the definition he included the additive identity and was interested mainly in coloring of commutative rings. Later Anderson and Livingston \cite{al} modified the definition of zero-divisor graphs and excluded the additive identity of the ring in the zero-divisor set. We denote the ring of integers modulo $n$ by $\mathbb{Z}_n$. The order of the zero-divisor graph $\Gamma(\mathbb{Z}_n)$ is $ n-\phi(n)-1 $, where $ \phi $ is Euler's totient function. The zero-divisor graphs have been extensively studied in \cite{akbari, al,ib,  ajaz}. The spectral properties of adjacency and Laplacian  matrices for zero-divisor graphs were discussed in \cite{sc,magi, my}.

For integers $d$ and $n$ with $ 1<d<n$, if $ d $ is a proper divisor of $ n $, we write $ d|n $. Let $ d_{1}, d_{2},\dots,d_{t} $ be the distinct proper divisors of $ n$. Let $ \Upsilon_{n} $ be the simple graph with vertex set $ \{d_{1}, d_{2},\dots,d_{t}\} $, in which two distinct vertices are adjacent if and only if $ n|d_{i}d_{j}$. Let the prime factorization of $n$ be $ n=p_{1}^{n_{1}}p_{2}^{n_{2}}\dots p_{r}^{n_{r}} $, where $ r,n_{1},n_{2},\dots,n_{r} $ are positive integers and $ p_{1},p_{2},\dots,p_{r} $ are distinct prime numbers. The order of $ \Upsilon_{n} $ is given by
\begin{equation*}
|V( \Upsilon_{n} )|=\prod\limits_{i=1}^{r}(n_{i}+1)-2.
\end{equation*}

This $  \Upsilon_{n}  $ is connected (see \cite{sc}) and plays a fundamental role in the present section.
For $ 1\leq i \leq t $, we consider the sets $ A_{d_{i}}= \{ x\in \mathbb{Z}_{n} : (x,n)=d_{i} \}$, where $(x,n)$ denotes the greatest common divisor of $x$ and $n$. We see that $ A_{d_{i}} \cap A_{d_{j}}=\phi$, when $ i\neq j $, implying that the sets $ A_{d_{1}}, A_{d_{2}}, \dots, A_{d_{t}} $ are pairwise disjoint and partitions the vertex set of $ \Gamma(\mathbb{Z}_{n}) $ as $V(\Gamma(\mathbb{Z}_{n}))=A_{d_{1}}\cup A_{d_{2}}\cup \dots \cup A_{d_{t}}.$ From the definition of $ A_{d_{i}} $, a vertex of $ A_{d_{i}} $ is adjacent to the vertex of $ A_{d_{j}} $ in $ \Gamma(\mathbb{Z}_{n}) $ if and only if $ n $ divides $ d_{i}d_{j} $ , for $ i,j\in \{1,2,\dots,t\} $ (see \cite{sc}).

\noindent We have the following observation.

\begin{lemma}\cite{my}\label{young} If $ d_{i} $ is the proper divisor of $ n $, then
$|A_{d_{i}}|=\phi\left( \frac{n}{d_{i}}\right)$, for $ 1\leq i \leq t .$
\end{lemma}

The induced subgraphs $ \Gamma(A_{d_i}) $ of $ \Gamma(\mathbb{Z}_{n}) $ are either cliques or their complements as can be seen in the following lemma.

\begin{lemma}\cite{sc}\label{sc1} Let $ n $ be a positive integer and $ d_{i} $ be its proper divisor. Then following hold.
\begin{itemize}
\item[\bf(i)] For $i\in\{1,2,\dots,t\} , $ the induced subgraph $ \Gamma(A_{d_i}) $ of $ \Gamma(\mathbb{Z}_{n}) $ on the vertex set $ A_{d_{i}} $ is either the complete graph $ K_{\phi\left (\frac{n}{d_{i}}\right )} $ or its complement $\overline{ K}_{\phi\left (\frac{n}{d_{i}}\right )} .$ Indeed, $ \Gamma(A_{d_i}) $ is $ K_{\phi\left (\frac{n}{d_{i}}\right )} $ if and only $ n $ divides $ d_{i}^{2} .$
\item[\bf(ii)] For $i,j\in\{1,2,\dots,t\} $ with $ i\neq j ,$ a vertex of $ A_{d_{i}} $ is adjacent to either all or none of the vertices in $ A_{d_{j}} $ of $ \Gamma(\mathbb{Z}_{n}) $.

\end{itemize}
\end{lemma}

The following lemma says that $ \Gamma(\mathbb{Z}_{n}) $ is the joined union of certain complete graphs and null graphs.
\begin{lemma}\label{sc2} \cite{sc} Let $ \Gamma(A_{d_i}) $ be the induced subgraph of $ \Gamma(\mathbb{Z}_{n}) $ on the vertex set $ A_{d_i} $ for $ 1\leq i\leq t .$ Then $ \Gamma(\mathbb{Z}_{n}) = \Upsilon_{n}[\Gamma(A_{d_1}),\Gamma(A_{d_2}),\dots,\Gamma(A_{d_t})]. $
\end{lemma}

In the following example, we illustrate the structure of $ \Gamma(\mathbb{Z}_{30}) $ and find its distance signless Laplacian eigenvalues.\\
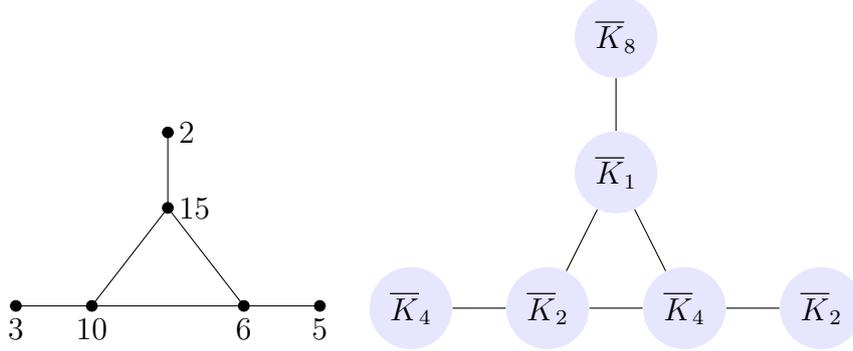
\begin{figure}\label{zero div graph of Z30}
\centering

\begin{tikzpicture}
\draw[fill=black] (0,0) circle (2pt);\draw[fill=black] (2,0) circle (2pt);\draw[fill=black] (1,1.3) circle (2pt);
\draw[fill=black] (1,2.3) circle (2pt);\draw[fill=black] (-1,0) circle (2pt);\draw[fill=black] (3,0) circle (2pt);

\draw[thin] (0,0)--(2,0); \draw[thin] (0,0)--(1,1.3);\draw[thin] (2,0)--(1,1.3);
\draw[thin] (0,0)--(-1,0); \draw[thin] (2,0)--(3,0);\draw[thin] (1,1.3)--(1,2.3);

\node at (-1,-.3) {$ 3 $};\node at (0,-.3) {$ 10 $};
\node at (1.35,1.3) {$ 15 $};\node at (1.25,2.3) {$ 2 $};
\node at (3,-.3) {$ 5 $};\node at (2,-.3) {$ 6 $};

\end{tikzpicture}\quad
\begin{tikzpicture}
[scale=.9, auto=left, every node/.style={circle, fill=blue!10, minimum width=0.1cm}]
\node (n1) at (0,0) {$ \overline{K}_{2} $};\node (1) at (-2,0) {$ \overline{K}_{4} $};
\node (n2) at (2,0) {$ \overline{K}_{4} $};\node (4) at (4,0) {$ \overline{K}_{2} $};
\node (n3) at (1,2) {$ \overline{K}_{1} $};\node (3) at (1,4) {$ \overline{K}_{8} $};

\foreach \from/\to in {n1/n2,n1/n3,n2/n3, n2/4, n1/1, n3/3}
\draw(\from)--(\to);
\end{tikzpicture}

\caption{Proper divisor graph $ \Upsilon_{30} $ and zero divisor graph $  \Gamma(\mathbb{Z}_{30})  $. }
\end{figure}
\noindent\textbf{Example} Let $ n=30 $. Then $ 2,3,5,6,10 $ and $ 15 $ are the proper divisors of $ n $ and $ \Upsilon_{n} $ is the graph $ G_{6} : 3\sim 10\sim 6 \sim 5, 10\sim 15\sim 2 $ and $ 6\sim 15 $, that is $ \Upsilon_{n} $ is the triangle graph having a pendent edge at each vertex of the triangle. Applying Lemma \ref{sc2}, we have
\begin{equation*}
\Gamma(\mathbb{Z}_{30})=\Upsilon_{30}[\overline{K}_{4},\overline{K}_{2},\overline{K}_{4},\overline{K}_{2},K_{1},\overline{K}_{8}].
\end{equation*}
By Theorem \ref{thm 22}, the distance signless Laplacian spectrum of $ \Gamma(\mathbb{Z}_{30}) $ consists of the eigenvalues $ \{46^{[4]}, 43^{[7]}, 33^{[3]}, 29\} $, and the remaining eigenvalues are given by

\begin{equation*}
\begin{pmatrix}
54 & 2 & 8 & 6 & 2 & 24\\
4 & 33 & 4 & 4 & 1 & 16\\
8 & 2 & 41 & 2 & 1 & 16\\
8 & 4 & 4 & 50 & 2 & 24\\
8 & 2 & 4 & 4& 26 & 8\\
12 & 4 & 8 & 6 & 1 & 59

\end{pmatrix}.
\end{equation*}
The characteristic polynomial of the above matrix is $ x^6-263x^{5}+27575x^{4}-1489941x^{3}+44016528x^{2}-676666908x+4239569664 $
 and it's approximated zeros are \begin{equation*} \{87.1555, 44.5461, 40.5727, 35.0098, 30.5597, 25.1562\}. \end{equation*}

Next we will discuss the distance signless laplacian spectrum of $ \Gamma(\mathbb{Z}_n) $ for $ n\in\{pq,p^{2}q,p^3,p^4\} $ with the help of Theorem \ref{thm 22}. Let $ n=pq $, where $ p $  and $ q $, $ p<q $ are distinct primes. Then  by Lemma \ref{sc1} and Lemma \ref{sc2}, we have
\begin{equation}\label{eq pq}
\Gamma(\mathbb{Z}_{pq})= \Upsilon_{pq}[\Gamma(A_{p}),\Gamma(A_{q})] =K_{2}[\overline{K}_{\phi(p)},\overline{K}_{\phi(q)}]=\overline{K}_{\phi(p)}\triangledown \overline{K}_{\phi(q)}=K_{\phi(p),\phi(q)}.\end{equation}

The following result gives the distance signless Laplacian eigenvalues of $ \Gamma(\mathbb{Z}_{n}) $, when $ n $ is product of two primes.
\begin{lemma}\label{pq}
The distance signless Laplacian eigenvalues of $ \Gamma(\mathbb{Z}_{n}) $ are $ 2n-q-3 $ and $ 2n-p-3 $ with multiplicities $ p-2 $ and $ q-2 $ respectively, and
\begin{equation*}
\frac{1}{2}\left( 5n-8\pm\sqrt{9(p-q)^{2}+4(p-1)(q-1)}\right). \end{equation*}
\end{lemma}
\noindent\textbf{Proof.} Let $ n=pq $, where $ p<q $ are distinct primes. Then the proper divisors of $ n $ are $ p $ and $ q $, so $ \Upsilon_{pq} $ is $ K_{2} $. By Equation \eqref{eq pq} and Theorem \ref{thm 22}, the  result follows.\qed

Next, result gives the distance signless Laplacian eigenvalues of $\Gamma(\mathbb{Z}_{n}) $ for $ n=p^{2}q $.
\begin{lemma}
The distance signless Laplacian spectrum of $ \Gamma(\mathbb{Z}_{p^{2}q}) $ is \begin{align*}
 \{(2pq+3p^2-4p-5)^{[pq-p-q]},&(pq+2p^2-2p-3)^{[p-2]},(2pq+p^2-2p-4)^{[q-2]},\\& (3pq+2p^2-3p-2q-4)^{p^2-p-1},x_{1}\geq x_{2} \geq x_{3} \geq x_{4}\} \end{align*} where $ x_{1}, x_{2}, x_{3}, x_{4} $ are the zeros of the characteristic polynomial of the matrix.

\end{lemma}
\noindent\textbf{Proof.} Let $ n=p^{2}q, $ where $ p $ and $ q $ are distinct primes. Since proper divisors of $ n $ are $ p, q, pq, p^2 $, so $ \Upsilon_{p^{2}q} $ is the path $ P_{4} : q\sim p^2 \sim pq\sim p$. By Lemma \ref{sc2}, we have \begin{equation*}
 \Gamma(\mathbb{Z}_{p^2q})= \Upsilon_{p^2q}[\Gamma(A_{q}),\Gamma(A_{P^2}),\Gamma(A_{pq}),\Gamma(A_{p})] =P_{4}[\overline{K}_{\phi(p^2)},\overline{K}_{\phi(q)},K_{\phi(p)},\overline{K}_{\phi(pq)}] .\end{equation*}
Now, by Theorem \ref{thm 22}, clearly $ 2pq+3p^2-4p-5,pq+2p^2-2p-3,2pq+p^2-2p-4$ and $3pq+2p^2-3p-2q-4$ are the distance signless Laplacian eigenvalues with multiplicities $pq-p-q, p-2, q-2$ and $ p^2-p-1$ respectively. The remaining four distance signless Laplacian eigenvalues are given by the equitable quotient matrix
\begin{equation*}
\begin{pmatrix}
p(4p+3q-5)-2q-4& p-1&2(q-1)& 3p(p-1)\\
pq-p-q+1 & p(p+2q-2)+2q-7 & q-1 & 2p(p-1)\\
2(pq-p-q+1) & p-1& p(2p+q-1)-4 & p(p-1)\\
3(pq-p-q+1)& 2(p-1)& q-1&d_{4}

\end{pmatrix},
\end{equation*} 
where $ d_{4}=p(3p+4q-6)-2q-3 $\qed

A graph $G$ is called distance signless Laplacian integral graph if all its distance signless Laplacian eigenvalues are integers. In order to see when the joined graph or a zero divisor graph $ \Gamma(\mathbb{Z}_n) $ is distance signless Laplacian integral, we have the following proposition.
\begin{proposition}
The zero divisor graph $ \Gamma(\mathbb{Z}_n) $ is distance signless Laplacian integral if and only if matrix $ Q $ of Theorem \ref{thm 22} is integral.
\end{proposition}

The following result say that $ \Gamma(\mathbb{Z}_{n}) $ is distance singnless Laplacian integral only for prime power.
\begin{proposition}
Let $ n=p^2 $, where $ p $ is any prime. Then  $ \Gamma(\mathbb{Z}_{n}) $  is distance signless Laplacian integral and  spectrum is given by $ \{2p-4, (p-3)^{[p-2]}\} .$
\end{proposition}
\noindent\textbf{Proof.} Since $ \Gamma(\mathbb{Z}_{p^2}) =\Gamma(A_{p})$ is the complete graph $ K_{p-1} ,$ the result follows according to $ p=2 $ or $ p>2 .$ \qed

In the next result, we find distance signless Laplacian eigenvalues of the complete split graph  of $ \Gamma(\mathbb{Z}_{n}). $
\begin{proposition}

Let $n=p^3  $. Then the distance signless Laplacian eigenvalues of $ \Gamma(\mathbb{Z}_{n}) $ are $ 2n-p-3 $ and $ 2n-p^2-1 $ with multiplicities $ p^2-p-1 $ and $ p-2 $ respectively, and
\begin{equation*}
 \left\lbrace \frac{1}{2}\left(5n-2(p-1)-6\pm \sqrt{(3(2p-2-n)-2(p-2))^{2}+4(p-1)(n-p+1)}\right)\right\rbrace. \end{equation*}
\end{proposition}
\noindent\textbf{Proof.} Since proper divisors of $ n $ are $ p $ and $ p^2 $, so $ \Upsilon_{n} $ is $ K_{2}: $ $ p\sim p^2 $. By Lemma \ref{sc2}, we have

\begin{equation*}
 \Gamma(\mathbb{Z}_{p^3}) =\Upsilon_{p^3}[\Gamma(A_{p}),\Gamma(A_{p^2})] =K_{2}[\overline{K}_{\phi(p^2)},\overline{K}_{\phi(p)}]=\overline{K}_{p(p-1)}\triangledown K_{p-1} .\end{equation*}
This implies that $ \Gamma(\mathbb{Z}_{p^3}) $ is a complete split graph of order $ p^2-1 $, having independent set of cardinality $ p(p-1) $ and clique of order $ p-1 $. By Theorem \ref{thm 22} and Corollary \ref{cp}, the result follows.\qed

Next, proposition gives the distance signless Laplacian eigenvalues of $ \Gamma(\mathbb{Z}_{p^{4}}) $.

\begin{proposition}
The distance signless Laplacian eigenvalues of $ \Gamma(\mathbb{Z}_{n}) $, where $ n=p^4 $, are $ n-2 $, $ 2n-p-3$ and $ 2n-p^2-1 $ with multiplicities $ p-2 $,  $ p^2(p-1)-1 $ and $ p(p-1)-1 $ respectively and three more distance signless Laplacian eigenvalues given by
\begin{equation*}
 \begin{pmatrix}
n+p-3 & p^2(p-1) & p(p-1)\\
p-1 & 2n+p(2p^2-2p-1)-3 & 2p(p-1)\\
p-1 & 2p^2(p-1) & 2n-p-1
\end{pmatrix}. \end{equation*}
\end{proposition}
\noindent\textbf{Proof.}
As proper divisors of $ n $ are $ p, p^2 $ and $ p^3 $, so $ \Upsilon_{n} $ is $ P_{3}: $ $ p\sim p^3\sim p^2 $. By Lemmas \ref{young}, \ref{sc1} and \ref{sc2}, we have
$ \Gamma(A_{p})=\overline{K}_{\phi(p^3)}=\overline{K}_{p^2(p-1)}, \Gamma(A_{p^2})=K_{\phi(p^2)}=K_{p(p-1)}, $ and $ \Gamma(A_{p^3})=K_{\phi(p)}=K_{p-1}. $ Therefore
\begin{align*}
\Gamma(\mathbb{Z}_{p^4}) &=\Upsilon_{p^3}[\Gamma(A_{p}),\Gamma(A_{p^3}),\Gamma(A_{p^2})] =P_{3}[\overline{K}_{p^2(p-1)},K_{{p-1}},K_{p(p-1)}]\\
&=K_{p-1}\triangledown (\overline{K}_{p^2(p-1)}\cup K_{p(p-1)}) .\end{align*}
Thus, by Theorem \ref{thm 22}, the distance signless Laplacian  spectrum of $ \Gamma(\mathbb{Z}_{p^4}) $ is
\begin{equation*}
\left\lbrace (n-2)^{[p-2]},(2n-p-3)^{[p^2(p-1)-1]},(2n-p^2-1)^{[p(p-1)-1]},x_{1}\geq x_{2}\geq x_{3} \right\rbrace,
\end{equation*}
where $ x_{1}, x_{2}, x_{3} $ are the zeros of the polynomial $ x^3- x^2(2 p^3-2p^2-p + 5 n-7)+x( p^5-p^4-6np^3+12p^3-17p^2-8np+17p-8n^2+38n-44) - p^7+ 4 p^6 - 7 p^5 -p^4(3 n-13) + p^3(4 n^2-10 n ) -p^2( 4 n^2-13n+9)+ p( 4 n^2- 16 n+16) + 4 n^3 - 28n^2+64n-48
$.\qed

\section{ Spectral radius of generalized distance matrix of   graphs}
The \textit{generalized distance matrix} $\mathcal{D}_{\alpha}(G)$ of $ G $ introduced by Cui et al. \cite{CHT} as a convex combinations of $Tr(G)$ and $\mathcal{D}(G)$, and is defined as $\mathcal{D}_{\alpha}(G)=\alpha Tr(G)+(1-\alpha)\mathcal{D}(G)$,  for $0\leq \alpha\leq 1$. Clearly,
\[\mathcal{D}_{0}(G)=\mathcal{D}(G), \qquad 2\mathcal{D}_{\frac{1}{2}}(G)=\mathcal{D}^{Q}(G), \qquad \mathcal{D}_{1}(G)=Tr(G)\]
and  \[\mathcal{D}_{\alpha}(G)-\mathcal{D}_{\beta}(G)=(\alpha-\beta)\mathcal{D}^{L}(G).\]
The matrix $ \mathcal{D}_{\alpha}(G) $ reduces to merging the distance spectral, distance Laplacian spectral and distance signless Laplacian spectral theories. Since the matrix $ \mathcal{D}_{\alpha}(G)$ is real symmetric, all its eigenvalues are real. Therefore, we can arrange them  as $ \partial_{1}\geq \partial_{2}\geq \dots \geq \partial_{n}$. The largest eigenvalue  $ \partial_{1} $ of the matrix $\mathcal{D}_{\alpha}(G)$ is called the \textit{generalized distance spectral radius} of $ G$ ( we will denote $\partial_1(G)$ by $\partial(G)$). As $ \mathcal{D}_{\alpha}(G) $ is non-negative and irreducible, by the Perron-Frobenius theorem, $ \partial(G)$ is the unique eigenvalue and there is a unique positive unit eigenvector $X$ corresponding to $ \partial(G),$ which is called the \textit{generalized distance Perron vector} of $G.$ Further results on generalized distance matrix can be seen in \cite{CTZ, dpr,pham,pbhr,bilalbrazil,bilalalgebra}.

For $ M\in \mathbb{M}_{mn}(\mathbb{R}),$ the \emph{Frobenius norm} is defined as 
\[ \parallel M \parallel_F =\sqrt{\sum\limits_{i=1}^{n}\sum\limits_{j=1}^{n}|m_{ij}|^2}= \sqrt{trace(M^TM)} ,\] where the \emph{trace} of a square matrix is defined as sum of the diagonal entries. Further, if $ MM^T=M^TM $, then $ \parallel M \parallel_F^2=\sum\limits_{i=1}^{n}|\lambda_i(M)|^2 $, where $ \lambda_i's $ are the eigenvalues of $ M$. Also, $ \|M\|_{1}=|\lambda_{1}| $ is known as the spectral norm of $ M. $\\

We have the following observations.
\begin{lemma}
Let $G$ be a connected graph of order $n$ having transmission degrees $Tr_1,Tr_2,\dots,Tr_n$ and Wiener index $W(G)$. Then
\begin{itemize}
\item[{\bf(1)}] $\sum\limits_{i=1}^{n}\partial_i=2\alpha W(G).$
\item[{\bf(2)}] $\sum\limits_{i=1}^{n}\partial_i^2=\alpha^2 \sum\limits_{i=1}^{n}Tr_i^2+(1-\alpha)^2\parallel \mathcal{D}(G)\parallel_F^2$.
\item[{\bf(3)}] $\sum\limits_{i=1}^{n}\left (\partial_i-\frac{2\alpha W(G)}{n}\right )^2=\alpha^2\sum\limits_{i=1}^{n}Tr_i^2+(1-\alpha)^2\parallel \mathcal{D}(G)\parallel_F^2-\dfrac{4\alpha^2W^2(G)}{n}$.
\item[{\bf(4)}] $\partial(G)\geq \dfrac{2W(G)}{n}$, equality holds if and only if $G$ is transmission regular graph.
\item[{\bf(5)}] $ \partial(G)\geq \sqrt{\dfrac{\sum\limits_{i=1}^{n}Tr_i^2}{n}} $, equality holds if and only if $G$ is transmission regular graph.
\end{itemize}
\end{lemma}
\noindent\textbf{Proof}. {\bf(1)}. We have $\sum\limits_{i=1}^{n}\partial_i=\alpha \sum\limits_{i=1}^{n}Tr_i+(1-\alpha)\sum\limits_{i=1}^{n}\rho_i^\mathcal{D}=\alpha\sum\limits_{i=1}^{n}Tr_i=2\alpha W(G)$.\\
{\bf(2)}. We have
\begin{align*}
&\sum\limits_{i=1}^{n}\partial_i^2 =\alpha^2\sum\limits_{i=1}^{n}Tr_i^2+(1-\alpha)^2\sum\limits_{i=1}^{n}(\rho_i^\mathcal{D})^2  =\alpha^2\sum\limits_{i=1}^{n}Tr_i^2+(1-\alpha)^2\sum\limits_{i=1}^{n}\sum\limits_{j=1,j\neq i}^{n}(d_{ij})^2\\
& =\alpha^2\sum\limits_{i=1}^{n}Tr_i^2+(1-\alpha)^22\sum\limits_{i=1}^{n-1}\sum\limits_{j=i+1}^{n}(d_{ij})^2 =\alpha^2\sum\limits_{i=1}^{n}Tr_i^2+(1-\alpha)^22|\sum\limits_{i=1}^{n-1}\sum\limits_{j=i+1}^{n}(d_{ij})^2|\\
& =\alpha^2\sum\limits_{i=1}^{n}Tr_i^2+(1-\alpha)^2\parallel \mathcal{D}(G)\parallel_F^2
\end{align*}
{\bf(3)}. This follows as a consequence of \textbf{(2)}.\\
{\bf(4)}. From the Raleigh-Ritz's theorem for Hermitian matrices \cite{hj}, we have $$\partial(G)=max\{\textbf{X}^t\mathcal{D}_\alpha(G)\textbf{X}\}:\textbf{X}^t\textbf{X}=1\}.$$
Then, for the unit vector \textbf{U}, we have\[  \partial(G)\geq \dfrac{\textbf{U}^t\mathcal{D}_\alpha(G)\textbf{U}}{\textbf{U}^t\textbf{U}}=\dfrac{\alpha\sum\limits_{i=1}^{n}Tr_i+
(1-\alpha)\sum\limits_{i=1}^{n}Tr_i}{n}=\dfrac{\sum\limits_{i=1}^{n}Tr_i}{n}=\dfrac{2W(G)}{n}.\]
Suppose $G$ is $k$ transmission regular. Then each row sum of $\mathcal{D}_\alpha(G)$ sums to a constant $k$. Therefore, by the Perron-Frobenius theorem \cite{hj}, $k$ is simple and largest eigenvalue of $\mathcal{D}_\alpha(G)$, thus $\partial(G)=k$ and equality holds. Conversely, assume equality holds. Then $\partial(G )= \dfrac{\textbf{U}^t\mathcal{D}_\alpha(G)\textbf{U}}{\textbf{U}^t\textbf{U}}$. Hence $\mathcal{D}_\alpha(G)\textbf{U} =\partial(G)\textbf{U}$. Therefore $Tr_i=\partial(G)$ for all $i$ and thus $G$ is transmission regular.\\
{\bf(5)}. Proof is similar to that of {\bf(4)}. \qed \\

Let $ H $ be a graph with  transmission degree sequence $ t_{i} $ such that there is a vertex $v_{i}$ having transmission degree equal to $t_{1}=Tr_{\max}$, transmission degree of every neighbor of $v_{i}$ equal to $t_{2}$ and transmission degree of every vertex non-adjacent to $v_{i}$ equal to $Tr_{\min}$.

We consider a column vector $X=[x_{1},x_{2},\ldots,x_{n}]^{T}\in \mathbb{R}^{n}$ to be a function defined on $ V(G) $ which maps vertex $ v_{i} $ to $ x_{i} $, that is,
$ X(v_{i})=x_{i} $ for $i$, $1\leq i \leq n$. The quadratic form $X^{T}\mathcal{D}_{\alpha}(G)X$ can be described as
\begin{eqnarray*}
X^{T}\mathcal{D}_{\alpha}(G)X= \alpha\sum_{i=1}^{n}t(v_{i})x_{i}^{2}+2(1-\alpha)\sum_{1\leq i<j\leq n}d(v_{i},v_{j})x_{i}x_{j},
\end{eqnarray*}
and
\begin{eqnarray*}
X^{T}\mathcal{D}_{\alpha}(G)X= (2\alpha-1)\sum_{i=1}^{n}t(v_{i})x_{i}^{2}+(1-\alpha)\sum_{1\leq i<j\leq n}d(v_{i},v_{j})(x_{i}+x_{j})^2.
\end{eqnarray*}
Also,  $\partial$ is an eigenvalue of $\mathcal{D}_{\alpha}(G)$ corresponding to the eigenvector $X$ if and only if
$X\neq \textbf{0}$ and
\begin{eqnarray*}
\partial x_v=\alpha Tr(v_{i})x_{i}+(1-\alpha)\sum_{j=1}^{n}d(v_{i},v_{j})x_{j}.
\end{eqnarray*}
These equations are called the $(\partial, x)$-\textit{eigenequations} of $G$. For a normalized column vector $ X\in \mathbb{R}^{n} $ with at least one
non-negative component, by the Rayleigh's  principle, we have
$$\partial(G)\geq X^{T}\mathcal{D}_{\alpha}(G)X,$$
with equality if and only if $X$ is the generalized distance Perron vector of $G$.

The following result gives the upper bound for the largest eigenvalue of $ \mathcal{D}_{\alpha}(G). $
\begin{theorem}
If $\partial(G)$ is the spectral radius of $\mathcal{D}_{\alpha}(G)$, then
\begin{align}\label{zero}
\partial(G) \leq \max\limits_{1\leq i,~j\leq n}& \frac{1}{2}\Big[\alpha t_{i}+t_{j}-(1-\alpha)d_{ij}\\\nonumber
&+\sqrt{(\alpha t_{i}-t_{j})^{2}+(1-\alpha)(1-\alpha-2t_{j}-4t_{i}-2\alpha t_{i})d_{ij}}\Big],
\end{align}
where $t_{max}=t_{1}\geq t_{2}\geq \dots \geq t_{n}=t_{min}$ are the vertex transmission degrees of $G$. Further, equality holds if and only if $G$ is a transmission regular graph or $G$ is a semi-regular graph $G(n-1,s)$ or $G$ is in $H$, where $H$ is defined above.
\end{theorem}
\noindent\textbf{Proof.}  Corresponding to the eigenvalue $\partial$ of $\mathcal{D}_{\alpha}(G)$, let $X=[x_{1},x_{2},\ldots,x_{n}]^{T}$ be the eigenvector, with $x_{i}=1$ and $x_{t}\geq 1$ for all $t\neq i$. Let $x_{j}=\max \{x_{t}: ~1\leq t\leq n ~and~ t\neq i \}$. Therefore, from the $i$-th and $j$-th equations of $\partial~X=\mathcal{D}_{\alpha}(G)~X$, we have
\begin{equation}\label{1}
\partial~x_{i}=\alpha t_{i}x_{i}+(1-\alpha)\sum\limits_{t=1, t\neq i}^{n} d_{it}x_{t}
\end{equation}
and
\begin{equation}\label{2a}
\partial~x_{j}=\alpha t_{j}x_{j}+(1-\alpha)\sum\limits_{t=1, t\neq j}^{n} d_{jt}x_{t}.
\end{equation}
(\ref{1}) implies that $\partial \leq \alpha t_{i}+(1-\alpha)t_{i}x_{j},$ so that
\begin{equation}\label{3a}
\partial -\alpha t_{i}\leq (1-\alpha)t_{i}x_{j}
\end{equation}
Also, (\ref{2a}) implies that $\partial~x_{j}\leq \alpha t_{j}x_{j}+(1-\alpha)d_{ij}+ (1-\alpha)(t_{j}-d_{ij})x_{j}$,
which on simplification gives
\begin{equation}\label{4a}
\left[ \partial -t_{j}+(1-\alpha)d_{ij}\right]x_{j}\leq (1-\alpha)d_{ij}.
\end{equation}
Combining (\ref{3a}) and (\ref{4a}), we get
\begin{equation*}
(\partial-\alpha t_{i})[\partial-t_{j}+(1-\alpha)d_{ij}]x_{j}\leq [(1-\alpha)t_{i}x_{j}][(1-\alpha)d_{ij}].
\end{equation*}
Since $x_{j}>0$, therefore we have
\begin{equation*}
(\partial-\alpha t_{i})[\partial-t_{j}+(1-\alpha)d_{ij}]\leq (1-\alpha)^{2}t_{i}d_{ij},
\end{equation*}
which on simplification yields
\begin{equation}\label{5}
\partial^{2}-\left[\alpha t_{i}+t_{j}-(1-\alpha)d_{ij}\right]\partial+\alpha t_{i}t_{j}-(1-\alpha) t_{i}d_{ij}\leq 0.
\end{equation}
The solution of (\ref{5}) is given by
\begin{equation*}
\partial \leq \frac{1}{2}\left[\alpha t_{i}+t_{j}-(1-\alpha)d_{ij}+\sqrt{[\alpha t_{i}+t_{j}-(1-\alpha)d_{ij}]^{2}-4t_{i}\left[\alpha t_{j}+(1-\alpha)d_{ij}\right]}\right]
\end{equation*}
and
\begin{equation*}
\partial \geq \frac{1}{2}\left[\alpha t_{i}+t_{j}-(1-\alpha)d_{ij}-\sqrt{[\alpha t_{i}+t_{j}-(1-\alpha)d_{ij}]^{2}-4t_{i}\left[\alpha t_{j}+(1-\alpha)d_{ij}\right]}\right].
\end{equation*}
Therefore, in all cases, it follows that \begin{align*}
\partial(G)\leq \max\limits_{1\leq i,~j\leq n} &\frac{1}{2}\Big[\alpha t_{i}+t_{j}-(1-\alpha)d_{ij}\\
&+\sqrt{[\alpha t_{i}+t_{j}-(1-\alpha)d_{ij}]^{2}-4t_{i}[\alpha t_{j}+(1-\alpha)d_{ij}]}\Big].
\end{align*}

Now, we characterize the graphs for which equality holds in (\ref{zero}). In this regard, first assume that equality holds in (\ref{zero}). Therefore, equality cases in (\ref{3a}) and (\ref{4a}) imply that $x_{t}=x_{j}$ for all $t$, with $t\neq i$. As $x_{j}\leq x_{i}=1$, we have the following two possibilities to consider.\\
\noindent{\bf Case 1.} $x_{j}= x_{i}=1$. Here, $\partial =t_{i}$ for all $i$, $1\leq i\leq n$. So $G$ is a transmission regular graph.\\
\noindent{\bf Case 2.} $x_{j}< x_{i}=1$. Therefore, either $d_{i}=n-1$ or $d_{i}<n-1$. We look at these two cases separately as follows.\\
\noindent{\bf Case 2.1.} If $d_{i}=n-1$, then clearly $v_{t}\in N(v_{i})$ for all those $v_{t}$ which are in $V(G)-\{v_{i}\}$. Now, for any $v_{t}\in V(G)-\{v_{i}\}$, we have $\partial x_{j}=(t_{t}+\alpha -1) x_{j} +(1-\alpha)$. This clearly indicates that transmission degree of every vertex $v_{t}\in V(G)-\{v_{i}\}$ is equal to $t_{t}$, with of course $t_{t}> n-1$. So, evidently $t_{2}=t_{3}=\dots =t_{n}=t_{min}$. Further, $d_{i}=n-1$ implies that $G$ is of diameter 2 and thus $t_{i}=2n-2-d_{i}$. Therefore, the vertex degrees of $G$ are $n-1$, $d_{2}=d_{3}=\dots =d_{n}=s$ (say) with clearly $n-1>s$. Hence $G$ is isomorphic to $G(n-1,s)$.\\
\noindent{\bf Case 2.2.}  Now, for the case $d_{i}<n-1$, we consider the vertex partition of $V(G)-\{v_{i}\}$ as $U=N(v_{i})$ and $W=V(G)-(U\cup \{v_{i}\})$. Then
\begin{equation}\label{11a}
\partial =t_{i}\left[\alpha +(1-\alpha)x_{j}\right] ~~for~~v_{i},
\end{equation}
\begin{equation}\label{12}
\partial x_{j}=\left[t_{t}-(1-\alpha)\right]x_{j}+(1-\alpha) ~~for~~v_{t}\in N(v_{i}),
\end{equation}
\begin{equation}\label{13a}
\partial =t_{k}~~~~for~~~~v_{k}\notin \left[N(v_{i})\cup \{v_{i}\}\right].
\end{equation}
Simplifying (\ref{12}) implies that $t_{t}=\partial +(1-\alpha)(x_{j}-1)\frac{1}{x_{j}}$, for all $v_{t}\in N(v_{i})$. Therefore, transmission degree of every vertex in $N(v_{i})$ is equal to $t_{t}$, and transmission degree of every vertex outside $N(v_{i})\cup \{v_{i}\}$ is equal to $t_{i}$. Since $\alpha +(1-\alpha)x_{j}\leq 1$, from (\ref{11a}) and (\ref{13a}), we have $t_{k}>t_{i}$. In a similar way, $t_{k}>t_{t}$ follows from (\ref{12}) and (\ref{13a}). Thus, $t_{i}=t_{max}$, and $t_{k}=t_{2}$ for all $v_{k}\in W$. Also, $t_{t}=t_{min}$ for all vertices in $U$.\\
\indent From the arguments given above, we conclude that the given connected graph contains a vertex $v_{i}$ with transmission degree $t_{1}=t_{max}$, transmission degree of every neighbor of $v_{i}$ is equal to $t_{2}$ and transmission degree of every vertex non-adjacent to $v_{i}$ is equal to $t_{min}$. This is clearly the graph $H$ defined above.\\
\indent Conversely, for the regular graph, or the graph isomorphic to $G(n-1,s)$, or graph $H$, it is easy to verify that equality holds. \qed

\section{Genralized distance energy}
Let the auxiliary eigenvalues $ \Theta_i $, corresponding to the generalized distance eigenvalues be
\begin{equation*}
\Theta_i=\partial_i-\frac{2\alpha W(G)}{n}.
\end{equation*}
Now, we define the generalized distance energy of $G$ as the mean deviation of the  values of the generalized distance eigenvalues of $G$, that is,
\begin{equation*}
E^{\mathcal{D}_{\alpha}}(G)=\sum_{i=1}^{n}\left|\partial_i-\frac{2\alpha W(G)}{n}\right|=\displaystyle\sum_{i=1}^{n}|\Theta_{i}|.
\end{equation*}
We note that $ \sum\limits_{i=1}^{n}
|\Theta_{i}| $ is the \emph{trace norm} of the matrix $ \mathcal{D}_{\alpha}(G)-\frac{2W(G)}{n}I_{n} $, where $ I_{n} $ is the identity matrix. It is easy to see that $\sum_{i=1}^{n}\Theta_{i}=0.$  Clearly, $E^{\mathcal{D}_{0}}(G)=E^{\mathcal{D}}(G)$ and $2E^{\mathcal{D}_{\frac{1}{2}}}(G)=E^{Q}(G).$ This shows that the concept of generalized distance energy of $G$ merges the theories of distance energy and the distance signless Laplacian energy of $G$. Therefore, it will be interesting to study the quantity $E^{\mathcal{D}_{\alpha}}(G)$ and explore some properties like the bounds, the dependence on the structure of graph $G$ and the dependence on the parameter $\alpha$. 

The following lemma can be found in \cite{CHT}.

\begin{lemma}\label{lem31}
Let $G$ be a connected graph of order $n$ and $ \frac{1}{2}\leq \alpha\leq 1 .$ If $ G'$ is a connected graph obtained from $ G $ by deleting an edge, then  $ \partial_i(\mathcal{D}_\alpha(G'))\geq\partial_i(\mathcal{D}_\alpha(G))$, for all $ 1\leq i\leq n$.
\end{lemma}

The following result shows that the complete bipartite graph has the minimum generalized distance energy among all connected bipartite graphs of order $n$.
\begin{theorem}
Let $G$ be a connected  bipartite graph of order $n\geq 2$ having partite sets of cardinality $a$ and $n-a$, with $n-a\geq a,a\geq 2$ and let $\alpha\in [\frac{1}{2},1)$.
\begin{itemize}
\item[{\bf(1)}] If $2a<n< \frac{2a^2}{a-2}$ and $\alpha< \gamma$ or $\frac{2a^2}{a-2}\le n< 4.877 a$, then for $t=n-a$, we have
 \[ E^{\mathcal{D}_\alpha}(G)\geq (\alpha+2)n-4+\sqrt{\Delta}+2(n-a-1)(\alpha(2n-a)-2)-\frac{4\alpha(n-a) W(G)}{n}, \]
and for $t=n-a+1$, we have
\[ E^{\mathcal{D}_\alpha}(G)\geq 2(\alpha+2)n-8+2(n-a-1)(\alpha(2n-a)-2)-\frac{4\alpha(n-a+1) W(G)}{n}. \]
\item[{\bf(2)}] If $n\ge 4.877 a$, then \[ E^{\mathcal{D}_\alpha}(G)\geq (\alpha+2)n-4+\sqrt{\Delta}+2(n-a-1)(\alpha(2n-a)-2)-\frac{4\alpha(n-a) W(G)}{n}. \]
\item[{\bf(3)}] If $n=2a$, then \[ E^{\mathcal{D}_\alpha}(G)\geq 2(2a\alpha+4a-4)-\frac{4\alpha W(G)}{n}, \] where $\Delta=(2a^2+n^2-2na)(\alpha-2)^2+2(na-a^2)(\alpha^2-2)$, $\gamma=\frac{2n}{2n+an-2a^2}$ and $t$ is the number of generalized distance eigenvalues greater than or equal to $\frac{2\alpha W(K_{a,n-a})}{n}$.
Equality occurs in each case if and only if $G\cong K_{a,n-a}$.
\end{itemize}
\end{theorem}
\noindent\textbf{Proof.}
Let $ G $ be a connected bipartite graph of order $n$  with vertex set $V(G)$. Let $V(G)=V_1\cup V_2$, with $|V_1|=a, |V_2|=n-a$ and $n-a\geq a,~ a\geq 2$, be the bipartition of the vertex set  of $G$. Clearly $ G $ is a spanning subgraph of the complete bipartite graph $ K_{a,n-a}.$ Therefore, by Lemma \ref{lem31} with $ \frac{1}{2}\leq \alpha\leq 1 $, we have $ \partial_i(G)\geq\partial_i(K_{a,n-a})$, for all $i=1,2,\dots,n$.
 Let $ s$ be a positive integer such that $ \partial_s\geq\frac{2\alpha W(G)}{n}$ and $ \partial_{s+1}<\frac{2\alpha W(G)}{n} $. Using $\sum\limits_{i=1}^{n}\partial_i=2\alpha W(G)$ and the definition of generalized distance energy, we have
\begin{align*}
E^{\mathcal{D}_\alpha}(G)&=\sum\limits_{i=1}^{n}\left|\partial_i-\frac{2\alpha W(G)}{n}\right|\\&=\sum\limits_{i=1}^{s}\left( \partial_i-\frac{2\alpha W(G)}{n}\right) +\sum\limits_{i=s+1}^{n}\left( \frac{2\alpha W(G)}{n}-\partial_i\right)\\
&=2\left(\sum\limits_{i=1}^{s}\partial_i-\frac{2s\alpha W(G)}{n} \right).
\end{align*}
We first show that
\begin{equation}\label{def}
E^{\mathcal{D}_\alpha}(G)=2\left(\sum\limits_{i=1}^{s}\partial_i-\frac{2s\alpha W(G)}{n} \right)=2 \max_{1\leq j\leq n}\left\lbrace \sum\limits_{i=1}^{j}\partial_i-\frac{2\alpha jW(G) }{n}\right\rbrace.
\end{equation}
If $ j>s$, we have
\begin{align*}
 \sum\limits_{i=1}^{j}\partial_i -\frac{2\alpha jW(G)}{n} &
 =\sum\limits_{i=1}^{s}\partial_i+\sum\limits_{i=s+1}^{j}\partial_i-\frac{2\alpha jW(G)}{n} \\
 & < \sum\limits_{i=1}^{s}\partial_i-\frac{2\alpha sW(G)}{n}
 \end{align*}
as $\partial_i<\frac{2\alpha W(G)}{n}$, for $i\geq s+1$. Similarly, for $ j\leq s$, we have
\begin{align*}
\sum\limits_{i=1}^{j}\partial_i -\frac{2\alpha jW(G)}{n}\leq \sum\limits_{i=1}^{s}\partial_i-\frac{2\alpha sW(G)}{n}.
\end{align*}
This proves \eqref{def}. Let $t$ be the number of generalized distance eigenvalues of $K_{a,n-a}$ which are greater than or equal to $\frac{2\alpha W(K_{a,n-a})}{n}$ Then $1\le t\le n$. Therefore, for $n\geq 2a$ and $\alpha\in [\frac{1}{2},1)$, from (\ref{def}), we have
\begin{align}\label{y1}
E^\mathcal{D}_{\alpha}(G)=&2\max_{1\leq j\leq n}\left\lbrace \sum\limits_{i=1}^{j}\partial_i-\frac{2\alpha jW(G) }{n}\right\rbrace\geq 2\sum\limits_{i=1}^{t}\partial_i(G)-\frac{4t\alpha W(G) }{n}\nonumber\\&
\geq   2\sum\limits_{i=1}^{t}\partial_i(K_{a,n-a})-\frac{4t\alpha W(G)}{n},
\end{align} where $t$ is a positive integer defined above.\\
\indent In \cite{pham}, it is shown that the generalized distance spectrum of $ K_{a,n-a} $ is
$$\left\lbrace  (\alpha(a +n)-2)^{[a-1]},(\alpha(2n-a)-2)^{[n-a-1]},\frac{\alpha n+2n-4 \pm \sqrt{\Delta}}{2} \right\rbrace ,$$
where $ \Delta = (2a^2+n^2-2na)(\alpha-2)^2+2(na-a^2)(\alpha^2-2)$ and $ \partial_1=\frac{\alpha n+2n-4 + \sqrt{\Delta}}{2} $ is the spectral radius. Also, we have $\frac{2\alpha W(K_{a,n-a})}{n}=\frac{\alpha(2a^2+2n^2-2n-2na) }{n}.$ When $\frac{1}{2}\leq\alpha< 1$, we will show that $t=1$ or  $2$ or $n-a$ or $n-a+1$. It is easy to see that $\alpha(2n-a)-2\geq \alpha(n+a)-2$, for all $n-a\geq a$. We claim that $\alpha(2n-a)-2$ is the second largest generalized distance eigenvalue of $K_{a,n-a}$ for $n\ge 3.5 a$. We have
\begin{align*}
\alpha(2n-a)-2\ge \frac{(\alpha+2)n-4- \sqrt{\Delta}}{2},
\end{align*} which implies that
\begin{align}\label{xa}
\sqrt{\Delta}\ge 2n-(3n-2a)\alpha.
\end{align} Since right hand side of \eqref{xa} is less or equal to zero  for $\frac{2n}{3n-2a}\leq \alpha\leq 1$, it follows that our claim is true for $\alpha\in \left [\frac{2n}{3n-2a},1\right ]$. So, suppose that $\alpha\in \left [\frac{1}{2},\frac{2n}{3n-2a}\right )$, and for this value of $\alpha$, from \eqref{xa}, we have $4(-2n^2+3an-a^2)\alpha^2+8\alpha(n^2-a^2)-12a(n-a)\geq 0,$ that is,
\begin{align*}
-4(n-a)\Big((2n-a)\alpha^2-2\alpha(n+a)+3a\Big)\geq 0.
\end{align*}
Consider the function $f(\alpha)=(2n-a)\alpha^2-2\alpha(n+a)+3a$, for all $\alpha\in \left [\frac{1}{2},\frac{2n}{3n-2a}\right )$. It is easy to see that this function is increasing for $\alpha\geq \delta$ and decreasing for $\alpha<\delta$, where $\delta=\frac{n+a}{2n-a}$. Clearly,  $\delta\in \left [\frac{1}{2},\frac{2n}{3n-2a}\right )$, for $n>2a$. Since $f(\alpha)$ is increasing for $\alpha\geq \delta$ and decreasing for $\alpha<\delta$, it follows that the sign of $f(\frac{1}{2})$ will decide the nature of $f(\alpha)$ in $\left [\frac{1}{2},\delta\right ]$ and sign of $f(\frac{2n}{3n-2a})$ will decide the nature for $f(\alpha)$ in $[\delta,\frac{2n}{3n-2a}]$. We have $f(\frac{1}{2})=-\frac{1}{4}(2n-7a)\le 0$, for all $n\ge \frac{7a}{2}$. This shows that Inequality  \eqref{xa} holds for all $n\ge 3.5 a$.  If $2a< n<3.5 a$, then for some values of $\alpha\in [\frac{1}{2},\delta)$ the eigenvalue $\alpha(2n-a)-2$ is the second largest, while as for some other values $\alpha\in [\frac{1}{2},\delta)$ the eigenvalue $\frac{(\alpha+2)n-4- \sqrt{\Delta}}{2}$ is the second largest. For $2a< n\le 3a$, it is easy to see that $\delta\ge 0.8$ implying that $0.8\in [\frac{1}{2},\delta)$. We have $f(0.8)=\frac{1}{25}(19a-8n)\le 0$, provided that $n\ge 2.375 a$. This gives that $f(0.8)> 0$ for  $2a< n< 2.375a$ and  $f(0.8)\le 0$ for  $2.375a\le n\le 3a$. Again $0.6, 0.7\in  [\frac{1}{2},\delta)$, for all $2a< n<3.5 a$ and we have $f(0.6)=\frac{1}{25}(36a-12n)$ and $f(0.7)=\frac{1}{100}(111a-42n)$, giving that $f(0.6)> 0$ for  $2a\le n< 3a$ and  $f(0.6)\le 0$ for  $3a\le n< 3.5a$, while as $f(0.7)> 0$ for  $2a\le n< \frac{37a}{14}$ and  $f(0.7)\le 0$ for  $\frac{37a}{14}\le n< 3.5a$.  Further, for $\alpha\in [\delta,\frac{2n}{3n-2a}]$, we have $f(\frac{2n}{3n-2a})=(2n-a)(\frac{2n}{3n-2a})^2-2(n+a)(\frac{2n}{3n-2a})+3a=\frac{1}{(3n-2a)^2}(n-a)(n-2a)^2(3a-4n)<0$, for all $n\ge 2a$.  Thus, it follows that for $\alpha\in [0.5,\delta)$ with $n\ge 3.5 a$ and for  $\alpha\in [\delta,1)$ with $n> 2a$, the eigenvalue  $\alpha(2n-a)-2$ is the second largest. For $\alpha\in [0.5,\delta)$ with $2a< n<3.5 a$, there are two choices for the second largest eigenvalue, depending upon the value of $\alpha$. In fact, if $\alpha\in [0.5,0.6]$, then for $2a<n\le 3a$, the eigenvalue $\frac{(\alpha+2)n-4- \sqrt{\Delta}}{2}$  is the second largest, while as the eigenvalue $\alpha(2n-a)-2$ is the third largest and for $3a<n< 3.5a$, the eigenvalue $\alpha(2n-a)-2$  is the second largest, while as the eigenvalue $\frac{(\alpha+2)n-4- \sqrt{\Delta}}{2}$  is the third largest. For $\alpha\in [0.6, \delta)$, we have both the possibilities. For the eigenvalue $\frac{(\alpha+2)n-4- \sqrt{\Delta}}{2}$ and $ \alpha(n+a)-2$, proceeding similarly it can be seen that $\frac{(\alpha+2)n-4- \sqrt{\Delta}}{2}\geq \alpha(n+a)-2$ for all $\alpha\in [\frac{1}{2},1]$ and $a\ne 1$. Thus it follows that for $\alpha\in [0.5,\delta)$ with $n\ge 3.5 a$ and for  $\alpha\in [\delta,1)$ with $n>2a$, $a\ne 1$, we have
\begin{align*}
\frac{\alpha n+2n-4 +\sqrt{\Delta}}{2}\geq (\alpha(2n-a)-2)\geq\frac{\alpha n+2n-4 -\sqrt{\Delta}}{2}\geq (\alpha(a +n)-2).
\end{align*} For  $\alpha\in [0.5,\delta)$ with $2a\le n<3.5 a$, we have
\begin{align*}
\frac{\alpha n+2n-4 +\sqrt{\Delta}}{2}\geq \phi_1\geq\phi_2\geq (\alpha(a +n)-2),
\end{align*} where $\phi_1=\max\{(\alpha(2n-a)-2),\frac{\alpha n+2n-4 -\sqrt{\Delta}}{2}\}$ and  $\phi_2=\min\{(\alpha(2n-a)-2),\frac{\alpha n+2n-4 -\sqrt{\Delta}}{2}\}$.
For the eigenvalue $\alpha(n+a)-2$, it is easy to see that $\alpha(n+a)-2< \frac{2\alpha W(K_{a,n-a})}{n}$ holds for all $n> 2a$ and $\alpha\ge 0.5 $. For the eigenvalues $(\alpha(2n-a)-2)$, with $n>2a$, we have
\begin{align*}
\frac{2\alpha W(K_{a,n-a})}{n}=\frac{\alpha(2a^2+2n^2-2n-2na) }{n}\leq (\alpha(2n-a)-2),
\end{align*}
which gives
\begin{align}\label{x1}
\alpha(2a^2-2n-an)+2n\leq 0.
\end{align}
This in turn implies that  (\ref{x1}) is true for  $\alpha\geq \gamma$, where $\gamma=\frac{2n}{2n+an-2a^2}$. For $a\geq 2$, it can be seen that $\gamma<\frac{1}{2}$, for all $n\ge \frac{2a^2}{a-2}$ implying that   (\ref{x1}) is true for all $\alpha\geq \frac{1}{2}$. Since $3.5 a\ge \frac{2a^2}{a-2}$, for all $a\ge 5$, it follows that Inequality (\ref{x1}) is true for all $n\ge 3.5 a$ with $a\ge 5$ and is true for $2\le a\le 4$, provided that $n\ge \frac{2a^2}{a-2}$. So, suppose $3.5a\leq n<\frac{2a^2}{a-2}$. For this value of $n$, it is clear that $\gamma\in [\frac{1}{2},1]$ and therefore Inequality (\ref{x1}) holds for $ \alpha\geq \gamma$ and does not hold for $\alpha<\gamma$. This shows that for $a\geq 5$ and  $n\ge 3.5 a$; for $2\le a\le 4$ and $n\ge \frac{2a^2}{a-2}$; for  $3.5a\leq n<\frac{2a^2}{a-2}$ and $\alpha\geq \gamma$, we have $t\geq n-a$; while for $3.5a\leq n<\frac{2a^2}{a-2}$ and $\alpha<\gamma$, we have $t=1$ or $2$.\\
\indent For  $n\ge 3.5a$, we have
\begin{align*}
\frac{\alpha n+2n-4 - \sqrt{\Delta}}{2}< \frac{\alpha(2a^2+2n^2-2n-2na) }{n},
\end{align*}
implying
\begin{align}\label{x2}
(2-3\alpha)n^2+4(\alpha-1)n-4a\alpha(a-n)< n\sqrt{\Delta}.
\end{align}
Consider the function $h(\alpha)=(2a^2+n^2-2na)(\alpha-2)^2+2(na-a^2)(\alpha^2-2)$, for $\alpha\in [\frac{1}{2},1]$. It is easy to see that the function $h(\alpha)$ is a decreasing function of $\alpha$. Therefore the minimum value of $h(\alpha)$ is attained at $\alpha=1$, that is, $h(\alpha)\geq h(1)=(n-2a)^2$. With this, from (\ref{x2}), we have
$(2-3\alpha)n^2+4(\alpha-1)n-4a\alpha(a-n)< n(n-2a)$, implies that $n^2+2an-4n<\alpha(3n^2-4n(a+1)+4a^2)$, 
which in turn implies that $\alpha >\frac{n^2+2an-4n}{3n^2-4n(a+1)+4a^2}=\gamma_1$. This implies that inequality \ref{x2} holds for $\alpha\in (\gamma_1,1)$. It is easy to see that for $n\ge 3.5 a$, we have $\gamma_1< \frac{3}{4}$. Since $h(\alpha)$ is a decreasing function of $\alpha$, therefore $h(\alpha)\ge h(\gamma_1)>h(\frac{3}{4})$, for all $\alpha\in [0.5,\gamma_1]$. We have $ h(\frac{3}{4})=\frac{1}{16}\Big(25n^2-96an+96a^2\Big)> \frac{1}{16}\Big(5n-10a\Big)^2$. With this it follows from inequality \ref{x2} that for $\alpha\in [0.5,\gamma_1]$, we have $(2-3\alpha)n^2+4(\alpha-1)n-4a\alpha(a-n)< n\Big(\frac{5n-10a}{4}\Big)$ implying that $\alpha> \frac{3n^2+10an-16n}{4(3n^2-4n(a+1)+4a^2)}=\gamma_2$. This shows that Inequality \eqref{x2} holds for $\alpha\in (\gamma_2,1)$. It is easy to see that for $n\ge 3.5 a$, we have $\gamma_2< \frac{2}{3}$. Again using the fact that $h(\alpha)$ is a decreasing function of $\alpha$, it follows that $h(\alpha)\ge h(\gamma_2)>h(\frac{2}{3})$, for all $\alpha\in [0.5,\gamma_2]$.  We have $ h(\frac{2}{3})=\frac{1}{9}\Big(16n^2-60an+60a^2\Big)$. From Inequality \eqref{x2}, it follows that for $\alpha\in [0.5,\gamma_2]$, we have $(2-3\alpha)n^2+4(\alpha-1)n-4a\alpha(a-n)< \frac{n}{3}\sqrt{16n^2-60an+60a^2}$ implying that $\alpha> \frac{6n^2-12n-n\sqrt{16n^2-60an+60a^2}}{3(3n^2-4n(a+1)+4a^2)}=\gamma_3$. This shows  that Inequality \eqref{x2} holds for $\alpha\in (\gamma_3,1)$. It is easy to see that $\gamma_3\le 0.5$, for $n\ge 4.877 a$, implying that the Inequality \eqref{x2} holds for $\alpha$, provided that $n\ge 4.877 a$. For $2a< n< 4.877 a$, the Inequality \eqref{x2} does not hold for some values of $\alpha$ as can been seen as follows. Since $h(\alpha)$ is a decreasing function of $\alpha$, it follows that $h(\alpha)\le h(\frac{1}{2})=\frac{1}{4}\Big(9n^2-32an+32a^2\Big)$. Therefore, from Inequality \eqref{x2}, it follows that $(2-3\alpha)n^2+4(\alpha-1)n-4a\alpha(a-n)< \frac{n}{2}\sqrt{9n^2-32an+32a^2}$ implying that $\alpha> \frac{2n^2-4n-0.5 n\sqrt{9n^2-32an+32a^2}}{3n^2-4n(a+1)+4a^2}=\gamma_4$. If $\gamma_4\le 0.5$, then 
\begin{align}\label{ya}
8n^4-(40a-8)n^3+(24a^2+32a-16)n^2+32na^2(a-1)-16a^4\ge 0.
\end{align} 
Consider the function $p(n)=nq(n)-16a^4$, where $q(n)=8n^3-(40a-8)n^2+(24a^2+32a-16)n+32a^2(a-1)$. It is easy to see that the function $q(n)$ is  decreasing in $(2a,\psi)$ and increasing in $[\psi, 4.877a)$. Since $\psi<3a$, so we have $q(2a)=-16a(a^2-4a+2)<0$, for $a\ge 2$, and $q(4a-3)=-8(a-3)(5a-5)\le 0$, for $a\ge 2$. For $2a< n< 4a-3$, we have $q(n)<0$ implying that $p(n)<0$ in this range of $n$. From this, it follows that for $2a< n< 4a-3$ with $\alpha \in [0.5,\gamma_4]$, the Inequality  \eqref{x2} does not hold.\\
\indent Thus, it follows that if $n\ge 4.877 a$, then $t=n-a$; if $\frac{2a^2}{a-2}\le n< 4.877 a$, then $t=n-a$ or $t=n-a+1$; if $2a< n< \frac{2a^2}{a-2}$ and $\alpha\ge \gamma$, then $t=n-a$ or $t=n-a+1$ and if $2a< n< \frac{2a^2}{a-2}$ and $\alpha< \gamma$, then $t=1$ or $t=2$. Using this information, part \textbf{(1)} and \textbf{(2)} now follows by direct calculation from (\ref{y1}).\\
\indent If $n=2a$, then the generalized distance spectrum of $ K_{a,n-a} $ is $\{3a-2, 2a\alpha+a-2,3a\alpha-2^{[2a-2]}\}$, with $3a-2\ge 2a\alpha+a-2\ge 3a\alpha-2$. Also, $\frac{2\alpha W(K_{a,a})}{n}=\alpha(3a-2)$. It is easy to see that $t=2$, in this case. With this the result now follows by direct calculation from (\ref{y1}). This completes the proof. 
\qed \\

The following lemma can be found in \cite{CTZ}.

\begin{lemma}\label{lam12}
 Let $ T $ be a tree of order $ n\geq 4.$ If $ 0\leq\alpha<1,$ then $$\partial(T)\geq \partial(S_{n})= \frac{\alpha n+2n-4+\sqrt{k}}{2},$$
where $ k=(n^2-2n+2)(\alpha-2)^2+2(n-1)(\alpha^2-2), $ with equality if and only if $ T\cong S_{n}. $
\end{lemma}

The following result gives the generalized distance energy of the star graph of order $n$.

\begin{theorem}\label{star}
For $\alpha\in (0,1)$, the generalized distance energy of the star graph $ S_{n}$ of order $n, n\geq 4$, is
\begin{equation*}
E^\mathcal{D}_{\alpha}(S_{n}) = \begin{cases}
\sqrt{k}+2n+\alpha(8-3n)-4-\frac{4\alpha}{n}, &0<\alpha<\frac{2n}{3n-2}\\
\sqrt{k}-n\alpha(4n-19)-22\alpha -2n+4+\frac{8\alpha}{n}, &\frac{2n}{3n-2}\leq \alpha<1
\end{cases},
\end{equation*}
where $ k=(n^2-2n+2)(\alpha-2)^2+2(n-1)(\alpha^2-2).$
\end{theorem}
\noindent\textbf{Proof.}
 Let $t$ be the greatest positive integer such that $\partial_t(S_{n})\geq \frac{2\alpha W(S_{n})}{n}$. The generalized distance spectrum of $ S_{n} $ is $$\left\lbrace  (\alpha(2n-1)-2)^{[n-2]}, \frac{\alpha n+2n-4\pm \sqrt{k}}{2}\right\rbrace,$$ where $k=(n^2-2n+2)(\alpha-2)^2+2(n-1)(\alpha^2-2)$. For $n\geq 4$, from Lemma \ref{lam12}, it follows that
 $ \frac{\alpha n+2n-4+ \sqrt{k}}{2}$ is the generalized distance spectral radius of $ S_{n}$, which is always greater or equal to $\frac{2\alpha W(S_{n})}{n}=\alpha(2n-4+\frac{2}{n})$. It is easy to verify that $ \frac{\alpha n+2n-4- \sqrt{k}}{2} $ is the least eigenvalue of $ S_{n}.$ For the eigenvalue $ \alpha(2n-1)-2$, we have two different cases.\\
{\bf Case (i)}. If $ 0<\alpha<\frac{2n}{3n-2}$, then it is easy to see that $ \alpha(2n-1)-2<\alpha(2n-4+\frac{2}{n})$ and so $t=1$. Therefore, using (\ref{def}), we have
 \begin{align*}
 E^{\mathcal{D}_\alpha}(S_{n})&=2\left(\sum\limits_{i=1}^{t}\partial_i-\frac{2t\alpha W(G)}{n} \right)\\&
 = \sqrt{k}+2n+\alpha(8-3n)-4-\frac{4\alpha}{n}.
 \end{align*}
{\bf Case (ii)}. If $ \frac{2n}{3n-2}\leq \alpha <1$, then it is easy to see that $ \alpha(2n-1)-2\geq\alpha(2n-4+\frac{2}{n})$ and so $t=n-1$. Therefore, using (\ref{def}), we have
  \begin{align*}
  E^{D_\alpha}(S_{n})&=2\left(\sum\limits_{i=1}^{t}\partial_i-\frac{2t\alpha W(G)}{n} \right)\\&
  =2\sum\limits_{i=1}^{n-1}\partial_i- 4\alpha(n-1)(2n-4+\frac{2}{n})\\&
  =\sqrt{k}+2n\alpha(n-2)+2\alpha-\alpha n+4-4\alpha(n-1)(2n-4+\frac{2}{n})\\&
  = \sqrt{k}-n\alpha(4n-19)-22\alpha -2n+4+\frac{8\alpha}{n}.
  \end{align*}
\qed

Now we obtain a lower bound for the generalized distance energy of a tree in terms of the order $n$, the Wiener index  and the parameter $\alpha$.

\begin{theorem} \label{thmst}
Let $T$ be a tree of order $ n$ and let $\alpha\in (0,1)$. Then
\begin{equation}\label{sl}
E^{D_\alpha}(T)\geq\alpha n+2n-4+\sqrt{k}- \frac{4\alpha W(G)}{n},
\end{equation}  where $ k=(n^2-2n+2)(\alpha-2)^2+2(n-1)(\alpha^2-2)$. Equality occurs for $ 0<\alpha<\frac{2n}{3n-2}$, if and only if $T\cong S_{n}$.
\end{theorem}
\noindent\textbf{Proof.}
Let $T$ be a tree of order $n$ and let $t$ be the largest positive integer such that $\partial_t(T)\geq \dfrac{2\alpha W(T)}{n}$.  Using  \eqref{def} and Lemma \ref{lam12}, we have
  \begin{align*}
E^{\mathcal{D}_\alpha}(T)&=2\left(\sum\limits_{i=1}^{t}\partial_i(T)-\frac{2t\alpha W(T)}{n} \right)=2\max_{1\leq j\leq n}\left\lbrace \sum\limits_{i=1}^{j}\partial_i(T)-\frac{2\alpha jW(T) }{n}\right\rbrace \\&
 \geq 2(\partial_1(T)-\dfrac{2\alpha W(T)}{n})
  \geq 2\partial_1(S_{n})- \frac{4\alpha W(T)}{n}\\
& \geq \alpha n+2n-4+\sqrt{k}- \frac{4\alpha W(T)}{n},
 \end{align*} where $ k=(n^2-2n+2)(\alpha-2)^2+2(n-1)(\alpha^2-2)$. Assume that equality occurs  in (\ref{sl}), then all the inequalities above occur as equalities. Since equality occurs in Lemma \ref{lam12} for $T\cong S_{n}$, from \eqref{def}, it follows  that  equality occurs in (\ref{sl}) if and only if $t=1$ and $T\cong S_{n}$. Now using Lemma \ref{star} the result follows.  \qed

For $ 0<\alpha<\frac{2n}{3n-2}$, Theorem \ref{thmst} shows that among all trees of order $n$ the star graph $S_{n}$ has the minimum generalized distance energy.

\newpage
\textbf{Papers Accepted/ Published:}
\begin{enumerate}
\item	S. Pirzada, Hilal A. Ganie, \textbf{Bilal A. Rather} and R. U. Shaban, On the generalized distance energy of graphs, \textit{Linear Algebra and its Applications} \textbf{603} (2020) 1--19.\\
\textbf{SCI, Scopus}, ISSN: 0024-3795.

\item	Hilal A. Ganie, S. Pirzada, \textbf{Bilal A. Rather} and V. Trevisan, Further developments on Brouwer's conjecture for the sum of Laplacian eigenvalues of graphs, \textit{Linear Algebra  and its Applications} \textbf{588} (2020) 1--18.\\
\textbf{SCI, Scopus}, ISSN: 0024-3795.

\item	S. Pirzada, \textbf{Bilal A. Rather}, M. Aijaz and T. A. Chishti, Distance signless Laplacian spectrum of graphs and spectrum of zero divisor graphs of  $ \mathbb{Z}_{n} $,  \textit{Linear and MultiLinear Algebra} (2020) DOI:10.1080/03081087.2020.1838425.\\
\textbf{SCI, Scopus}, ISSN: 03081087

\item	S. Pirzada, \textbf{Bilal A. Rather}, Hilal A. Ganie and R. U. Shaban, On Brouwers conjecture and spectrally threshold graphs, 
accepted for publication in\textit{ Journal of the Ramanujan Mathematical Society} (January 2021).\\
\textbf{SCI, Scopus}, ISSN: 23203110 (online), ISSN: 09701249 (print).

\item	S. Pirzada, \textbf{Bilal A. Rather}, Hilal A. Ganie, and R. U. Shaban, On generalized distance spectral radius of a bipartite graph, \textit{Matemati\u{c}ki Vesnik} \textbf{72(4)} (2020) 327--336.\\
Web of Science \textbf{(ESCI), Scopus}, ISSN: 2406-0682 (online), ISSN: 0025-5165 (print).

\item	S. Pirzada, \textbf{Bilal A. Rather}, R. U. Shaban and T. A. Chishti, On the sum of powers of $ A_{\alpha} $ eigenvalues of graphs and $ A_{\alpha} $-energy like invariant, \textit{Boletim da Sociedade Paranaese de Matem\'atica} (2020) DOI:10.5269/bspm.52469.\\
Web of Science \textbf{(ESCI), Scopus}, ISSN: 2175-1188 (online), ISSN: 0037-8712 (print).

\item	S. Pirzada, \textbf{Bilal A. Rather} and R. U. Shaban, On graphs with minimal distance signless Laplacian energy, accepted for publication in \textit{ Acta Universitatis Sapientiae Mathematica} ( August 2020).\\
Web of Science \textbf{(ESCI), Scopus}, ISSN: 2066-7752 (online), ISSN: 1844-6094 (print).

\item	T. A. Naikoo, U. Samee, S. Pirzada and \textbf{Bilal A. Rather}, On degree sets in k-partite graphs, \textit{Acta Universitatis Sapientae Informatica} \textbf{12}(2) (2020) 251--259.\\
Web of Science \textbf{(ESCI)}, ISSN: 2066-7760 (online), ISSN: 1844-6086 (print).

\item	S. Pirzada, \textbf{Bilal A. Rather} and R. U. Shaban and Merajuddin, On signless Laplacian spectrum of zero divisor graphs of the ring $ \mathbb{Z}_{n} $, accepted for publication in  \textit{ The Korean Journal of Mathematics} ( Jan 2021).\\
Web of Science \textbf{(ESCI)}, ISSN: 2288-1433 (online), ISSN: 1976-8605 (print).

\item	T. A. Naikoo, \textbf{Bilal A. Rather} and S. Pirzada, On Zagreb index of tournaments, accepted for publication in  \textit{ Kragujevac Journal of Mathematics } ( Jan 2021).\\
Web of Science \textbf{(ESCI), Scopus}, ISSN: 2406-3045 (online), ISSN: 1450-9628 (print).
\end{enumerate}
\newpage
\textbf{Communicated publications}

\begin{enumerate}
\item	R. U. Shaban, Bilal A. Rather, S. Pirzada, On distance signless Laplacian spectral radius of power graphs of cyclic and dihedral groups, Submitted (Nov 2020).
\item	Bilal A. Rather, S. Pirzada, M Aouchiche, On eigenvalues and energy of geometric-arithmetic index of graphs, Submitted (Nov 2020).
\item	S. Pirzada, Bilal A. Rather, T. A. Naikoo, On general Zagreb index of tournaments, Submitted (Nov 2020).

\item	Bilal A. Rather, S. Pirzada, T. A. Naikoo, On Laplacian eigenvalues of zero divisor graphs of finite commutative ring, Submitted  (Oct 2020).
\item	Bilal A. Rather, Hilal A. Ganie and S. Pirzada, On the $ A_{\alpha} $-spectrum of joined union and its applications to power graphs of certain finite groups, Submitted  (Aug 2020).
\item	Hilal A. Ganie, Bilal A. Rather and S. Pirzada, Laplacian energy of trees of diameter four and beyond, Submitted (July 2020).
\item	Hilal A. Ganie, R. U. Shaban, Bilal A. Rather and S. Pirzada, On distance Laplacian energy, vertex connectivity and  independence number of graphs, Submitted (April 2020).
\item	S. Pirzada, Bilal A. Rather, Hilal A. Ganie and R. U. Shaban, On $ \alpha $-adjacency energy of graphs, Submitted (May 2020).
\item	S. Pirzada, Bilal A. Rather and T. A. Chishti , On normalized Laplacian spectrum of zero divisor graphs of commutative ring  $ \mathbb{Z}_{n} $, Submitted  (Febuary 2020).

\item	Bilal A. Rather, S. Pirzada and G. Zhou, On distance Laplacian spectra of power graphs of certain finite groups, Submitted  (March 2020).
\end{enumerate}

\end{document}